# University of Aix-Marseille

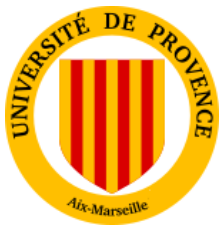

**Department of Mathematics**

**Memoire**

# CALABI CONJECTURE

**Research Master**
**In**

**Fundamental of Mathematics (M2)**

**Submitted by Hassan Jolany**

**Supervisor: Prof. Nader Yeganefar**

Marseille, June 2012

## Acknowledgment

This Memoire marks the end of a period and opens a new chapter of my academic life. I would like to take this opportunity to express my appreciation to all persons who fruitfully helped me along my journey of academic pursuit. First and foremost, I would like to express my heartfelt gratitude to my advisor Professor Nader Yeganefar for all his encouragement. I also wish to thank of my professors in Centre de Mathématiques et Informatique (CMI) and Luminy Institute of Mathematics (at Aix-Marseille University), in particular Professors Alexander Borichev (for his two lectures about Holomorphic space I,II), Professor Ctirad Klimčík (for his lecture about Poisson Manifolds and Poisson Lie structures ) and Professor Patrick Delorme (for his lecture about Lie Algebra and Lie Groups ). I also wish to thanks for the scholarship of council of city of Marseille and Association SARA-GHU à Marseille, for its hospitality and support during the writing of this Memoire.

## Abstract

The main propose of this memoire is understanding the alternative proof of Calabi conjecture which posted by H.-D.-Cao. Now here we explain this conjecture and its history. Let $(M, \omega)$ is a Kahler manifold and $\rho$ be its associated Ricci form. The form $\rho / 2\pi$ represents the first Chern class $C_1(M)$ of $M$. E.Calabi asks if given a closed $(1,1)$ $-$form $\tilde{\rho}$ which represents the first Chern class of $M$, one can find a Kahler metric $\tilde{\omega}$ on $M$ whose Ricci tensor is $\tilde{\rho}$. Also he shows that this metric is unique on each Kahler class (i.e. a cohomology class of type $(1,1)$ that contains a form that is positive definite). The existence of $\tilde{\omega}$ on each Kahler class is known as Calabi's conjecture and it has been solved by S.T.Yau in 1976. Calabi's Conjecture has some immediate consequence of the existence and uniqueness of Kahler Einstein manifold. Let $\rho = \lambda \omega$ where $\lambda \in R$. According to the sign of $\lambda$, the first chern class $C_1(M)$, must either vanish or have a representative which is negative or positive. When $C_1(M)$ vanishes, $M$ admits an unique Ricci-flat metric on each Kahler class. When $C_1(M) < 0$ , S.-T.-Yau prove that there exists a unique (up to homotheties) Kahler-Einstein metric depending only on the complex structure $M$. When $C_1(M) > 0$ the existence of Kahler-Einstein metric is not guaranteed. The problem has been largely studied by G.Tian. Y.Matsushima proves that when $C_1(M) > 0$ . If two Einstein forms $\omega_1$, $\omega_2$ are cohomologous then they are isometric, namely there exists $F \in Aut(M)$ such that $F^* \omega_2 = \omega_1$. H.-D.-Cao proved a conjecture about deformation of Kahler Einstein metrics and by his theorem easily Calabi's Conjecture could be verified. Now we summarize our result as follows. This Memoire consists of two main results. In the first one we describe Ricci flow theory and we give an educative way for proving Elliptization Conjecture and then we prove Poincare conjecture which is well-known the second proof of Perelman for Poincare conjecture

(of course our proof was trivial in his second paper. But we verify that why of his theorems we can prove Poincare conjecture). In the second one which is the main propose of our memoire we exhibit a complete proof of Calabi-Yau conjecture by using Cao's method.

# Chapter 1
# Chapter 1

# • Ricci flow

In this chapter we give some elementary properties of Ricci flow theory and after by using Perelman's results we will proof the Poincare conjecture. The study of geometric flows, especially Ricci flow, is one of the most central topics in modern geometric analysis. By applying the theory of geometric flow, several breakthroughs were made in the last decade in resolving several log-standing conjectures. For instance, Poincare conjecture, Thurston's geometrization conjecture, or a conjecture about Differentiable Sphere Theorem and etc. All the above-mentioned achievements involve seeking canonical metrics by geometric deformations.

The Ricci flow, which evolves a Riemannian metric by its Ricci curvature, is a natural analogue of the heat equation for metrics. Firstly, we introduce flows $t \to (M(t), g(t))$ on Riemannian manifolds $(M, g)$, which are important for describing smooth deformations of such manifolds over time, and derive the basic first variation formulae for how various structures on such manifolds change by such flows. In this section we get a one-parameter family of such manifolds $t \to (M(t), g(t))$, parameterized by a "time" parameter $t$. In the usual manner the time derivatives $\dot{g}(t) = \frac{d}{dt} g(t)$ is

$$\frac{d}{dt} g(t) = \lim_{dt \to 0} \frac{g(t + dt) - g(t)}{dt}$$

and also the analogue of the time derivative $\frac{d}{dt} g(t)$ is then the Lie derivative $\mathcal{L}_{\partial_t} g$ .

**Definition (Ricci flow)**: A one-parameter family of metrics $g(t)$ on a smooth manifold $M$ for all time t in an interval I is said to obey Ricci flow if we have

$$\frac{d}{dt} g(t) = -2Ric(t)$$

Note that this equation makes tensorial sense since $g$ and $Ric$ are both symmetric rank *2* tensor. The factor of 2 here is just a notational convenience and is not terribly important, but the minus sign$-$ is crucial. Note that, due to the minus sign on the right-hand side of the equation, a solution to the equation shrinks in directions of positive Ricci curvature and it expands in direction of negative Ricci curvature (you can check the first example of this chapter). For more details see [1-9].

As geometrically Ricci flow equation is good guy. Because it is invariant under the full diffeomorhism group, the reason is: If $\varphi: M \to \tilde{M}$ is a time-independent diffeomorphism such that $g(t) = \varphi^* \tilde{g}(t)$ and $\tilde{g}$ is a solution of the Ricci flow, we get

$$\frac{\partial}{\partial t} g = \varphi^* \left( \frac{\partial}{\partial t} \tilde{g} \right) = \varphi^* (-2Ric(\tilde{g})) = -2Ric(\varphi^* \tilde{g}) = -2Ric(g)$$

Where the second last equality come from of the fact that if $\varphi: (M, g) \to (\tilde{M}, \tilde{g})$ is local isometry then $\varphi^* \tilde{R} = R$ .

Let us give a quick introduction of what the Ricci flow equation means . In harmonic coordinates $(x^1, x^2, \dots, x^n)$ about $p$, that is to say coordinates where $\Delta x^i = 0$ for all $i$ ,we have

$$Ric_{ij} = Ric \left( \frac{\partial}{\partial x^i}, \frac{\partial}{\partial x^j} \right) = -\frac{1}{2} \Delta g_{ij} + Q_{ij}(g^{-1}, \partial_t g)$$

Where $Q$ is a quadratic form in $g^{-1}$ and $\partial g$ .See lemma (3.32) on page (92) of [40] for more details

So in these coordinates, the Ricci flow equation is actually a heat equation for the Riemannian metric

$$\frac{\partial g}{\partial t} = \Delta g + 2Q(g^{-1}, \partial_t g)$$

In local coordinates

$$-2R_{ij} = \sum_{k,l=1}^{n} g^{kl} \left\{ \frac{\partial^2 g_{kl}}{\partial x^i \partial x^j} - \frac{\partial^2 g_{jl}}{\partial x^i \partial x^k} - \frac{\partial^2 g_{il}}{\partial x^j \partial x^k} + \frac{\partial^2 g_{ij}}{\partial x^k \partial x^l} \right\}$$

$$+ \left( \text{quadratic terms in } \frac{\partial}{\partial x^k} g_{ij} \right)$$

- **Special solutions of flows**

The Ricci flow $\frac{\partial g}{\partial t} = -Ric$ introduced by Hamilton is a degenerate parabolic evolution system on metrics (In reality Ricci flow is a weakly parabolic system where degeneracy comes from the gauge invariance of the equation under diffeomorphism but understanding this issue is far of the main propose of my memoire! But the reader can check it in [12])

***Theorem*** (**Hamilton** [1]) Let $(M, g_{ij}(x))$ be a compact Riemannian manifold .Then there exists a constant $T > 0$ Such that the Ricci flow $\frac{\partial g}{\partial t} = -2Ric$, $with$ $g_{ij}(x,0) = g_{ij}(x)$ ,admits a unique smooth solution $g_{ij}(x,t)$ for all $x \in M$ and $t \in [0, T)$

- *Example 1(Einstein metric)*

Before of starting our example, we need to a corollary which we state here

Corollary: If $\tilde{g} = cg$, are two Riemannian metrics on an n-manifolds $M^n$, related by a scaling factor $C$ , then various geometric quantities scale as follows:

- $\tilde{g}^{ij} = c^{-1} g^{ij}$

- $\tilde{\Gamma}_{ij}^{\ k} = \Gamma_{ij}^{\ k}$

- $\tilde{R}_{ij} = R_{ij}$

- $\tilde{R} = c^{-1} R$

- $d\tilde{\mu} = c^{\frac{n}{2}} d\mu$

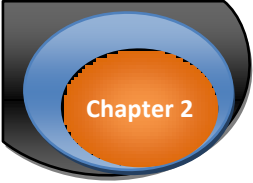

Recall that a Riemannian metric $g_{ij}$ is Einstein if $R_{ij} = \lambda g_{ij}$ for some for some constant $\lambda$ .Now since the initial metric is Einstein with positive Einstein constant $\lambda > 0$ so that

$$R_{ij}(x, 0) = \lambda g_{ij}(x, 0)$$

*Let* $g_{ij}(x, t) = \rho^2(t) g_{ij}(x, 0)$.

Now from the previous corollary for the Ricci tensor, one sees that

$$R_{ij}(x, t) = R_{ij}(x, 0) = \lambda g_{ij}(x, 0)$$

Thus the equation $\frac{\partial g_{ij}}{\partial t} = -2 Ric_{ij}$ corresponds to

$$\frac{\partial \left( \rho^2(t) g_{ij}(x, 0) \right)}{\partial t} = -2\lambda g_{ij}(x, 0)$$

this gives the *ODE* $\frac{d\rho}{dt} = -\frac{\lambda}{\rho}$ .

Whose solution is given by $\rho^2(t) = 1 - 2\lambda t$ .Thus the evolving metric $g_{ij}(x, t)$ shrinks homothetically to a point as $\rightarrow T = \frac{1}{2\lambda}$ .

By contrast, if the initial metric is an Einstein metric of negative scalar curvature , the metric will expand homothetically for all times .Indeed if $R_{ij}(x, 0) = -\lambda g_{ij}(x, 0)$ with $\lambda > 0$ and $g_{ij}(x, t) = \rho^2(t) g_{ij}(x, 0)$ ,the $\rho(t)$ satisfies the *ODE* $\frac{d\rho}{dt} = \frac{\lambda}{\rho}$ with solution $\rho^2(t) = 1 + 2\lambda t$ , hence the evolving metric $g_{ij}(x, t) = \rho^2(t) g_{ij}(x, 0)$ exists and expands homothetically for all times, and the curvature fall back to zero like $\frac{-1}{t}$ .

**Definition** *(see [18])*.Let $X \in T_p M$ be a unit vector. Suppose that X is contained in some orthonormal basis for $T_p M$. $Rc(X, X)$ is then the sum of the sectional curvatures of planes spanned by $X$ and other elements of the basis.

- **Example 2**

On an n-dimensional sphere of radius $r$ *(where $n > 1$)*, the metric is given by $g = r^2 \overline{g}$ where $\overline{g}$ is the metric on the unit sphere.The sectional curvatures are all $\frac{1}{r^2}$ (because $S^n(r)$ has constant sectional curvature $\frac{1}{r^2}$) . Thus for any unit vector $v$, the result of previous definition tells us that $Rc(v,v) = \frac{(n-1)}{r^2}$. Therefore $Rc = \frac{(n-1)}{r^2} g = (n-1)\overline{g}$

So the Ricci flow equation becomes an ODE

$$\frac{\partial}{\partial t} g = -2Rc$$

$$\Rightarrow \frac{\partial}{\partial t}(r^2\overline{g}) = -2(n-1)\overline{g}$$

$$\Rightarrow \frac{\partial}{\partial t}(r^2) = -2(n-1)$$

We have the solution $(t) = \sqrt{R_0^2 - 2(n-1)t}$ , where $R_0$ is the initial radius of the sphere . The manifold shrinks to a point as $\rightarrow \frac{R_0^2}{2(n-1)}$ .

Similarly, for hyperbolic $n$-space $H^n$ (where n > 1),The Ricci flow reduces to the ODE

$$\frac{\partial}{\partial t}(r^2) = 2(n-1)$$

Which has the solution $(t) = \sqrt{R_0^2 + 2(n-1)t}$ . So the solution expands out to infinity. In reality we used of this fact that $H^n(r)$ has constant sectional curvature $-\frac{1}{r^2}$.

- **Example3 (Ricci Soliton)**

$A$ Ricci soliton is a Ricci flow $(M, g(t))$ , $0 \leq t < T \leq \infty$ ,with the property that for each $t \in [0, T)$ there is a diffeomorphism $\varphi_t : M \rightarrow M$ and a constant $\sigma(t)$ such that $\sigma(t)\varphi_t^* g(0) = g(t)$ . That is to say , in a Ricci soliton all the Riemannian manifold $(M, g(t))$ are isomorphic up to a scale factor that is allowed to vary with t .The soliton is said to be shrinking if $\frac{\partial \varphi}{\partial t}(t) < 0$ for all $t$ (note that $\varphi_0 = id$ and $\sigma(0) = 1$ )

Taking the derivative of the equation $g(t) = \sigma(t)\varphi_t^* g(0)$ and evaluating $at\ t = 0$ yields

$$\frac{\partial}{\partial t} g(t) = \frac{\partial\, \sigma(t)}{\partial t} \varphi_t^* g(0) + \sigma(t)\frac{\partial}{\partial t} \varphi_t^* g(0)$$

$$-2Rc\big(g(0)\big) = \sigma'(0)g(0) + \mathcal{L}_v g(0)$$

where $= \frac{d\varphi_t}{dt}$. *Let us set* $\sigma'(0) = 2\lambda$. Now because on a Riemannian manifold $(M, g)$, we have

$$(\mathcal{L}_X g)_{ij} = \nabla_i X_j + \nabla_j X_i$$

we get

$$-2R_{ij} = 2\lambda g_{ij} + \nabla_i X_j + \nabla_j X_i$$

As a special case we can consider the case that $v$ is the gradient vector field of some scalar function $f$ on $M^n$, *i.e.* $v_i = \nabla_i f$, the equation then becomes $R_{ij} + \lambda g_{ij} + \nabla_i \nabla_j f = 0$ such solutions are known as gradient Ricci solutions .

- **Short time existence on heat equation of the Ricci flow type**

Now we discuss on short time existence for solutions of heat equation of the Ricci flow type

Consider the heat equation

$$\begin{cases} \frac{\partial}{\partial t} u(x, t) = \Delta u(x, t) \quad , \quad x \in \mathbb{R}^n , t > 0 \\ u(x, 0) = u_0(x) \quad x \in \mathbb{R}^n \end{cases}$$

It is well known that

- If $|u_0(x)| \leq c_1 e^{c_2 |x|^2}$ on $\mathbb{R}^n$, then the heat equation admits a unique solution $u(x, t)$ on some short-time interval $[0, T]$ with $|u(x, t)| \leq \tilde{c}_1 e^{\tilde{c}_2 |x|^2}$ on $\mathbb{R}^n \times [0, T]$
- Even if $u_0(x) \equiv 0$, the heat equation has a nontrivial solution with big oscillation

Now we the Ricci flow on non-compact manifolds

$$\frac{\partial}{\partial t} g_{ij} = -2R_{ij}$$

where in Kähler manifold, the right hand side (after we will proof it) is given by

$$-2R_{i\bar{j}} = 2\frac{\partial^2}{\partial z^i \partial \bar{z}^j} \log \det(g_{i\bar{j}})$$

Thus "curvature bounded condition in Ricci flow is analogous the "the growth condition" $|u(x,t)| \leq c_1 e^{c_2 |x|^2}$ in the heat equation.

In 1989, Shi generalized the theorem of Hamilton about short-time existence and uniqueness theorem for the Ricci flow to complete non-compact manifolds with bounded curvature.

**Theorem** (Shi [17]) Let $\left(M, g_{ij}(x)\right)$ be a complete non compact Riemannian manifold of dimension n with bounded curvature .Then there exists a constant T > 0 such that the initial value problem

$$\begin{cases} \dfrac{\partial}{\partial t} g_{ij}(x,t) = -2R_{ij}(x,t) \ on \ M \\ g_{ij}(x,0) = g_{ij}(x) \ on \ M \end{cases}$$

admits a smooth solution $g_{ij}(x,t), t \in [0,T]$, with bounded curvature

Recently Chen and Zhu proved the following uniqueness Theorem.

***Theorem*** (Chen-Zhu[17]) Let $\left(M, g_{ij}(x)\right)$ be a complete non-compact Riemannian manifold of dimension n with bounded curvature. Let $g_{ij}(x,t)$ and $\overline{g_{ij}}(x,t)$ be two solutions to the Ricci flow on $M \times [0,T]$ with $g_{ij}(x)$ as the initial data and with bounded curvatures. Then $g_{ij}(x,t) = \overline{g_{ij}}(x,t)$ for $(x,t) \in M \times [0,T]$.

- **Time evolving metrics:**

**Definition**: a Normal coordinates about the point p are defined by

a)  $\gamma_v(t) = (tv^1, tv^2, \ldots, tv^n)$ is a geodesic $(v \in R^n)$
b)  $g_{ij}(p) = \delta_{ij}$
c)  $\Gamma_{ij}^k(p) = 0, \partial_i g_{jk}(p) = 0$

Suppose that $g_{ij}(t)$ is a time-dependent Riemannian metric and

$$\frac{d}{dt} g_{ij}(t) = h_{ij}(t)$$

Then various geometric quantities evolve according to the following equations:

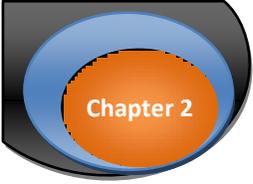

1) Metric inverse

$$\frac{\partial}{\partial t} g^{ij}(t) = -h^{ij} = -g^{ik} g^{jl} h_{kl}$$

2) Christoffel symbols

$$\frac{\partial}{\partial t} \Gamma_{ij}^{k} = \frac{1}{2} g^{kl} (\nabla_i h_{jl} + \nabla_j h_{il} - \nabla_l h_{ij})$$

3) Riemann curvature tensor

$$\frac{\partial}{\partial t} R_{ijk}^{l} = \frac{1}{2} g^{lp} \begin{Bmatrix} \nabla_i \nabla_j h_{kp} + \nabla_i \nabla_k h_{jp} - \nabla_i \nabla_p h_{jk} \\ -\nabla_j \nabla_i h_{kp} + \nabla_j \nabla_k h_{ip} - \nabla_j \nabla_p h_{ik} \end{Bmatrix}$$

4) Ricci tensor

$$\frac{\partial}{\partial t} R_{ij} = \frac{1}{2} g^{pq} (\nabla_p \nabla_i h_{jp} + \nabla_q \nabla_j h_{ip} - \nabla_q \nabla_p h_{ij} - \nabla_i \nabla_j h_{qp})$$

5) Scalar curvature

$$\frac{\partial}{\partial t} R = -\Delta H + \nabla^p \nabla^q h_{pq} - h^{pq} R_{pq}$$

Where $H = g^{pq} h_{pq}$

6) Volume element

$$\frac{\partial}{\partial t} d\mu = \frac{H}{2} d\mu$$

7) Volume of manifold

$$\frac{\partial}{\partial t} \int_M d\mu = \int_M \frac{H}{2} d\mu$$

8) Total scalar curvature on a closed manifold $M$

$$\frac{\partial}{\partial t} \int_M R d\mu = \int_M \left( \frac{1}{2} RH - h^{ij} R_{ij} \right) d\mu$$

We only prove 2), 5) and 6)

Proof2) We know

$$\frac{d}{dt} \Gamma_{ij}^{k} = \frac{d}{dt} \frac{1}{2} g^{kl} (\partial_i g_{jl} + \partial_j g_{il} - \partial_l g_{ij})$$

So we get

$$\partial_t \Gamma_{ij}^{k} = \frac{1}{2} (\partial_t g^{kl}) (\partial_i g_{jl} + \partial_j g_{il} - \partial_l g_{ij}) + \frac{1}{2} g^{kl} (\partial_i \partial_t g_{jl} + \partial_j \partial_t g_{il} - \partial_l \partial_t g_{ij})$$

Now we work in normal coordinates about a point $p$ .So according to properties b) and c) of normal coordinates we get $\partial_i g_{jk} = 0$, $\partial_i A = \nabla_i A$ at $p$ for any tensor $A$. Hence

$$\partial_t \Gamma_{ij}^k(p) = \frac{1}{2} g^{kl}\big(\nabla_i h_{jl} + \nabla_j h_{il} - \nabla_l h_{ij}\big)(p)$$

Now although the Christoffel symbols are not the coordinates of a tensor quantity, their derivative is. (This is true because the difference between the Christoffel symbols of two connections is a tensor .Thus, by taking a fixed point connection with Christoffel symbols $\widetilde{\Gamma_{ij}^k}$, we have $\partial_t \Gamma_{ij}^k = \partial_t\left(\Gamma_{ij}^k - \widetilde{\Gamma_{ij}^k}\right)$ and the right hand side is clearly a tensor).

Hence both sides of this equation are the coordinates of tensorial quantities, so it does not matter in what coordinates we evaluate them. In particular, the equation is true for any coordinates, not just normal coordinates, and about any point $p$.

Proof 5) from 4) and 1) we get

$$\partial_t R = \partial_t\big(g^{ij} R_{ij}\big) = \partial_t\big(g^{ij}\big) R_{ij} + g^{ij}\big(\partial_t R_{ij}\big)$$

$$= -h^{ij} R_{ij} + g^{ij}\left(\frac{1}{2} g^{pq}\big(\nabla_q \nabla_i h_{jp} + \nabla_q \nabla_j h_{ip} - \nabla_q \nabla_p h_{ij} - \nabla_i \nabla_j h_{qp}\big)\right)$$

$$= -\Delta H + \nabla^p \nabla^q h_{pq} - h^{pq} R_{pq}$$

(Note that $\nabla g = 0$ $and$ $\Delta = g^{ij} \nabla_i \nabla_j$).

Proof 6) we know

$$d\mu = \sqrt{\det g_{ij}}\, dx^1 \wedge dx^2 \wedge \ldots \wedge dx^n$$

and

$$\frac{d}{dt} \det A = (A^{-1})^{ij}\left(\frac{dA_{ij}}{dt}\right) \det A$$

Now by the chain rule formula we obtain

$$\partial_t d\mu = \partial_t \sqrt{\det g_{ij}}\, dx^1 \wedge dx^2 \wedge \ldots \wedge dx^n = \frac{1}{2\sqrt{\det g_{ij}}} g^{ij} h_{ij} \det g\, dx^1 \wedge dx^2 \wedge \ldots \wedge dx^n = \frac{H}{2} d\mu$$

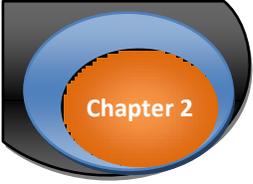

- **The normalized Ricci flow**

In case of 3-manifolds with positive Ricci curvature and higher dimensional manifolds with positive curvature operator, works by Hamilton, Huisken , Bohm and Wilking showed that the normalized Ricci flow will evolve the metric to one with constant curvature . One way to avoid collapsing is to add the condition that volume be preserved along the evaluation .To preserve the volume we need to change a little bit the equation for the Ricci flow.

We know

$$d\mu = \sqrt{det g_{ij}} \, dx^1 . dx^2 \ldots dx^n$$

Let us define the average scalar curvature

$$r = \frac{\int_M R d\mu}{\int_M d\mu}$$

Let $\tilde{g}(t) = \psi(t)g(t)$ with $\psi(0) = 1$. Let us choose $\psi(t)$ so that the volume of the manifold with respect to $\tilde{g}$ is constant. So

$$vol\big(\tilde{g}(t)\big) = vol\big(\tilde{g}(0)\big)$$

But we know if $\tilde{g} = cg$ then by applying the first corollary of volume we have $d\tilde{\mu} = c^{\frac{n}{2}}d\mu$, so

$$\psi(t)^{\frac{n}{2}} vol\big(g(t)\big) = vol\big(g(0)\big)$$

then we obtain

$$\psi(t) = \left(\frac{\int_{M^n} d\mu(t)}{\int_{M^n} d\mu(0)}\right)^{\frac{2}{-n}}$$

But according to the relation $\frac{\partial}{\partial t}\int_{M^n} d\mu = -\int_{M^n} R d\mu$ we get.

$$\frac{d}{dt}\psi(t) = \frac{2}{-n}\left(\frac{\int_{M^n} d\mu(t)}{\int_{M^n} d\mu(0)}\right)^{\frac{2}{-n}-1}\frac{\frac{d}{dt}\int_{M^n} d\mu(t)}{\int_{M^n} d\mu(0)}$$

$$= \frac{2}{n} \frac{\psi(t)}{\left(\frac{\int_{M^n} d\mu(t)}{\int_{M^n} d\mu(0)}\right)} \frac{\int_{M^n} R d\mu(t)}{\int_{M^n} d\mu(0)} = \frac{2r}{n} \psi(t)$$

The normalized average scalar curvature , $\bar{r}$, is defined by

$$\bar{r} = \frac{\int_{M^n} \tilde{R} d\mu}{\int_{M^n} d\mu} = \frac{r}{\psi(t)} \quad \left( \text{Because if } \tilde{g} = cg \text{ Then } \tilde{R} = c^{-1}R \right)$$

Where $\widetilde{R} := R(\tilde{g})$. So

$$\frac{d}{dt} \psi(t) = \frac{2\bar{r}}{n} \big(\psi(t)\big)^2$$

So we get

$$\frac{\partial}{\partial t} \tilde{g} = \left(\frac{\partial}{\partial t} g\right) \psi(t) + g \frac{d}{dt} \psi(t)$$

$$= -2\psi(t) Rc\big(g(t)\big) + \frac{2\bar{r}}{n} \big(\psi(t)\big)^2 g$$

$$= \psi(t) \left( -2Rc(\tilde{g}(t)) + \frac{2\bar{r}}{n} \tilde{g} \right)$$

We define a rescaling of time to get rid of the $\psi(t)$ terms in this evolution equation:

$$\tau = \int_0^t \psi(u) du$$

So $\frac{\partial}{\partial t} \tilde{g} = -2\widetilde{Rc} + \frac{2\bar{r}}{n} \tilde{g}$, $\tilde{g}(t)$ is called the normalized Ricci flow.

Therefore we checked that $\tilde{g}(t)$ is a solution to the normalized Ricci flow.

**Remark**: Note that $\frac{\partial}{\partial t} g_{ij} = -2R_{ij} + \frac{2}{n} r g_{ij}$ , called normalized Ricci flow, may not have solutions even for short time .[Hamilton 1-5].

- **Laplacian spectrum under Ricci flow**

Let $(M^n, g)$ be a Riemannian manifold, we know

$$\Delta f = g^{ij} \nabla_i \nabla_j f$$

$$= g^{ij} \left( \frac{\partial^2 f}{\partial x_i \partial x_j} - \Gamma_{ij}^k \frac{\partial f}{\partial x_k} \right)$$

by an eigenvalue $\lambda$ of $\Delta$ we mean there exists a non-zero function $f$ such that

$$\Delta f + \lambda f = 0$$

The set of all eigenvalues is called the spectrum of the operator it is well-known that the spectrum of the Laplacian –Beltrami operator is discrete on a compact manifold and non-negative .So we can write for the Laplacian eigenvalues is the following ways

$$0 = \lambda_0 < \lambda_1 \leq \lambda_2 \leq \lambda_3 \leq \cdots$$

Note that if $M^n$ is closed, then assuming $f \not\equiv 0$ we have

$$\lambda = \frac{\int_{M^n} |\nabla f|^2 d\mu}{\int_{M^n} f^2 d\mu}$$

Here we assume that if we have a smoothly varying one-parameter family of metrics $(t)$ , each Laplacian eigenvalue $\lambda_\alpha(g(t))$ will also vary smoothly.

***Lemma*** On a closed manifold we have

$$\int_{M^n} Rc(\nabla f, \nabla f) d\mu \leq \frac{n-1}{n} \int_{M^n} (\Delta f)^2 d\mu$$

**Theorem** . Suppose $f$ is an eigenfunction of the Laplacian with eigenvalue $\lambda$:

$$\Delta f + \lambda f = 0$$

If $Rc \geq (n-1)kg$, where $k > 0$ is a constant, then

$$\lambda \geq nk$$

Proof.  By applying lemma 2.4, we have

$$(n-1)k \int_{M^n} |\nabla f|^2 d\mu \leq \int_{M^n} Rc(\nabla f, \nabla f) d\mu \leq$$

$$\leq \frac{n-1}{n} \int_{M^n} (\Delta f)^2 d\mu = \frac{n-1}{n} \lambda^2 \int_{M^n} f^2 d\mu$$

and because $\lambda = \frac{\int_{M^n} |\nabla f|^2 d\mu}{\int_{M^n} f^2 d\mu}$, the proof is complete.

- ## **Proof of Elliptization conjecture by using Ricci flow theory**

In this subsection we give an educative proof of Elliptization conjecture by using Perlman's results and Ricci flow theory. In reality Perelman proved elliptization conjecture so Poincare conjecture of two ways and most of the efforts of mathematicians was understanding on his first proof and here we present his other proof of elliptization conjecture with complete details.

**Definition**(**connected sum**) If a closed 3-manifold $M$ contains an embedded sphere $S^2$ separating $M$ into two components, we can split $M$ along this $S^2$ into manifolds $M_1$ and $M_2$ with boundary $S^2$ .We can then fill in these boundary spheres with 3-balls to produce two closed manifolds $N_1$ and $N_2$ .One says that $M$ is the connected sum of $N_1$ and $N_2$ ,and one writes $M = N_1 \# N_2$ .This splitting operation is commutative and associative .

There is also a strong relationship between the topology of a connected sum and that of its components.

**Van Kampen theorem for connected sum**: Let X, Y be connected manifolds of dimension n. Then their connected sum $X \# Y$ is naturally decomposed into two open sets $U \cup V$ with $U \cap V \cong I \times S^{n-1} \approx S^{n-1}$. If $n > 2$ then $\pi_1(S^{n-1}) = 0$ , and hence $\pi_1(X \# Y) = \pi_1(X) * \pi_1(Y)$.

**Remark**. Let $M$ and $M'$ be connected manifolds of the same dimension

1) $M \# M'$ is compact if and only if $M$ and $M'$ are both compact

2)  $M \# M'$ is Orientable  if and only if $M$ and $M'$ are both Orientale

3) $M \# M'$ is simply connected  if and only if $M$ and $M'$ are both simply connected

**Remark**. One rather trivial possibility for the decomposition of $M$ as a connected sum is $M = M\#S^3$.

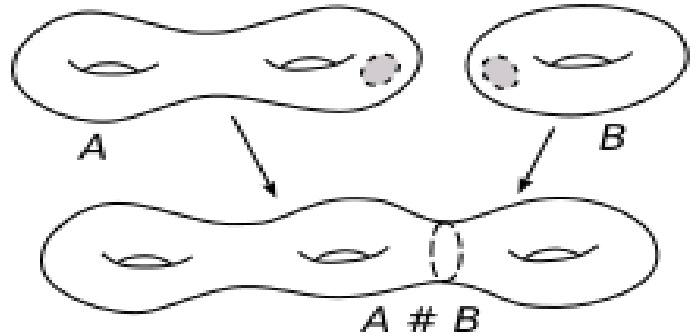

Illustration of connected sum.

**Theorem** (Perelman [16]) Let $(M, g)$ is a closed and oriented Riemannian manifold. Then there exists a Ricci flow with surgery such that every $t \in [0, \infty)$ corresponds to a Riemannian manifold $(M(t), g(t))$ and has following properties

   i.    $g(0) = g$ , $M(0) = M$

   ii.   There exists a discrete subset $T \subset (0, \infty)$ of surgery times such that if $I$ be a connected component of $[0, \infty)\setminus T$ (so an interval), and if $t_I = \inf I$ then for every $t \in \{t_I\} \cup I$ , $t \to (M(t), g(t))$ is a smooth Ricci flow and $M(t)$ is constant on this interval and we have

$$\frac{\partial g}{\partial t} = -2Rc \;\; \text{for } t \in \{t_I\} \cup I$$

   iii.  If $t \in T$, $\varepsilon > 0$ is sufficiently small then every connected component of $M(t - \varepsilon)$ is diffeomorphic with the connected sum of finite connected components of $(t)$ , $S^2$-Bundles on $S^1$ and $S^3/_\Gamma$ , where $\Gamma$ is a finite subgroup of isometries of $S^3$ .

   Moreover every connected component of $M(t)$ appears exactly in one of the connected components of $M(t - \varepsilon)$ as a summand.

**Remark**: Here we should note that according to *ii)* , surgery times are a discrete subset of $(0, \infty)$ . So in every finite time interval time we have only finite surgery and so surgery times cannot accumulate in a point.

Also for proving Elliptization conjecture we need to the following result of Perelman

**Theorem** (G.Perelman[31]) (Finite time extinction) Let $(M, g)$ be an oriented closed 3-manifold which the fundamental group of $M$ is free product of finite groups and infinite circle groups then for $t$ sufficiently large, $M(t)$ is empty.

Now we are ready to start the proof of elliptization and Poincare conjecture.

Remark : Every simply connected manifolds are orientable

**Elliptization conjecture**: Every closed 3-manifold with finite fundamental group is diffeomorphic to $S^3/\Gamma$ .

Proof : Let $(M, g)$ be a closed and oriented 3-manifold with finite fundamental group, So by applying the second theorem of Perelman, there exists a $T > 0$ (we can suppose $T$ is as surgery time) such that $M(t) = \phi$ , $\forall t \geq T$ . So because the surgery times are discrete in $[0, T]$, we have only finite surgery times which we call $t_1, t_2, ..., t_n = T$ . But by using part iii) of the first theorem in this subsection, we have

$$M(t_i) = M(t_{i+1}) \,\#\, N_1 \,\#\, \cdots \,\#\, N_k$$

such that $N_j$,s are $S^3/\Gamma$ or $S^2$-bundles on $S^1$ (which is $S^2 \times S^1$ or $S^2 \times S^1/Z_2$ ) . Now, we have

$$M(t_{n-1}) = \emptyset \,\#\, N_1 \,\#\, \cdots \,\#\, N_k$$

it means that $M(t_{n-1})$ is a connected sum of $S^3/\Gamma$ and $S^2$-bundles on $S^1$, so by continuing this process we get

$$M = M(0) = M(t_1) \,\#\, N_1 \,\#\, \cdots \,\#\, N_l \quad (*)$$

that is connected sums of $S^3/\Gamma$ and $S^2$-bundles on $S^1$ . But we know that the fundamental group of an $S^2$-bundles on $S^1$ is $\mathbb{Z}$ , so by applying Van Kampen theorem we obtain that $M$ is a connected sum of $S^3/\Gamma$. But by assumption, fundamental group of $M$ is finite and also the fundamental group of $M$ is free product of non-trivial infinite groups (from (*) and assumption). So again by applying Van Kampen theorem the number of summands should be one, so

$$M \cong S^3/\Gamma$$

So we get the desired result and the proof is complete.

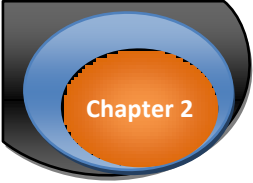

Before of proving Poincare conjecture, we need the following theorem which is essential in the proof of Poincare conjecture.

**Theorem**: Let $G$ be a connected, simply connected compact Lie group, and let $\Gamma$ be a finite subgroup of $G$. Then the homogeneous space $X = G/\Gamma$ has fundamental group $\Gamma$, which acts by right multiplication on the universal covering space $G$.

**Poincare conjecture**. If $M$ is a closed $3-$manifold with trivial fundamental group, then M is diffeomorphic to $S^3$.

Proof **:** If M is a closed manifold with trivial fundamental group then by the previous result we have $M \cong S^3/\Gamma$ , but $\Gamma$ is fundamental group of $M$ (by previous theorem and the fact that $S^3$ is simply connected) , So $\Gamma$ is the trivial group and $M \cong S^3$ and this proves the Poincare conjecture.

# Chapter 2

# PDE

The literature about elliptic equations, parabolic equations and Monge-Ampere equations are immense and it is very difficult to have a complete picture of results. Very nice books such as [] are attempt to gather and order the most significant advances in these wide fields. In this chapter we present a short review on these subjects which are important in the main chapter for our study on Calabi-Yau Theorem. Firstly we start with a short introduction on Sobolev space and Hölder space which are signification in PDE theory

- ## Sobolev Space

Broadly speaking, a Sobolev Space is a vector space of functions equipped with a norm that is a combination of $L^p -$norms of the function itself as well as its derivatives up to a given order. The derivatives are understood in a suitable week sense to make the space complete, thus a Banach Space

Definition: We say a domain $\Omega' \subset \Omega \subset R^n$ is a strictly interior subdomain of $\Omega$ and write $\Omega' \subset\subset \Omega$, if $\overline{\Omega'} \subset \Omega$ (here an open connected set $\Omega \subset R^n$ is called a domain. By $\overline{\Omega}$ we denote the closure of $\Omega$ ; $\partial\Omega$ is the boundary)

If $\Omega'$ is bounded and $\Omega' \subset\subset \Omega$, then $dist\{\Omega', \partial\Omega\} > 0$. We use the following notation:

$$x = (x_1, x_2, \dots, x_n) \in R^n, \ \ \partial_j u = \frac{\partial u}{\partial x_j}$$

$$\alpha = (\alpha_1, \alpha_2, \dots, \alpha_n) \in Z_+^n \ \text{ is a multi- index}$$

$$|\alpha| = \alpha_1 + \alpha_2 + \cdots + \alpha_n \ , \ \ \partial^\alpha u \ (or \ D^\alpha u) = \frac{\partial^{|\alpha|} u}{\partial x_1^{\alpha_1} \partial x_2^{\alpha_2} \dots \partial x_n^{\alpha_n}}$$

Next, $$\nabla u = (\partial_1 u, \dots, \partial_n u), |\nabla u| = \left(\sum_{j=1}^n |\partial_j u|^2\right)^{\frac{1}{2}}$$

Definition: $L_q(\Omega)$, (or sometimes $L^q(\Omega)$), $1 \leq q < \infty$, is the set of all measurable functions $u(x)$ in $\Omega$ such that the norm

$$\|u\|_{q,\Omega} = \left(\int_\Omega |u(x)|^q dx\right)^{\frac{1}{q}}$$

is finite . $L_q(\Omega)$ is a Banach space.

Definition: $L_{q,loc}(\Omega)$ , $1 \leq q < \infty$, is the set of all measurable functions $u(x)$ in $\Omega$ such that $\int_{\Omega'} |u(x)|^q dx < \infty$, for any bounded strictly interior subdomain $\Omega' \subset\subset \Omega$, (note that $L_{q,loc}(\Omega)$ is a topological space but not a Banach space)

We say $u_k \rightarrow u$ (when k $\rightarrow \infty$) in $L_{q,loc}(\Omega)$ , if $\| u_k - u\|_{q,\Omega'} \rightarrow 0$ (when k $\rightarrow \infty$ ) for any bounded $\Omega' \subset\subset \Omega$ .

## Weak derivatives:

Let $\alpha$ be a multi-index. Suppose that $u, v \in L_{1,loc}(\Omega)$, and

$$\int_\Omega u(x)\partial^\alpha \eta\,(x)dx = (-1)^{|\alpha|} \int_\Omega v(x)\eta\,(x)dx,$$

for all infinitely differentiable functions with compact support $\eta$ .

Then $v$ is called the weak (or distributional) partial derivative of $u$ in $\Omega$ , and denoted by $\partial^\alpha u$. If $u(x)$ is sufficiently smooth to have continuous derivative $\partial^\alpha u$, we can integrate by parts:

$$\int_\Omega u(x)\partial^\alpha \eta\,(x)dx = (-1)^{|\alpha|}\int_\Omega \partial^\alpha u(x)\eta\,(x)dx,$$

So, the classical derivative $\partial^\alpha u$ is also the weak derivative. Of course, $\partial^\alpha u$ may exist in the weak sense without existing in the classical sense.

**Remark** :

1) To define the weak derivative $\partial^\alpha u$, we don't need the existence of derivatives of the smaller order (Like in the classical definition)
2) The weak derivative is defined as an element of $L_{1,loc}(\Omega)$, So we can change it on some set of measure zero
3) Weak derivative is unique and linear ($u \to \partial^\alpha u$).

Another equivalent definition of the weak derivative.

**Definition**:  Suppose that $u, v \in L_{1,loc}(\Omega)$ and there exists a sequence $u_m \in C^l(\Omega)$ , $m \in N$ , such that $u_m \to u,\ (m \to \infty)$ and $\partial^\alpha u_m \to v\ (m \to \infty)$, in $L_{1,loc}(\Omega)$. Here $\alpha$ is a multi-index and $|\alpha| = l$ . Then $v$ is called the weak derivative of $u$ in $\Omega$: $\partial^\alpha u = v$.

For better understanding of this definition we give an example:

**Example**: Take $n = 1, \Omega = [-,1]$ and $f(x) = 1 - |x|$. We show that the weak derivative is given by

$$g(x) := f(x) = \begin{cases} 1, & x < 0 \\ -1, & x > 0 \end{cases}$$

To show this, we break the interval $[-1,1]$ into the two parts in which $f$ is smooth, and integrate by parts

$$\int_{-1}^{1} f(x)\,\varphi'(x)dx = \int_{-1}^{o} f(x)\,\varphi'(x)dx + \int_{0}^{1} f(x)\,\varphi'(x)dx$$

$$= f\varphi|_{-1}^{0} - \int_{-1}^{0} (+1)\,\varphi(x)dx + f\varphi|_{0}^{1} - \int_{0}^{1} (-1)\,\varphi(x)dx$$

$$= -\int_{-1}^{1} f(x)\,\varphi(x)dx + f(0^-)\varphi(0^-) - f(0^+)\varphi(0^+) = -\int_{-1}^{1} g(x)\,\varphi(x)dx$$

because $f$ is continuous at 0.

**Definition** : For $p \in [1, \infty]$ , $k \in N$ and $\Omega$ an open subset of $R^l$, let

$$W_{loc}^{k,p}(\Omega) := \left\{ f \in L_{loc}^p(\Omega) : \partial^\alpha f \in L_{loc}^p(\Omega) \text{ (weakly) for all } |\alpha| \leq k \right\},$$

$$W^{k,p}(\Omega) := \left\{ f \in L^p(\Omega) : \partial^\alpha f \in L^p(\Omega) \text{ (weakly) for all } |\alpha| \leq k \right\}$$

$$\|f\|_{W^{k,p}(\Omega)} := \left( \sum_{|\alpha| \leq k} \|\partial^\alpha f\|_{L^p(\Omega)}^p \right)^{\frac{1}{p}} \quad if \ p < \infty$$

$$\|f\|_{W^{k,p}(\Omega)} := \left( \sum_{|\alpha| \leq k} \|\partial^\alpha f\|_{L^\infty(\Omega)} \right) \quad if \ p = \infty$$

In the special case of $p = 2$, we write $W_{loc}^{k,2}(\Omega) =: H_{loc}^k(\Omega)$ and $W^{k,2}(\Omega) =: H^k(\Omega)$ in which case $\|.\|_{W^{k,2}(\Omega)} = \|.\|_{H^k(\Omega)}$ is a Hilbertian norm associated to the inner product

$$(f,g)_{H^k(\Omega)} = \sum_{|\alpha| \leq k} \int_\Omega \partial^\alpha f . \overline{\partial^\alpha g} \, dm$$

Also the function, $\|.\|_{W^{k,p}(\Omega)}$, is a norm which makes $W^{k,p}(\Omega)$ into a Banach space.

# PDE

PDE is a differential equation that contains unknown multivariable functions and their partial derivatives. PDEs are used to formulate problems involving functions of several variables. When one assume that a solution to a PDE lies in a specific function class one is in fact imposing a regularity condition on solutions to the PDE. Proving an a priori estimate amounts to showing that under the basic assumption that a solution to a PDE belongs to a given function class that we can in fact find a uniform bound for all solutions of that function class. For more details see [12], [22] and [33].

Let $U$ be an open subset of $R^n$ and $u: U \to R$ be a function ($U$ assumed to be connected and bounded) we define operator $Du$ as follows

$$Du = (D_1 u, D_2 u, \ldots, D_n u) \text{ and } D_i u = \frac{\partial u}{\partial x_i}, D_{ij} u = \frac{\partial^2 u}{\partial x_i \partial x_j}, \ldots$$

Also for any $k \in N$, denote

$$D^k u(x) = \{D^\alpha u(x) : |\alpha| = k\}$$

The set of all partial derivatives of order k. And we also regard $D^k u(x)$ as a point in $R^{n^k}$ and

$$\left| D^k u(x) \right|^2 = \sum_{|\alpha|=k} |D^\alpha u(x)|^2$$

It is clear that when $k = 1$, $Du = (D_1 u, D_2 u, \ldots, D_n u)$ is the gradient operator.

**Definition** (Hessian of $u$): When $n = 2$, $D^2 u = (D_{ij} u)$ as an $n \times n$ matrix, is called Hessian of $u$.

## • Second order elliptic PDEs

The application of linear and non-linear second order elliptic PDE's are numerous in many fields, in particular in geometric analysis and differential geometry. We are mainly concerned with general properties of solutions of linear equations as follows

$$Lu \equiv a^{ij}(x)D_{ij}u + b^i(x)D_i u + c(x)u = f(x) \quad i,j = 1,2 \ldots, n \quad (2.1)$$

Definition: When $f = 0$, we call $Lu = 0$ a homogeneous equation and the ellipticity of these equations is expressed by the fact that the coefficient matrix $[a^{ij}]$ is positive definite in the domain.

Definition: Let $\Omega \subset R^n$ be an open (bounded) subset and $u \in C^2(\Omega)$. The Laplacian of $u(x)$, denoted by $\Delta u$, is defined by

$$\Delta u = \sum D_{ii} u = div(Du)$$

Note that $\Delta u = Tr(D^2 u)$ .

Also if $u$ is the solution of the Laplace equation. $\Delta u = 0$, we say $u$ is harmonic. We recall that $L$ defined in (2.1) is elliptic at a point $x \in \Omega$ if the matrix $\left(a^{ij}(x)\right)$ is positive definite; that is there exists $\lambda(x), \Lambda(x)$ such that

$$0 < \lambda(x)|\xi|^2 \leq a^{ij}(x)\xi_i\xi_j \leq \Lambda(x)|\xi|^2.$$

For any nonzero vector $\xi$. $L$ is called strictly elliptic if $\lambda(x) \geq \lambda_0 > 0$ for some positive constant $\lambda_0$ . Also $L$ is said to be uniformly elliptic if

$$\lambda_0|\xi|^2 \leq a^{ij}(x)\xi_i\xi_j \leq \Lambda_0|\xi|^2$$

for positive constants $\lambda_0, \ \Lambda_0$ .

Now we introduce some function classes that will be of importance in our analysis.

- **Linear Parabolic equations**

The standing example of linear parabolic equations with constant coefficients is the heat equation

$$\frac{\partial u}{\partial t} - \Delta u = f$$

➤ **Hölder Space**

**Definition**: A function $f$ is uniformly Hölder continuous with exponent $0 < \alpha < 1$ if

$$[f]_{\alpha;\Omega} = \underset{\substack{x,y \in \Omega \\ x \neq y}}{\text{Sup}} \frac{|f(x) - f(y)|}{|x - y|^\alpha} < \infty$$

We define the Hölder space $C^{k,\alpha}(\overline{\Omega})$ $(C^{k,\alpha}(\Omega))$ to be the space consisting of functions whose k$^{th}$ order partial derivatives are uniformly Hölder continuous. Also $C^{k,0}(\overline{\Omega}) = C^k(\overline{\Omega})$ where $C^k(\overline{\Omega})$ is the space of continuous functions whose k$^{th}$ order partial derivatives are continuous up to the boundary. We define the following norms.

$$\|u\|_{C^k(\overline{\Omega})} = \sum_{j=0}^{k} [D^j u]_{0,\overline{\Omega}} \ ,$$

$$\|u\|_{C^{k,\alpha}(\overline{\Omega})} = \|u\|_{C^k(\overline{\Omega})} + [D^k f]_{\alpha,\overline{\Omega}}$$

When $= 0$ , we usually use $C^\alpha$ for $C^{0,\alpha}$ since there is no ambiguity for $0 < \alpha < 1$. Also we have the following inequality

$$|f_1 f_2|_{C^\alpha} \leq \left(\sup_\Omega |f_1|\right) |f_2|_{C^\alpha} + \left(\sup_\Omega |f_2|\right) |f_1|_{C^\alpha}$$

In short, $C^\alpha$ is an algebra.

Now, here we state a fundamental existence and uniqueness result for linear parabolic equation with Holder continuous coefficients.

**Theorem**. If $\Omega$ is a $C^{2,\alpha}$ domain and the coefficients $A, b, c \in C^\alpha((0,T) \times \Omega$ ) and $f \in C^\alpha([0,+\infty) \times R^d)$, $g \in C^{2+\alpha}((0,T) \times \Omega)$, $h \in C^{2,\alpha}(R^d)$, and $g$ and $h$ are compatible (The fact that data $g$ and $h$ are compatible has to do with conditions ensuring that a solution which is regular up to the boundary can be constructed), then there exists a unique solution $u \in C^{2,\alpha}(Q)$ (which $Q \subseteq (0,T) \times \Omega$ ) of

$$\begin{cases} \dfrac{\partial u}{\partial t} - \Delta u = f & \text{in } (0,T) \times \Omega \\ u = g & \text{on} \quad (0,+\infty) \times \partial\Omega \\ u = h & \text{on} \quad \{0\} \times \overline{\Omega} \end{cases}$$

# Maximum principle in PDE

The Maximum principle in PDE says that the maximum of a harmonic function in a domain is to be found on the boundary of that domain. Roughly speaking, the strong maximum principle says that if achieves its maximum in the interior of the domain, the function is uniformly a constant. The weak maximum principle says that the maximum of the function is to be found on the boundary, but may re-occur in the interior as well.

Recall

$$L(u) \equiv \sum_{i,j=1}^{n} a^{ij}(x) \frac{\partial^2 u}{\partial x_i \partial x_j} + \sum_{i=1}^{n} b^i(x) \frac{\partial u(x)}{\partial x_i} + c(x)u(x) = f(x)$$

In some domain $\Omega \subset R^n$.

Lemma (Maximum principle [12], [32-33]) Assume $(x) \equiv 0$ , and let $u$ satisfy in $\Omega$,

$$L(u) \geq 0$$

that is

$$\sum_{i,j=1}^{n} a^{ij}(x) \frac{\partial^2 u}{\partial x_i \partial x_j} + \sum_{i=1}^{n} b^i(x) \frac{\partial u(x)}{\partial x_i} \geq 0$$

Then $\sup_{x \in \Omega} u(x) = \max_{x \in \partial\Omega} u(x)$ . In the case $L(u) \leq 0$ , a corresponding result holds with sup/max replaced by inf/min . A consequence is the uniqueness of solutions when $(x) \equiv 0$ .

So by letting $\Omega^+ = \{x \in \Omega: \ u(x) \geq 0\}$ and applying the previous lemma we get following corollary.

**Corollary**: Suppose $c(x) \leq 0$ in $\Omega$ . Let $u \in C^2(\Omega) \cap C^0(\overline{\Omega})$ satisfy $L(u) \geq 0$ in $\Omega$. Write $u^+(x) \equiv \max(u(x), 0)$, we then have

$$\sup_{\Omega} u^+ \leq \max_{\partial\Omega} u^+$$

Theorem: Suppose $c(x) \equiv 0$, let $u$ satisfy

$$L(u) = 0 \quad \text{in } \Omega$$

if $u$ attains its maximum in the interior of $\Omega$, then it has to be constant. If $c(x) \leq 0$, then $u$ has to be a constant if it attains a nonnegative interior maximum .

**Lemma**([32]): Suppose $c(x) \leq 0$ and $L(u) \geq 0$ in $\Omega \subseteq R^n$, and let $x_0 \in \partial\Omega$ , Moreover, assume

    I.    $u$ is continuous at $x_0$

    II.    $u(x_0) \geq 0$ if $c(x) \not\equiv 0$

    III.    $u(x_0) > u(x)$ for all $x \in \Omega$

    IV.    There exists an open ball $B_r(y) \subset \Omega$ with $x_0 \in \partial B_r(y)$ and we have

$$\frac{\partial u}{\partial n}(x_0) > 0$$

where $n$ is the outer normal of the ball $B_r(y)$ at $x_0$ , provided that this derivative exists

- **Maximum principle for parabolic equations**

Consider the operator

$$L(u) \equiv \sum_{i,j=1}^{n} a_{ij}(x,t) \frac{\partial^2 u}{\partial x_i \partial x_j} + \sum_{i=1}^{n} b_i(x,t) \frac{\partial u}{\partial x_i} + c(x,t)u - \frac{\partial u}{\partial t}$$

$$A(u) = \sum_{i,j=1}^{n} a_{ij}(x,t) \frac{\partial^2 u}{\partial x_i \partial x_j} + \sum_{i=1}^{n} b_i(x,t) \frac{\partial u}{\partial x_i} + c(x,t)u$$

in $\Omega_T = \Omega \times (0, \text{T}]$, with $T > 0$, and $\Omega$ domain in $R^n$, (open and bounded)

A) We say that $L$ is parabolic in $\Omega_T$, if there exists $\lambda > 0$ such that for every $(x,t) \in \Omega_T$ and for any real vector $\xi \neq 0$ ,

$$\sum_{i,j=1}^{n} a_{ij}(x,t)\xi_i\xi_j > \lambda|\xi|^2$$

B) We assume that the coefficients in $L$ are bounded functions in $\Omega_T$

**Definition**: We define the set $C^{(2,1)}(\Omega_T)$ as follows

$$C^{(2,1)}(\Omega_T) = \left\{ u : \Omega_T \to R \; ; u, u_t, \frac{\partial u}{\partial x_i}, \frac{\partial^2 u}{\partial x_i \partial x_j} \in C(\Omega_T) \right\}$$

**Theorem([12],[17] and [18])**: Let A) and B) hold and $= 0$ . If $u \in C^{(2,1)}(\Omega_T) \cap C(\overline{\Omega_T})$ satisfies $(u) = A(u) - u_t \geq 0$ , then

$$\sup_{\Omega_T} u = \max_{\overline{\Omega_T}} u = \max_{\partial\Omega_T \setminus \Omega \times \{T\}} u$$

**Theorem ([12],[17] and [18])**: (weak maximum principle for the parabolic equation): Let A) and B) hold and $\leq 0$ . If $u \in C^{(2,1)}(\Omega_T) \cap C(\overline{\Omega_T})$ satisfies $(u) = A(u) - u_t \geq 0$ , then

$$\sup_{\Omega_T} u = \max_{\overline{\Omega_T}} u \leq \max_{\partial_p \Omega_T \setminus \overline{\Omega} \times \{T\}} u^+$$

where $u = u^+ - u^-$, $u^+ = \max(u, 0)$

For equations of parabolic type, a simple version says the following.

Theorem (**weak maximum principle for super solutions of the heat equation**). Let $g(t)$ be a family of metrics on a close manifold $M^n$ and let $u : M^n \times [0, T) \to R$ satisfy

$$\frac{\partial u}{\partial t} \geq \Delta_{g(t)} u.$$

Then if $u \geq c$ at $t = 0$ for some $c \in R$, then $u \geq c$ for all $t \geq 0$.

Proof. The idea is simply that given a time $t_0 \geq 0$, if the spatial minimum of $u$ is attained at a point $x_0 \in M$, then

$$\frac{\partial u}{\partial t}(t_0, x_0) \geq \Delta_{g(t)} u(t_0, x_0) \geq 0$$

so that the minimum should be non decreasing. Note that at $(t_0, x_0)$ we actually have $(\nabla_i \nabla_j u) \geq 0$. More rigorously, we proceed as follows. Given any $\varepsilon > 0$, define $u_\varepsilon : M \times [0, T) \to R$,

$$u_\varepsilon = u + \varepsilon(1 + t)$$

since $u \geq c$ at $t = 0$, we have $u_\varepsilon > c$ at $t = 0$. Now suppose for some $\varepsilon > 0$ we have $u_\varepsilon \leq c$ somewhere in $M \times (0, T)$. Then since M is closed, there exists $(t_1, x_1)$ such that $u_\varepsilon(t_1, x_1) = c$ and $u_\varepsilon(t, x) > c$ for all $x \in M$ and $t \in M$ and $t \in [0, t_1)$. we then at $(t_1, x_1)$ have

$$0 \geq \frac{\partial u_\varepsilon}{\partial t} \geq \Delta_{g(t)} u_\varepsilon + \varepsilon > 0$$

which is a contradiction. Hence $u_\varepsilon > c$ on $M \times [0, T)$ for all $\varepsilon > 0$ and by taking the limit as $\varepsilon \to 0$ we get $u \geq c$ on $\times [0, T)$ . So we obtain the desired result and the proof is complete.

## ➢ **Monge Ampere equation**

Monge-Ampere equations are partial differential equations whose leading term involves the determinant of the Hessian of the unknown function. Historically, the study of Monge-Ampere equation is very much motivated by the following two problems: Minkowski problem and Weyl problem. One is of prescribing Gauss curvature type, another is of embedding type. These equations play an important role in geometry, because fundamental geometric notions (such as volumes and curvatures) are given in various contexts by determinants of Hessians. In particular, the complex Monge-Ampere equation is of considerable interest in Kähler geometry and Kähler ricci flow. One of the important connections between Monge-Ampere equations and geometry is exhibited by the following face. If $M = \{(x, u(x)) | x \in \Omega \subseteq Rn$ is a graph of a function $u(x)$, the Gauss curvature of the (graph) surface $k(x, u(x))$ satisfies

$$\det(D^2 u) = k(x, u(x))(1 + |Du|^2)^{\frac{n+2}{2}}$$

**Definition** (**Real Monge-Ampere equation**) Let $\varphi: R^n \to R$, be a twice differentiable function, $\frac{\partial^2 \varphi}{\partial x_i \partial x_j}$, its Hessian matrix. The real Monge- Ampere equation is given by

$$\det\left(\frac{\partial^2 \varphi}{\partial x_i \partial x_j} - A(x, \varphi, d\varphi)\right) = F(x, \varphi, d\varphi)$$

where $\varphi$ is unknown, and $F$ a given function.

**Remark**: It is elliptic, if $\varphi$ is convex, and the matrix $A$ is positive definite.

- **Solving the Monge-Ampere equation**

1. Uniqueness of solutions(on compacts or with prescribed boundary conditions)
2. Existence of weak solutions (solutions which are generalized functions, that is, with singularities)
3. Elliptic regularity (every weak solution is in fact smooth and real analytic)

One of the tricks for solving the Monge-Ampere equation is the continuity method of S.-T.Yau which we explain here without details and in last chapter we will use of this method with complete details for proving short time existence on Kähler Ricci flow. For more details see [25] and [30].

- **Continuity method of S.-T.-Yau**

A) Suppose we have a Monge-Ampere equation $MA(\varphi) = F_t$ depending on $t \in [0,1]$. Solve $MA(\varphi) = F_t$ for $= 0$. We have this fact that the set of all $t$ for which one can solve $MA(\varphi) = F_t$ is open and closed

B) A limit of solutions of $MA(\varphi) = F_t$ is a weak solution.

- **Complex Monge-Ampere equation**

**Definition**: Let $\varphi$ be a function on $\mathbb{C}^n$, and $dd^c\varphi$ its complex Hessian, $dd^c\varphi = Hess(\varphi) + I(Hess(\varphi))$. It is a Hermitian form.

We recall that A Kähler manifold is a complex manifold with a Hermitian metric $g$ which is locally represented as $g = dd^c\psi$.

**Definition**: Let $(M, g)$ be a Kähler manifold. The complex Monge-Ampere equation is

$$\det(\mathrm{g} + dd^c\varphi) = e^f$$

Also we after prove that on a complex Kähler manifold along a Kähler Ricci flow, The complex Monge-Ampere equation has a unique solution, for any smooth function $f$ subject to constraint

$$\int_M e^f Vol_g = \int_M Vol_g$$

# Chapter 3

## Almost Complex Manifolds

In this section, some basic definitions and facts about Hermitian geometry and Kähler geometry are stated. More information regarding this can be found in [],[].

A real differentiable manifold is a topological space which locally like an open subset of $R^n$ . Now we will look at extra structure which we may impose on differentiable manifolds, namely that of a complex structure. Differentiable manifolds with a complex structure are topological spaces which locally like a part of $\mathbb{C}^n$ rather than $\mathbb{R}^n$ . For more details see [19], [20] and [28].

**Definition:** A complex manifold of complex dimension $m$ is a topological manifold $(M, U)$ whose atlas $(\phi_U)_{u \in U}$ satisfies the following compatibility condition:

$$\phi_{UV} := \phi_U o \phi_{V^{-1}}$$

is holomorphic . A pair $(U, \phi_U )$ is called a chart and the collection of all charts is called a holomorphic structure.

**Definition**: A function $F = f + ig : U \subset \mathbb{C} \to \mathbb{C}$ is called holomorphic, if it satisfies the Cauchy-Riemann equation

$$\frac{\partial f}{\partial x} = \frac{\partial g}{\partial y} \quad \text{and} \quad \frac{\partial f}{\partial y} + \frac{\partial g}{\partial x} = 0$$

Now we define holomorphic in $\mathbb{C}^n$.

We identify $\mathbb{C}^m$ with $R^{2m}$ via

$$(z_1, z_2, \ldots, z_m) = (x_1 + iy_1, x_2 + iy_2, \ldots, x_m + iy_m) \rightarrow (x_1, \ldots, x_m, y_1, \ldots y_m)$$

and denote by $J_m$ the endomorphism of $R^{2m}$ corresponding to multiplication by $i$ on $\mathbb{C}^m$.

$$J_m := \begin{pmatrix} 0 & -I_m \\ I_m & 0 \end{pmatrix}$$

A map $F: U \subset \mathbb{C}^n \rightarrow \mathbb{C}^m$ is holomorphic if and only if the differential $F_*$ of $F$ as real map $F: U \subset R^{2n} \rightarrow R^{2m}$ satisfies $J_m o (F_*)_p = (F_*)_p o J_n$ , $\forall p \in U$ .( here $(F_*)_p$ is derivative of $F$ in p).

Since every holomorphic map between open sets of $\mathbb{C}^m$ is in particular a smooth map between open sets of $R^{2m}$ every complex manifold $M$, of complex dimension $m$, defines a real $2m$-dimensional smooth manifold $M_\mathbb{R}$ , which is the same as $M$ as topological space. The converse does not hold.

**Definition** :Let $V$ be a vector space . A complex structure on $V$ is a Linear map $J: V \rightarrow V$ with $J^2 = -Id$

For every $\in T_x M_R$ , choose $U \in \mathcal{U}$ containing $U$ and define

$$J_U(X) = (\phi_U)_*^{-1} o J_m o (\phi_U)_* \ (X)$$

If we take some other $V \in \mathcal{U}$ containing $U$, then $\phi_{VU} = \phi_V o \phi_U^{-1}$, is holomorphic, and $\phi_V = \phi_{VU} o \phi_U$, So

$J_V(X) = (\phi_V)_*^{-1} o J_m o (\phi_V)_* \ (X) = (\phi_V)_*^{-1} o J_m o (\phi_{VU})_* . o (\phi_U)_* (X) = (\phi_V)_*^{-1} o (\phi_{VU})_* o J_m o (\phi_U)_* (X) = o J_m o (\phi_U)_* \ (X) = J_U(X),$

Showing that $J_U$ does not depend on $U$. Their collection is thus a well-defined tensor $J$ on $M_R$ satisfying $J^2 = -Id$.

**Definition:** Let $M$ be a real $2n$-dimensional manifold, and we can define a smooth field of complex structures on the tangent bundle $J: T_p M \rightarrow T_p M$ where $J^2 = -Id_{T_p M}$, we make the underlying vector space into a complex vector space by setting $(a + ib).v = a.v + b.J(v),$ for $a, b \in R$ and $v \in TM$. Then $J$ is an almost complex structure on the manifold and $(M, J)$ is an almost complex manifold.

Note that the condition $J^2 = -Id$ requires the real dimension of $M$ to be even because suppose $M$ is $n$-dimensional and let $J: TM \to TM$ be an almost complex structure. Then $det(J - xId)$ is a polynomial in $x$ of degree $n$. If $n$ is odd, then it has a real root . Then $det(J - zId) = 0$, so there exists a vector $v$ in $TM$ with $Jv = zv$. Hence $JJv = z^2v$ which is clearly not equal to $-v$, since $z$ is real thus $n$ must be even.

**Example**: Here we give three different examples of complex manifolds

❖ Let , $M = \mathbb{C}^n$ then $(U, \psi) = (\mathbb{C}^n, Id_{\mathbb{C}^n})$ is a chart on $\mathbb{C}^n$

❖ Let us consider the manifold given by the projective space $\mathbb{C}P^n$ . Let $U_k = \{[z_0 : z_1 : \ldots : z_n] \in CP^n : z_k \neq 0\}$ for each $= 0,1, \ldots n$ . Then $U_k$ is open in $\mathbb{C}P^n$. We define a map $\psi_k : U_k \to \mathbb{C}^n$ by

$$\psi_k([z_0 : z_1 : \ldots : z_{k-1} : z_k : z_{k+1} : \ldots : z_n]) = \left(\frac{z_0}{z_k}, \ldots, \frac{z_{k-1}}{z_k}, \frac{z_{k+1}}{z_k}, \ldots, \frac{z_n}{z_k}\right)$$

and $\{(U_k, \psi_k) : k = 0,1, \ldots, n\}$ is a complex chart on $\mathbb{C}P^n$ .

❖ Let $n > 1$, and $a \in \mathbb{C}^*$, $0 < |a| < 1$. Let $\mathbb{Z}$ act on $\mathbb{C}^n - \{0\}$ by holomorphic transformation $(z_1, \ldots, z_n) \to (a^k z_1, \ldots, a^k z_n)$ with $k \in \mathbb{Z}$ . Then the quotient $\left(\mathbb{C}^n - \{0\}\right)/\mathbb{Z}$ is diffeomorphic to $S^1 \times S^{2n-1}$. Hence $S^1 \times S^{2n-1}$ is a complex manifold. This construction is due to Hopf and the obtained manifold is named the Hopf manifold . After him Calabi and Eckmann extended this result to show that the product of odd dimensional spheres are complex manifolds.

It it interesting that $S^6$ admits an almost complex structure. In fact $S^2$ and $S^6$ are the only even dimensional spheres which admit almost complex structures. Also the following conjecture is a long-standing unsolved question.

**Conjecture : (Existence of complex structure on $S^6$ ) :** Does $S^6$ admits a complex structure.? More generally, One may ask which closed even dimensional manifolds admit almost complex structures but not complex structures . see [19]

**Example:** Let $M = \mathbb{C}P^n = \{[z_0 : z_1 : \ldots : z_n]; 0 \neq (z_0, \ldots, z_n) \in \mathbb{C}^{n+1}\}$ and let $U_0 = \{[1, z_1 : \ldots : z_n]\} \simeq \mathbb{C}^n$ be an open subset of $\mathbb{C}P^n$ . Set

$$g_{i\bar{j}} = \frac{\partial^2}{\partial z_i \partial \bar{z}_j} \log(1 + |z_1|^2 + \cdots + |z_n|^2)$$

or equivalently

$$\omega_g = \frac{\sqrt{-1}}{2}\partial\bar{\partial}\log(1+|z_1|^2+\cdots+|z_n|^2) = \frac{\sqrt{-1}}{2}\left(\frac{dz_i \wedge d\bar{z}_i}{1+|z|^2} - \frac{\bar{z}_i dz_i \wedge z_j d\bar{z}_j}{(1+|z|^2)^2}\right)$$

So we need to show that $\omega_g$ is globally defined on $\mathbb{C}P^n$. Let

$$U_1 = \{(\omega_0, 1, \omega_2, \ldots, \omega_n)\} \subset \mathbb{C}p^n$$

Then

$$U_0 \cap U_1 = \{[z_0:z_1:\ldots:z_n] = [\omega_0:1:\omega_2:\ldots:\omega_n]\}.$$

there exists a non-zero constant $\lambda$ such that $\lambda\omega_0 = 1$, $\lambda = z_1$ and $\lambda\omega_i = z_i$ for $i = 2, \ldots, n$, so $\lambda = \frac{1}{\omega_0}$ and therefore $z_i = \frac{\omega_i}{\omega_0}$ for all $i \neq 1$ and $z_1 = \frac{1}{\omega_0}$. So we get

$$\omega_g = \frac{\sqrt{-1}}{2}\partial\bar{\partial}\log(1+|z_1|^2+\cdots+|z_n|^2) = \frac{\sqrt{-1}}{2}\partial\bar{\partial}\log\left(1+\frac{1}{|\omega_0|^2}+\sum\frac{|\omega_i|^2}{|\omega_0|^2}\right)$$

$$= \frac{\sqrt{-1}}{2}\partial\bar{\partial}\log(1+|\omega|^2) - \partial\bar{\partial}\log(|\omega_0|^2) = \frac{\sqrt{-1}}{2}\partial\bar{\partial}\log(1+|\omega|^2)$$

since $\omega_0$ is holomorphic. So $\omega_g$ is globally defined and the corresponding metric on $\mathbb{C}P^n$ is called the Fubini-Study metric.

**Remark**: naturally $SU(n+1)$ acts on $\mathbb{C}P^n$ transitively

## ❖ Kähler Structure and the complexified tangent bundle

If we compare Riemannian manifolds with Kähler manifolds we see when we construct a new structure on Riemannian manifolds our results is completely different on this new structure, more precisely , a well–known theorem due to John Nash states that any Riemannian manifold can be isometrically immersed into the real Euclidean space $R^N$ for sufficiently large $N$, a Kähler manifold does not always admit an isometric and holomorphic (from now on Kähler) immersion into the complex Euclidean space $C^N$ (not even if $N$ is allowed to be infinite). For example, a compact manifold cannot be holomorphically immersed into $C^N$ for any value of $N$, since every holomorphic function from a connected compact set into $C$ is constant. But even if we consider noncompact manifolds, there are still many obstructions to the existence of such an immersion. So we are ready to construct our structure.

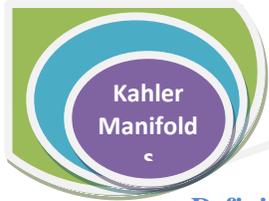

**Definition**. A Kähler structure on a Riemannian manifold $(M^n, g)$ is given by a 2-form $\Omega$, called Kähler form and a field of endomorphism $J$ of the tangent bundle, satisfying the following

- Algebraic conditions

  a. $J$ is an almost complex structure : $J^2 = -Id$

  b. The metric is almost Hermitian with respect to $J$ :
  $$g(X, Y) = g(JX, JY), \quad \forall X, Y \in TM$$

  c. $\Omega(X, Y) = g(JX, Y)$

- Analytic conditions

  d. The 2-form $\Omega$ is closed $d\Omega = 0$

  e. $J$ is integrable in the sense that its Nijenhuis tensor vanishes (to every almost complex structure $J$ one can associate a (2,1)-tensor $N^J$ called the Nijenhuis tensor , defined by $N^J(X, Y) = [X, Y] + J[JX, JY] + J[X, JY] - [JX, JY]$)

**Remark**: Let $G$ be the matrix representation of $g$. The $J-$invariance implies $G = J^t G J$, i.e.

$$G := \begin{bmatrix} A & B \\ B^t & D \end{bmatrix} = \begin{bmatrix} i.I & 0 \\ 0 & -i.I \end{bmatrix} \begin{bmatrix} A & B \\ B^t & D \end{bmatrix} \begin{bmatrix} i.I & 0 \\ 0 & -i.I \end{bmatrix} = \begin{bmatrix} -A & B \\ B^t & -D \end{bmatrix}$$

hence, $A = D = 0$.

- **Identifying the compact Kähler manifolds:**

In complex dimension one every manifold is Kahler, the reason is because the exterior derivative of any 2-form is zero, so every Hermitian metric is Kahler. In complex dimension two a surface is Kahler if and only if its first Betti number is even (Informally, the Betti number is the maximum number of cuts that can be made without dividing a surface into two separate pieces. Formally, the n-th Betti number is the rank of the n- th homology group of a topological space. The following table gives the Betti number of some common surfaces.). This is known since the 80's by classification of surfaces and hard results of Siu. The condition that the odd Betti numbers of the manifold be even for it to be Kahler is necessary by the Hodge decomposition theorem. Interestingly this is sufficient in dimension two, implying that the Kahler condition is topological there. This fails in higher dimensions, for there is an example due to Hironaka of Kahler manifolds of dimension 3 that can be deformed to a non-Kahler manifold. Since the underlying smooth manifolds of these complex manifolds are all diffeomorphic, a sufficient condition for being a Kahler manifold cannot be read off the topological or smooth structure of a given complex manifold. The list of necessary conditions that a manifold must satisfy to be Kahler is long and growing. First among these are the multiple conditions that the Hodge decomposition and the hard Lefschetz theorem impose on its cohomology ring. The properties of this ring are not yet well understood, considering how recently Voisin constructed examples demonstrating that there are Kahler manifolds whose homotopy type differs from that of projective manifolds. Serre show that any finite group can be the fundamental group of a Kahler manifold (projective surface, even). Finer results are available under additional hypotheses, for example Paun has shown that the fundamental group of a Kahler manifold with nef anticanonical bundle has polynomial growth.

**Remark**: From the definition of Kahler manifold, one can see Kahler geometry is a class of Hermitian geometry with one extra condition $d\omega = 0$. Also note that if $\omega$ is Kahler form on a compact complex manifold $M$ then the $(p,p)$ −form $\omega^p$, is not exact, because if $\omega^p = d\alpha$ for some $\alpha$, then

$$\int_M \omega^n = \int_M d(\alpha \wedge \omega^p) = 0$$

which is a contradiction. Since $\omega^p$ is real closed $2p$-form, it follows that for compact Kahler Manifolds $H^{2p}(M, \mathbb{R}) \neq 0$.

 In order to emphasize this, the following is a well known example of non-Kahler Hermitian manifold.

**Example** (**Hopf Surface**). Let $\phi \colon \mathbb{C}^2 \backslash \{0\} \to \mathbb{C}^2 \backslash \{0\}$ defined by $\phi(z) = 2z$. Denote $<\phi>$ to be the group generated by the automorphism $\phi$ of $\mathbb{C}^2 \backslash \{0\}$. One can verify that the quotient $\mathbb{C}^2 \backslash \{0\}/<\phi>$ has the complex manifold structure. This manifold is called Hopf surface. It can be proven that Hopf surface does not admit any Kähler structure. In reality Hopf surface, topologically equivalent to $S^1 \times S^3$. Therefore $H^2(M, \mathbb{R}) = 0$, so $M$ is not Kähler.

Let $(M, J)$ be an almost complex manifold. We would like to diagonalize the endomorphism. In order to do so, we have to complexify the tangent space.

Let $V$ be a real vector space, let $J$ be a complex structure on. We extend $J$ complex linearly to the complexification $V_C = V \otimes_R C$ of, $(v \otimes \alpha) = J(v) \otimes \alpha$ .

Then we still have $J^2 = -Id$, hence $V_C$ is the sum

$$V_C = V' \oplus V'' \qquad (3.1)$$

Of the eigenspace $V'$ and $V''$ for the eigenvalues $Id$ and $-Id$, respectively . The maps

$$V \to V' \qquad\qquad\qquad V \to V''$$

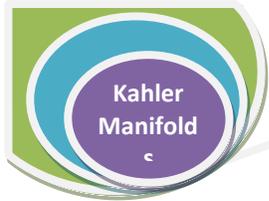

$$v \to v' := \frac{1}{2}(v - iJv) \qquad\qquad v \to v'' := \frac{1}{2}(v + iJv)$$

Are complex linear and conjugate linear isomorphism, respectively

Now, let $M$ be a smooth manifold with an almost complex structure $J$ and $T_C M = TM \otimes_R C$ be the complexified tangent bundle . As in (4.1) we have the

$$TM \to T'M \qquad , \qquad TM \to T''M$$

$$v \to v' \qquad\qquad v \to v''$$

Proposition: If $M$ has an almost complex structure $J$, then $J$ induces a splitting of the complexified tangent and cotangent bundle $TM \otimes_R C$ and $T^*M \otimes_R C$ such that we can decompose the complex tangent bundle as

$$T_p M \otimes_R C \cong T_p^{1,0} M \oplus T_p^{0,1} M$$

Where $T^{1,0}M = \{v \in TM \otimes C | Jv = iv\}$ and $T^{0,1}M = \{v \in TM \otimes C | Jv = -iv\}$ are eigenspaces for the eigenvalues $i$ and $-i$ of $J$ respectively

**Definition**: Suppose $M$ is a *2n*-dimensional manifold and $J$ an almost complex structure on M. Let $f: M \to C$ be a smooth function and write $f = u + iv$. We call $f$, $J$−holomorphic if $du = Jdv$, i.e., $dfoJ = idf$ .

Note that the splitting of the complexified tangent bundles of a manifold $M$ produced by the almost complex structure $J$ does not guarantee the existence of an atlas of complex charts $(U, \psi)$, such that all the transition functions are holomorphic. Indeed the previous definition will only be satisfied under certain quite restrictive conditions.

**Theorem** (**Newlander-Nirenberg** [20]) The almost complex structure $J$ gives each tangent space $T_p M$ the structure of a complex vector space. A necessary and sufficient condition for there to exist a holomorphic chart around each point of $M$, is the vanishing of the Nijenhuis tensor $N_J$ of $J$, where we write $N_J(v, w) = [v, w] + J([Jv, w] + [v, Jw]) - [Jv, Jw]$ for all smooth vector fields $v, w$ on $M$.

The Nijenhuis tensor $N_J$ represents an observation to the existence of holomorphic functions on $M$. The equations that a function $f$ must satisfy in order to be holomorphic on $M$ , and given in previous definition, are called the Cauchy-Riemann equations. For $n = 1$, the manifold always admits holomorphic coordinates, but for $n > 1$ the Cauchy-Riemann equations are over determined, and $N_J$ is an obstruction to the existence of holomorphic functions on $M$. If $N_J \neq 0$, then there can still exist some holomorphic functions on $M$, but not enough to construct holomorphic coordinates. Thus, an almost complex manifold only admits many, i.e., enough, holomorphic functions if the Nijenhuis tensor $N_J$ vanishes.

The Newlander –Nirenberg Theorem implies a new definition of complex manifold as a *2n*-dimensional manifold with an almost complex structure $J$ such that $N_J \equiv 0$. See [19] for more details

## de Raham Cohomology

Definition: Let $C^\infty(\Lambda^k T^*M \otimes_R C)$ denote the space of smooth complex differential forms of degree k on $M$, then we can write the exterior differential $d$ as

$$d: C^\infty(\Lambda^k T^*M \otimes_R C) \to C^\infty(\Lambda^{k+1} T^*M \otimes_R C)$$

Where $T^*M$ is the cotangent bundle, $\Lambda^k T^*M \otimes_R C$ is the complex vector bundle of complex valued $k-$form over $M$ and $C^\infty(\Lambda^k T^*M \otimes_R C)$ is the space of smooth sections of $\Lambda^k T^*M \otimes_R C$. The exterior differential $d$ satisfies the Leibniz rule and $= 0$ .

Since $dod = 0$, the chain of operators

$$0 \xrightarrow{d} C^\infty(\Lambda^0 T^*M \otimes_R C) \xrightarrow{d} C^\infty(\Lambda^1 T^*M \otimes_R C) \xrightarrow{d} ... \xrightarrow{d} C^\infty(\Lambda^n T^*M \otimes_R C) \xrightarrow{d} 0 \quad (3.2)$$

forms a complex.

Definition: The kernel of $d$ are the closed forms and the image of $d$ are exact forms.

**Definition** (**de Rham cohomology Group**) For $k = 0,1,...,n$ we define the $k$-th de Rham cohomology group of $M$ by

$$H^k_{dR}(M,C) = \frac{\ker\big(d: C^\infty(\Lambda^k T^*M \otimes_R C) \to C^\infty(\Lambda^{k+1} T^*M \otimes_R C)\big)}{Im(d: C^\infty(\Lambda^{k-1} T^*M \otimes_R C) \to C^\infty(\Lambda^k T^*M \otimes_R C))}$$

Note that the de Rham cohomology of a smooth manifold  , $H^k_{dR}(M,C)$ is isomorphic to the cohomology $H^k(M,C)$ of $M$ as topological space.

Definition: The Hodge star $*$ is an isomorphism of vector bundles $*: \Lambda^k T^*M \to \Lambda^{n-k} T^*M$, which is defined as follows. For $\alpha$ and $\beta$ k-forms , $* \beta$ is the unique $(n-k)-$form that satisfies the equation $\alpha \wedge (* \beta) = (\alpha, \beta) dV_g$ for all k-forms $\alpha$ on $M$ (here $M$ is a manifold with Riemannian metric $g$, and $dV_g$ is Volume form on $M$).

Let $(M, g)$ be a Riemannian manifold, then we can define the corresponding formal adjoint (we can introduce an $L^2$-metric on $\Lambda^k T^*M \otimes_R C$ ) differential operator $d^*$ where

$$d^*: C^\infty(\Lambda^{k+1} T^*M \otimes_R C) \to C^\infty(\Lambda^k T^*M \otimes_R C)$$

by

$$d^* \alpha = (-1)^{kn+n+1} * d(* \alpha)$$

Also note that if $\alpha$ is a $k-$form and $\beta$ be a $(k+1)-$form, then $\langle d\alpha, \beta \rangle_{L^2} = \langle \alpha, d^* \beta \rangle_{L^2}$ thus $d^*$ has the formal properties of the adjoint of . As $d^2 = 0$ we find that $(d^*)^2 = 0$, so that the corresponding chain of operator form a complex , similar to the expression given in (3.2) , of which we can compute the cohomology group as

$$H^k_{dR} = \frac{\ker\big(d^*: C^\infty(\Lambda^k T^*M \otimes_R C) \to C^\infty(\Lambda^{k-1} T^*M \otimes_R C)\big)}{Im(d^*: C^\infty(\Lambda^{k+1} T^*M \otimes_R C) \to C^\infty(\Lambda^k T^*M \otimes_R C))}$$

**Theorem** : For $(M, g)$ a complex Riemannian manifold, if we define the Laplacian $\Delta_d = dd^* + d^*d$ then $\Delta \alpha = 0 \Leftrightarrow d\alpha = 0$ and $d^*\alpha = 0$

## Dolbeault cohomology

Given an almost complex structure $J$ on a manifold $M$, the decomposition of the complex tangent bundle as $TM \otimes_R C = T^{1,0}M \oplus T^{0,1}M$ induces a similar decomposition on the bundle of complex differential forms ;

$$\left(\Lambda^k T^*M\right) \otimes_R C = \bigoplus_{p+q=k} \Lambda^{p,q} M$$

where $\Lambda^{p,q} M$ is the bundle $\Lambda^p T^{*1,0}M \otimes_C \Lambda^q T^{*0,1}M$ .

A section of $\Lambda^{p,q} M$ is called a $(p,q)$ −form which is a complex-valued differential form which can be expressed in local holomorphic coordinates $(z_1, z_2, \ldots, z_n)$ as

$$\sum_{\substack{1 \le i_1 < i_2 < \cdots < i_p \le n \\ 1 \le j_1 < j_2 < \cdots < j_p \le n}} f_{i_1 \ldots i_p, j_1 \ldots j_q} dz_{i_1} \wedge \ldots \wedge dz_{i_p} \wedge d\overline{z}_{j_1} \wedge \ldots \wedge d\overline{z}_{j_q}$$

The exterior differential $d$ splits informally as

$$d\alpha^{p,q} = \left(N_J.\,\alpha\right)^{p+2,q-1} + \partial\alpha^{p+1,q} + \overline{\partial}\alpha^{p,q+1} + \left(N_J.\,\alpha\right)^{p+2,q-1} \qquad (3.3)$$

where $\partial$ and $\overline{\partial}$ are operators such that

$$\partial\colon C^\infty(\Lambda^{p,q}M) \to C^\infty(\Lambda^{p+1,q}M)$$

$$\overline{\partial}\colon C^\infty(\Lambda^{p,q}M) \to C^\infty(\Lambda^{p,q+1}M)$$

And $N_J$ is the Nijenhuis tensor. So if the almost complex structure $J$ defined on the manifold $M$ is integerable, i.e., $M$ is a complex manifold and $N_J \equiv 0$, then we can rewrite (3.3) as

$$d\alpha^{p,q} = \partial\alpha^{p+1,q} + \overline{\partial}\alpha^{p,q+1}$$

Which we rewrite more easily as

$$d = \partial + \overline{\partial}$$

Also by a bit computation if $M$ is a complex manifold, then

$\partial o \partial = 0, \overline{\partial} o \overline{\partial} = 0$ and $\overline{\partial} + \overline{\partial} o \partial = 0$ . Also since $\overline{\partial}^2 = 0$ for each $p = 0, \ldots, n$ the chain of operators

$$0 \xrightarrow{\overline{\partial}} C^\infty(\Lambda^{p,0}M) \xrightarrow{\overline{\partial}} C^\infty(\Lambda^{p,1}M) \xrightarrow{\overline{\partial}} \ldots \xrightarrow{\overline{\partial}} C^\infty(\Lambda^{p,n}M)$$

forms a complex and for $p, q = 0, \ldots n$ we define Dolbeault Cohomology groups of $M$ by

$$H_{\overline{\partial}}^{p,q}(M;C) = \frac{\ker\left(\overline{\partial}\colon C^\infty(\Lambda^{p,q}M) \to C^\infty(\Lambda^{p,q+1}M)\right)}{Im\left(\overline{\partial}\colon C^\infty(\Lambda^{p-1,q}M) \to C^\infty(\Lambda^{p,q}M)\right)}$$

Also if the complex manifold $(M, J)$ carries a Hermitian metric $g$, then the coadjoint operator to $d$ given by $d^*$ splits similarly as $d^* = \partial^* + \overline{\partial^*}$  where

$$\partial^*: C^\infty(\Lambda^{p+1,q}M) \to C^\infty(\Lambda^{p,q}M)$$

$$\overline{\partial^*}: C^\infty(\Lambda^{p,q+1}M) \to C^\infty(\Lambda^{p,q}M)$$

The corresponding Laplacian operators are

$$\Delta_\partial = \partial \partial^* + \partial^* \partial$$

$$\Delta_{\overline{\partial}} = \overline{\partial}\,\overline{\partial^*} + \overline{\partial^*}\,\overline{\partial}$$

Also we have a fact that for $(M, J, \omega)$ a Kähler manifold the Laplacian splits as

$$\Delta_d = \Delta_\partial + \Delta_{\overline{\partial}}$$

and we have this fact that $\Delta_\partial = \Delta_{\overline{\partial}} = \frac{1}{2}\Delta_d$ .

**Proposition**: (The Global $\partial\overline{\partial}$-Lemma [19]) : If $M$ is a complex manifold and $d\alpha = 0$, then the following are equivalent

1. $\alpha$ is $d$ −exact
2. $\alpha$ is $\partial$ −exact
3. $\alpha$ is $\overline{\partial}$ −exact
4. $\alpha$ is $\partial\overline{\partial}$ −exact

- ## Chern Classes

In differential geometry the Chern classes are characteristic classes associated to a complex vector bundle . They are topological invariant associated to vector bundle on a smooth manifold. The question of whether two ostensibly different vector bundles are the same can be quite hard to answer. The Chern classes provide a simple test. If the chern classes of a pair of vector bundles do not agree, then the vector bundles are different. The converse however is not true.

The following proposition can be taken as a definition

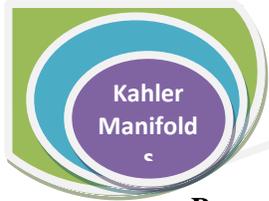

**Proposition[19]**: To every complex vector bundle $E$ over a smooth manifold $M$ one can associate a cohomology class $C_1(E) \in H^2(M, Z)$ called the first Chern classes of $E$ satisfying the following axioms:

- *(Naturality)* For every smooth $f: M \to N$ and complex vector bundle $E$ over $N$, one has $f^*(C_1(E)) = C_1(f^*E)$, where the left term denotes the pull-back in cohomology and $f^*E$ is the pull-back bundle defined by $f^*E_x = E_{f(x)}$, $\forall x \in M$.

- (**Whitney sum formula**) For every bundles $E, F$ over $M$ one has $C_1(E \oplus F) = C_1(E) + C_1(F)$, where $E \oplus F$ is the Whitney sum defined as the pull-back of the bundle $E \times F \to M \times M$ by the diagonal inclusion of $M$ in $M \times M$

- (**Normalization**) The first Chern class of the tautological bundle of $CP^1$ is equal to $-1$ in $H^2(CP^1, Z) \cong Z$, which means that the integral over $CP^1$ of any representative of this class equals -1.

**Remark**: Given any Kahler metric $\omega$, which its Ricci form $Ric(\omega)$ represents the first Chern class of $(M; J)$ (in the cohomology space $H^2(M; Z)$), if $M$ is closed, then the cohomology space could be

$$H^{1,1}(M; \mathbb{C}) \cap H^2(M; \mathbb{Z})$$

**Theorem (Gauss-Bonnet[19])**: The real Chern class of a Vector bundle is the image of the integer Chern class $C_1(B, Z)$ ($B$ is a line bundle) under the natural homomorphism $H^2(M, Z) \to H^2(M, R)$

**Definition**: The first Chern class of a complex manifold is $C_1(\Lambda^{n,0}(M))$

**Definition**: A Calabi-Yau manifold is a complex Kähler manifold with $C_1(M, Z) = 0$.

S.S. Chern pointed out that Chern classes can be represented by curvature form since both of them measure how far the vector bundle is away from a trivial product structure. Now, we will give the definition of Chern classes by curvature form by the language of de Rham cohomology theory.

**Definition: (Invariant Polynomial)**. Let $A = \left(A_{ij}\right)_{r \times r}$ be a matrix, a homogeneous polynomial defined by $P(A) = P(A_{11}, \dots, A_{ij}, \dots, A_{rr})$ is called invariant if

$$P(A) = P(gAg^{-1})$$

for any $g \in GL(r, \mathbb{C})$.

In the following context we consider the invariant polynomial $P(A) = \sigma_k(A)$ where $\sigma_k$ is the coefficient of $t^k$ in the polynomial

$$\det(I + tA) = \sum_{k=0}^{n} \sigma_k(A) t^k$$

$\sigma_k(A)$ for $k = 1, \dots, n$ is an invariant polynomial.

Now we present a theorem of Madson [],

**Theorem:** The invariant polynomial $P(\Theta)$ is a closed form and hence defines a cohomology class $[P(\Theta)]$ in the de Rham cohomology group of the vector bundle. Moreover the cohomology class $[P(\Theta)]$ is independent of choices of connections.

**Definition** (**Chern Forms, Chern Classes**). Define the Chern forms $c_i(\Theta)$ of the curvature $\Theta$ by

$$c_i(\Theta) = \sigma_i\left(\frac{\sqrt{-1}}{2\pi}\Theta\right) \qquad , \qquad i = 0, \dots, r$$

Suppose $E$ is a vector bundle over a complex manifold M . Define the $i-th$ Chern classes of $E$ by

$$c_i(E) = \sigma_i\left(\frac{\sqrt{-1}}{2\pi}\Theta\right) \in H_{DR}^{2i}(E) \qquad , \qquad i = 0, \dots, r$$

For $M$ a complex manifold of dimension $n$, define the $i-th$ Chern classes $c_i(M)$ of $M$ to be the *i-th* Chern classes of its tangent bundle $TM$ for $i = 1, \dots, n$.

**First Chern Classes and second Chern classes.**

By previous definition, the Chern classes are

$$c_1(E) = \left[\sigma_1\left(\frac{\sqrt{-1}}{2\pi}\Theta\right)\right] \qquad , \qquad c_2(E) = \left[\sigma_2\left(\frac{\sqrt{-1}}{2\pi}\Theta\right)\right]$$

The Chern forms are

$$c_1(\Theta) = \sigma_1\left(\frac{\sqrt{-1}}{2\pi}\Theta\right) = \left(\frac{\sqrt{-1}}{2\pi}\right)\text{Trace}(\Theta) = \left(\frac{\sqrt{-1}}{2\pi}\right)\left(\sum_{\alpha}\Theta_\alpha^\alpha\right)$$

and

$$c_2(\Theta) = \frac{-1}{4\pi^2} \sigma_2(\Theta) = \frac{-1}{4\pi^2} \sum_{|I|=2} \det(\Theta_{I,I}) = \frac{-1}{4\pi^2} \cdot \frac{1}{2} \cdot \sum_{\alpha,\beta} \Theta_\alpha^\alpha \wedge \Theta_\beta^\beta - \Theta_\alpha^\beta \wedge \Theta_\beta^\alpha,$$

Where $|I|$ stands for the number of indices in $I$.

Suppose that $M$ is a Hermitian manifold and E is the tangent bundle $TM$, then the first and second Chern forms are given by:

$$c_1(\Theta) = \frac{\sqrt{-1}}{2\pi} \sum_i \Theta_i^i = \frac{\sqrt{-1}}{2\pi} \sum_i \sum_{k,l} \Theta_{i\,kl}^i \, dz^k \wedge d\bar{z}^l$$

$$c_2(\Theta) = \frac{-1}{4\pi^2} \cdot \frac{1}{2} \sum_{i,j} \Theta_i^i \wedge \Theta_j^j - \Theta_i^j \wedge \Theta_j^i$$

$$= \frac{-1}{8\pi^2} \sum_{i,j} \sum_{p,q,r,s} \left( \Theta_{i\,p\bar{q}}^i \Theta_{j\,r\bar{s}}^j - \Theta_{i\,p\bar{q}}^j \Theta_{j\,r\bar{s}}^i \right) dz^p \wedge d\bar{z}^q \wedge dz^r \wedge d\bar{z}^s$$

So, we have the first Chern form of $M$ is the Ricci form of the tangent bundle of $M$:

$$c_1(\Theta) = \frac{\sqrt{-1}}{2\pi} \sum_{k,l} Ric_{k\bar{l}} \, dz^k \wedge d\bar{z}^l$$

That is, the Ricci form is in the cohomology class $c_1(M)$. If $M$ is further Einstein, then the Kahler form which equals to the Ricci form up to a constant, will also in the class $c_1(M)$.

**Definition**: We say $c_1(M) > 0$ , $(< 0)$ if $c_1(M)$ can be represented by a positive (negative) form. In local coordinates, this means that

$$\phi = \sqrt{-1}\phi_{k\bar{l}} dz_k \wedge d\bar{z}_l$$

where $\phi_{k\bar{l}}$ is positive (negative) definite. We say $c_1(M) = 0$, if the first Chern class $c_1(M)$ is cohomologous to zero. (here $c_1(M)$ represented by $\phi$ )

**Definition**: We say that $g$ is a Kahler-Einstein metric if there exists a real constant $\lambda$ such that $Ric(g) = \lambda\omega_g$ . A Kahler manifold $(M, g)$ is Kahler-Einstein if $g$ is a Kahler-Einstein metric.

Note that, $= \frac{s}{2m}$, where $m$ is the complex dimension and $s$, is the scalar curvature. It can be easily proven that the value of $s$, is as follows

$$V.s^m = \frac{(4\pi m)^m}{m!} c_1^m.$$

where, $c_1^m$ denotes the Chern number associated with the $m$-th power of the first Chern class of $M$, depending on the complex structure only, and $V$, total volume.

**Example**: In this example we compute the first Chern class of $\mathbb{C}P^n$. We know

$$\omega = \frac{\sqrt{-1}}{2} \partial \bar{\partial} \log(1 + \sum |z_i|^2) > 0,$$

So we can write

$$\omega^n = \left(\frac{\sqrt{-1}}{2}\right)^n \frac{n! \, dz_1 \wedge \ldots \wedge d\bar{z}_n}{(1 + \sum |z_i|^2)^{n+1}}$$

This implies that

$$Ric(\omega) = \frac{n+1}{\pi} \omega \quad \text{and so} \quad c_1(M) = \frac{n+1}{\pi} [\omega]$$

Where $\frac{1}{\pi}[\omega]$ is the positive generator of the cohomology $H^2(\mathbb{C}P^n, \mathbb{Z}) = \mathbb{Z}$. Also by some manipulations one can show that

$$\frac{2(n+1)}{n} c_2(M) = c_1(M)^2$$

here $c_1(M)^2 = c_1(M) \wedge c_1(M)$.

# Chapter 4

# Kahler -Ricci flow

Kähler geometry is the meeting place of Riemannian geometry, complex geometry and symplectic geometry, and the study of Kähler manifolds brings together the techniques developed in a wide range of mathematical disciplines including partial differential equation , notably the complex Monge- Ampere equations, algebraic geometry and Cohomology theory. In this chapter, we focuses on the study of the Ricci flow on Kähler manifolds, an important class of manifolds in complex differential geometry.

The Kähler Ricci flow is simply an abbreviation for the Ricci flow on Kähler manifolds. For more details see [10-11, 19 and 34]

**Definition**: Let $(M, J)$ be an almost complex manifold. A vector field $V$ is an infinitesimal automorphism of the almost complex structure if the Lie derivative of $J$ with respect to $V$ is zero, i.e.,

$$\mathcal{L}_V J = 0 \qquad\qquad (4.1)$$

Note that (4.1) is equivalent to $J([V, W]) = [V, JW]$ for any vector $W$, because

$$(\mathcal{L}_V J)(W) = \mathcal{L}_V (JW) - J(\mathcal{L}_V W) = [V, JW] - J([V, W])$$

There are various equivalent ways to define a Kähler manifold .

**Definition**: We say that a Riemannian manifold $(M, g)$ with an almost complex structure $J: TM \to TM$ is a Kähler manifold if the metric $g$ is $J-$invariant (or sometimes we say Hermitian):

$$g(JX, JY) = g(X, Y)$$

and $J$ is parallel :

$$\nabla J = 0$$

or equivalently , $\nabla_X(JY) = J\nabla_X(Y)$ for all $X, Y$ and the metric $g$ is called a Kähler metric.

**Lemma:** Almost complex structure which yield Kähler manifolds are necessarily integrable

Proof: By applying the Newlander-Nierenberg theorem, we only need to check that the Nijenhuis tensor vanishes for a Kähler manifold:

$$N_J(X, Y) = \nabla_{JX}JY - \nabla_{JY}JX - J(\nabla_{JX}Y - \nabla_Y JX) - J(\nabla_X JY - \nabla_{JY}X) - \nabla_X Y + \nabla_Y X$$

$$= -[J(\nabla_{JX}Y) - \nabla_{JX}JY] + [J(\nabla_{JY}X) - \nabla_{JY}JX] + [J(\nabla_Y JX) - \nabla_Y J(JX)]$$

$$- [J(\nabla_X JY) - \nabla_X J(JY)] = 0$$

So the proof is complete.

Let $\{z^\alpha\}$ be local holomorphic coordinates. We may write $z^\alpha := x^\alpha + iy^\alpha$, where $x^\alpha$ and $y^\alpha$ are real valued functions. Define $dz^\alpha := dx^\alpha + idy^\alpha$, $d\overline{z}^\alpha := dx^\alpha - idy^\alpha$ and

$$\frac{\partial}{\partial z^\alpha} := \frac{1}{2}\left(\frac{\partial}{\partial x^\alpha} - i\frac{\partial}{\partial y^\alpha}\right)$$

$$\frac{\partial}{\partial \overline{z}^\alpha} := \frac{1}{2}\left(\frac{\partial}{\partial x^\alpha} + i\frac{\partial}{\partial y^\alpha}\right)$$

So that

$$dz^\alpha\left(\frac{\partial}{\partial z^\beta}\right) = \delta_\beta^\alpha, \quad d\overline{z}^\alpha\left(\frac{\partial}{\partial \overline{z}^\beta}\right) = \delta_\beta^\alpha \quad , \quad d\overline{z}^\alpha\left(\frac{\partial}{\partial z^\beta}\right) = dz^\alpha\left(\frac{\partial}{\partial \overline{z}^\beta}\right) = 0$$

Note that

$$\overline{\frac{\partial}{\partial z^\beta}} = \frac{\partial}{\partial \overline{z}^\beta}, \quad \overline{\frac{\partial}{\partial \overline{z}^\beta}} = \frac{\partial}{\partial z^\beta}$$

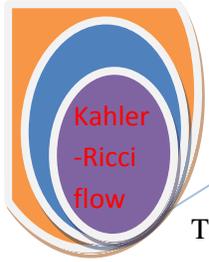

The holomorphic tangent bundle $T^{1,0}M$ is locally the span of the vectors $\left\{\frac{\partial}{\partial z^\alpha}\right\}_{\alpha=1}^n$ and also the tangent bundle $T^{0,1}M$ is locally the span of $\left\{\frac{\partial}{\partial \overline{z}^\beta}\right\}_{\beta=1}^n$.

We extend the Riemannian metric $g$ complex linearly to define

$$g_C : T_C M_p \times T_C M_p \to C$$

Similarly

$$\nabla_C : T_C M \times C^\infty(T_C M) \to C^\infty(T_C M)$$

is the complex linear extension of $\nabla : TM \times C^\infty(TM) \to C^\infty(TM)$ with the convention that $\nabla_X Y := \nabla(X,Y)$ and $(\nabla_C)_X Y := \nabla_C(X,Y)$.

The complex linear extension of $Rm$ is denoted by $Rm_C$.

Let $g_{\alpha\overline{\beta}} := g_C\left(\frac{\partial}{\partial z^\alpha}, \frac{\partial}{\partial \overline{z}^\beta}\right)$, $g_{\overline{\beta}\alpha} := g_C\left(\frac{\partial}{\partial \overline{z}^\beta}, \frac{\partial}{\partial z^\alpha}\right)$, since $\overline{g_C(V,W)} = g_C(\overline{V}, \overline{W}) = g_C(\overline{W}, \overline{V})$, these coefficients satisfy the Hermition condition :

$$g_{\alpha\overline{\beta}} = \overline{g_{\beta\overline{\alpha}}} = g_{\overline{\beta}\alpha}$$

Similarly, we define

$$g_{\alpha\beta} := g_C\left(\frac{\partial}{\partial z^\alpha}, \frac{\partial}{\partial z^\beta}\right), \quad g_{\overline{\alpha}\overline{\beta}} := g_C\left(\frac{\partial}{\partial \overline{z}^\alpha}, \frac{\partial}{\partial \overline{z}^\beta}\right)$$

**Remark:** We have

$$g_{\alpha\beta} = g_{\overline{\alpha}\overline{\beta}} = 0$$

Proof: If $X, Y \in T^{1,0}M$, then

$$g_C(X,Y) = g_C(JX, JY) = g_C\left(\sqrt{-1}X, \sqrt{-1}Y\right) = -g_C(X,Y)$$

Which implies $g_C(X,Y) = 0$, Similarly, if $X, Y \in T^{0,1}M$, then we also have $g_C(X,Y) = 0$. Thus, in local holomorphic coordinates, the Kähler metric takes the form

$$g_C = g_{\alpha\overline{\beta}}\left(dz^\alpha \otimes d\overline{z}^\beta + d\overline{z}^\beta \otimes dz^\alpha\right)$$

**Definition:** We say that a $(p,p)-$tensor (or form) $\eta$ is real if $\overline{\eta} = \eta$.

**Definition :** The Kähler form $\omega$, on a Riemannian manifold $(M, J, g)$ with an almost complex structure and whose metric is Hermitian, is defined to be the 2-form associated to $g$:

$$\omega(X, Y) := g(JX, Y)$$

which is a real (1,1)-form.

In holomorphic coordinate chart $(U, (z^1, \dots, z^n))$ with $p$ centered at origin, the Kähler form is

$$\omega = \sqrt{-1} \sum_{i,j} g_{i\overline{j}} \, dz^i \wedge d\overline{z}^j$$

**Theorem**: If $g$ is a Kähler metric, then $\omega$ is a closed 2-form. In fact, $\omega$ is parallel. (The converse is also true)

Proof: We compute

$$(\nabla_X \omega)(Y, Z) = X(\omega(Y, Z)) - \omega(\nabla_X Y, Z) - \omega(Y, \nabla_X Z) = X\big(g(JY, Z)\big) - g(J(\nabla_X Y), Z) - g(JY, \nabla_X Z)$$

$$= (\nabla_X g)(JY, Z) = 0$$

Where here, we used the definition of $\omega$, $J(\nabla_X Y) = \nabla_X (JY)$ , and $g$ is parallel. So proof is complete.

Note that in our convention

$$dz^\alpha \wedge d\overline{z}^\beta \left( \frac{\partial}{\partial z^\gamma}, \frac{\partial}{\partial \overline{z}^\delta} \right) = \frac{1}{2} \delta^\alpha_\gamma \delta^\beta_\delta.$$

Remark: We have $d\omega = 0 \Leftrightarrow \frac{\partial}{\partial z^\gamma} g_{\alpha\overline{\beta}} = \frac{\partial}{\partial z^\alpha} g_{\gamma\overline{\beta}}$ .

**Definition**: The Real cohomology class $[\omega] \in H^{1,1}(M. R) \subseteq H^2(M. R)$, is called the Kähler class of $\omega$.

Definition: The christoffel symbols of the Levi-Civita connection , defined by

$$(\nabla_C)_{\frac{\partial}{\partial z^\alpha}} \frac{\partial}{\partial z^\beta} := \sum_{\gamma=1}^{n} \left( \Gamma^\gamma_{\alpha\beta} \frac{\partial}{\partial z^\gamma} + \Gamma^{\overline{\gamma}}_{\alpha\beta} \frac{\partial}{\partial \overline{z}^\gamma} \right)$$

$$(\nabla_C)_{\frac{\partial}{\partial z^\alpha}} \frac{\partial}{\partial \overline{z}^\beta} := \sum_{\gamma=1}^{n} \left( \Gamma^\gamma_{\alpha\overline{\beta}} \frac{\partial}{\partial z^\gamma} + \Gamma^{\overline{\gamma}}_{\alpha\overline{\beta}} \frac{\partial}{\partial \overline{z}^\gamma} \right)$$

etc. they are zero unless all the indices are unbarred or all the indices are barred .

Let $g^{\alpha\overline{\beta}}$ be defined by $g_{\gamma\overline{\beta}}g^{\alpha\overline{\beta}} = \delta_\gamma^\alpha$. Then in holomorphic coordinates, we have

$$\Gamma_{\alpha\beta}^\gamma = \frac{1}{2}g^{\gamma\overline{\delta}}\left(\frac{\partial}{\partial z^\alpha}g_{\beta\overline{\delta}} + \frac{\partial}{\partial z^\beta}g_{\alpha\overline{\delta}} - \frac{\partial}{\partial\overline{z}^\delta}g_{\alpha\beta}\right) = g^{\gamma\overline{\delta}}\frac{\partial}{\partial z^\alpha}g_{\beta\overline{\delta}}$$

and

$$\Gamma_{\alpha\beta}^\gamma = \Gamma_{\beta\alpha}^\gamma$$

By same method

$$\Gamma_{\alpha\overline{\beta}}^\gamma = \frac{1}{2}g^{\gamma\overline{\delta}}\left(\frac{\partial}{\partial z^\alpha}g_{\overline{\beta}\overline{\delta}} + \frac{\partial}{\partial\overline{z}^\beta}g_{\alpha\overline{\delta}} - \frac{\partial}{\partial\overline{z}^\delta}g_{\alpha\overline{\beta}}\right) = 0$$

and also we have

$$\Gamma_{\overline{\alpha}\overline{\beta}}^\gamma = 0$$

**Definition**: The components of the curvature $(3,1)-$tensor $Rm_C$ are defined by

$$Rm_C\left(\frac{\partial}{\partial z^\alpha}, \frac{\partial}{\partial\overline{z}^\beta}\right)\frac{\partial}{\partial z^\gamma} := R_{\alpha\overline{\beta}\gamma}^\delta\frac{\partial}{\partial z^\delta} + R_{\alpha\overline{\beta}\gamma}^{\overline{\delta}}\frac{\partial}{\partial z^{\overline{\delta}}}$$

$$Rm_C\left(\frac{\partial}{\partial z^\alpha}, \frac{\partial}{\partial\overline{z}^\beta}\right)\frac{\partial}{\partial z^{\overline{\gamma}}} := R_{\alpha\overline{\beta}\overline{\gamma}}^\delta\frac{\partial}{\partial z^\delta} + R_{\alpha\overline{\beta}\overline{\gamma}}^{\overline{\delta}}\frac{\partial}{\partial z^{\overline{\delta}}}$$

etc,. Also we can write

$$R_{\alpha\overline{\beta}\gamma\overline{\delta}} = Rm_C\left(\frac{\partial}{\partial z^\alpha}, \frac{\partial}{\partial\overline{z}^\beta}, \frac{\partial}{\partial z^\gamma}, \frac{\partial}{\partial\overline{z}^\delta}\right)$$

Which is satisfies in following relations

$$R_{\alpha\overline{\beta}\gamma\overline{\delta}} = g_{\eta\overline{\delta}}R_{\alpha\overline{\beta}\gamma}^\eta, \quad R_{\alpha\overline{\beta}\overline{\gamma}\delta} = g_{\delta\overline{\eta}}R_{\alpha\overline{\beta}\overline{\gamma}}^{\overline{\eta}}$$

Note that the only non-vanishing components of the $(3,1)-$tensor are $R_{\alpha\overline{\beta}\gamma}^\delta$, $R_{\alpha\overline{\beta}\overline{\gamma}}^{\overline{\delta}}$, $R_{\overline{\alpha}\beta\gamma}^\delta$, and $R_{\overline{\alpha}\beta\overline{\gamma}}^{\overline{\delta}}$.

So the only non-vanishing components of the curvature (4,0)-tensor are

$$R_{\alpha\overline{\beta}\gamma\overline{\delta}}, \ R_{\alpha\overline{\beta}\overline{\gamma}\delta}, \ R_{\overline{\alpha}\beta\gamma\overline{\delta}}, R_{\overline{\alpha}\beta\overline{\gamma}\delta}$$

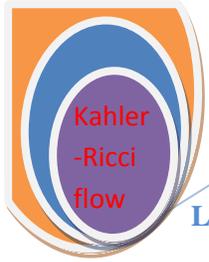

**Lemma**: For a Kähler Manifold $Rm(X,Y)$ is $J-$invariant

Proof : We have to show that $Rm(X,Y)JZ = J(Rm(X,Y)Z)$ so we have

$$Rm(X,Y)JZ =$$

$$= \nabla_X \nabla_Y JZ - \nabla_Y \nabla_X JZ - \nabla_{[X,Y]} JZ$$

$$= \nabla_X (J\nabla_Y Z) - \nabla_Y (J\nabla_X Z) - J(\nabla_{[X,Y]} Z)$$

$$= J(\nabla_X \nabla_Y Z) - J(\nabla_Y \nabla_X Z) - J(\nabla_{[X,Y]} Z)$$

$$= J(Rm(X,Y)Z)$$

**Definition**: The components of $Rc_C$ (or $Ric_C$) are defined by

$$Rc_C \left( \frac{\partial}{\partial z^\alpha}, \frac{\partial}{\partial \overline{z}^\beta} \right) := R_{\alpha\overline{\beta}}, \ \ Rc_C \left( \frac{\partial}{\partial z^\alpha}, \frac{\partial}{\partial z^\beta} \right) := R_{\alpha\beta}$$

Note that $R_{\alpha\overline{\beta}} = R^\delta_{\delta\overline{\beta}\alpha}$.

**Definition**: The scalar curvature of a metric $\omega$ is defined to be

$$s(\omega) = R_{i\overline{i}j\overline{j}} = -g^{iJ} \frac{\partial^2}{\partial z_i \partial \overline{z}_i} \log \det g_{k\overline{l}}$$

so the trace of Ricci curvature is scalar curvature.

Now by using first Chern class we give some facts on Kahler Einstein manifolds.

**Lemma**. We have the following relation between scalar curvature and first Chern class.

$$\pi c_1(M) [\omega_g]^{n-1} = \frac{1}{n} \int_M s(\omega_g) \omega_g^n$$

Proof: We have,

$$\omega_g = \frac{\sqrt{-1}}{2} \sum dz_i \wedge d\overline{z}_i \ \ , \ \ Ric(g) = \frac{\sqrt{-1}}{2} \sum R_{i\overline{i}} \, dz_i \wedge d\overline{z}_i$$

so it follows that,

$$Ric(g) \wedge \omega_g^{n-1} = (n-1)! \left(\frac{\sqrt{-1}}{2}\right)^n \sum_i R_{i\bar{i}} dz_1 \wedge d\bar{z}_1 \wedge ... \wedge dz_n \wedge d\bar{z}_n = \frac{1}{n} s(\omega_g) \omega_g^n$$

So we get

$$\pi c_1(M) [\omega_g]^{n-1} = \frac{1}{n} \int_M s(\omega_g) \omega_g^n$$

so we get the desired result and proof is complete.

Previous lemma learn to us that the average of scalar curvature is depends only on $[\omega_g]$, and $c_1(M)$.

So we are ready to state following theorem. But at first we give a proposition which is very useful for our manipulations in next chapter.

**Proposition**: Let $(M, g)$ be a Kahler Manifold and let $\phi_1, \phi_2 \in H^{1,1}(M, \mathbb{C})$ and suppose that $\phi_1$ is cohomologous to $\phi_2$. Then there exists a function $f \in C^\infty(M, \mathbb{R})$ such that $\phi_1 - \phi_2 = \partial\bar{\partial}f$.

**Theorem.** If $\pi c_1(M) = \lambda[\omega]$ and $\omega$ is a Kahler metric with constant scalar curvature, then $\omega$ is Kahler-Einstein. (Here we assume $M$ is compact).

**Proof**. By definition of first Chern class, we know $\frac{1}{\pi} Ric(\omega)$ represents the first Chern class of $M$ and that it is a (1,1)- form. So by previous lemma we have that $Ric(\omega) - \lambda\omega = \partial\bar{\partial}f$. But we know the Trace of Ricci curvature is scalar curvature. So by taking the trace on both side of the equation we get $s(\omega) - n\lambda = \Delta f$. (Note that the trace of $\partial\bar{\partial}f$ is $g^{i\bar{j}} \frac{\partial^2}{\partial z_i \partial \bar{z}_j}$ which is exactly $\Delta f$). But we showed that the scalar curvature is only depends on $[\omega_g]$, and $c_1(M)$, so, $s(\omega) = n\lambda$, so $f$ is a harmonic function on a compact manifold and therefore $f$ is constant and we get the desired result.

**Theorem.** For a Kähler Manifold $(M, g, J)$, we have that

$$R_{i\bar{j}} = -\partial_{\bar{j}} \partial_i \log \det g \quad (4.2)$$

Proof: We use of following identity for Christoffel symbols

$$\Gamma_{jk}^l = (\partial_j g_{k\bar{m}}) g^{\bar{m}l}$$

So by applying well-known formula for the derivative of the determinant of Hermitian matrix, (let $A = \left(A_{i\bar{j}}\right)$ with inverse $\left(A^{\bar{j}i}\right)$, then $\frac{\partial}{\partial s} \det A = A^{\bar{j}i} \left(\frac{d}{ds} A_{i\bar{j}}\right) \det A$). So by using this fact we calculate

$$R_{i\bar{j}} = -\partial_{\bar{j}}\Gamma_{ki}^{k} = -\partial_{\bar{j}}\left(g^{\bar{q}k}\partial_{i}g_{k\bar{q}}\right) = -\partial_{\bar{j}}\partial_{i}\log\det g$$

So we obtain the desired result and proof is complete.

Also we can show that for a Kähler Manifold $Rc$ is $J$-invariant, i.e.,

$$Rc(JX, JY) = Rc(X, Y)$$

**Definition:** The Ricci form $\rho$ is the 2-form associated to $Rc$.

$$\rho(X, Y) = \frac{1}{2}Rc(JX, Y)$$

which is a real $(1,1)$ −form. Also form $\rho$ is closed 2-form and in holomorphic coordinates we can write the Ricci form $\rho$ as follows

$$\rho = \sqrt{-1}R_{\alpha\bar{\beta}}dz^{\alpha}\wedge d\bar{z}^{\beta}.$$

And we may express $\rho$ as $\rho = -\sqrt{-1}\partial\bar{\partial}\log\det\left(g_{\gamma\bar{\delta}}\right)$.

**Definition**: (**First Chern classes**): The Real de Raham cohomology class $\left[\frac{1}{2\pi}\rho\right] := C_1(M)$ is the first chern class of $M$ (also it is only depends on the complex structure of $M$ )

**Definition**: The holomorphic sectional and bisectional curvature defined by

$$K_C(Z) = \frac{Rm_C(Z, \bar{Z}, Z, \bar{Z})}{|Z|^4}, \quad Z \in T^{1,0}M$$

$$K_C(Z, W) = \frac{Rm_C\left(Z, \bar{Z}, W, \bar{W}\right)}{|Z|^2|W|^2}, \quad Z, W \in T^{1,0}M - \{0\}$$

Respectively, Clearly we have $K_C(Z, Z) = K_C(Z)$. We say that bisectional curvature is positive if $K_C(Z, W) > 0$. Also bisectional curvature can be interpreted as follows:

If $|x| = |y| = 1$, and set

$$u = \frac{1}{\sqrt{2}}\left(x - \sqrt{-1}Jx\right) \text{ and } v = \frac{1}{\sqrt{2}}\left(y - \sqrt{-1}Jy\right),$$

Then

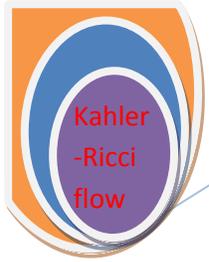

$$R(u, \bar{u}, v, \bar{v}) = R(x, y, y, x) + R(x, Jy, Jy, x)$$

Also a Kahler manifold $(M, g)$, is said to be of constant bisectional curvature if there exists a constant $\lambda$ such that in any local coordinates of $M$,

$$R_{i\bar{j}k\bar{l}} = \lambda(g_{i\bar{j}}g_{k\bar{l}} + g_{i\bar{l}}g_{k\bar{j}})$$

**Theorem:** The Kahler manifold $\mathbb{C}P^n$ is a manifold of constant bisectional curvature $+1$.

Proof: Note that on $\mathbb{C}P^n$, as we computed in previous chapter, we may write

$$\omega_g^n = \left(\frac{\sqrt{-1}}{2}\right)^2 \frac{(dz_i \wedge d\bar{z}_i)^n}{(1 + |z|^2)^{n+1}}$$

therefore by applying equation (4.2) we get,

$$Ric_{i\bar{j}} = -\frac{\partial^2}{\partial z_i \partial \bar{z}_j} \log\left(\frac{1}{(1 + |z|^2)^{n+1}}\right) = (n+1)g_{i\bar{j}}$$

But before we showed $g_{i\bar{j}} = \frac{\partial^2 \log(1+|z|^2)}{\partial z_i \partial \bar{z}_j}$ , so

$$-\frac{\partial^2 g_{i\bar{j}}}{\partial z_k \partial \bar{z}_l}\bigg|_{z=0} = -\frac{\partial^4 \log(1 + |z|^2)}{\partial z_k \partial \bar{z}_l \partial z_i \partial \bar{z}_j}\bigg|_{z=0} = \frac{1}{2}\frac{\partial^4 |z|^4}{\partial z_k \partial \bar{z}_l \partial z_i \partial \bar{z}_j}\bigg|_{z=0} = \frac{\partial^3}{\partial z_k \partial \bar{z}_l \partial z_i}\left(|z|^2 z_j\right)\bigg|_{z=0}$$

$$= \frac{\partial^2}{\partial z_k \partial z_i}(z_j z_l)\bigg|_{z=0} = \delta_{ij}\delta_{kl} + \delta_{il}\delta_{kj}$$

but, because isometry group of $g$ i.e, $(SU(n + 1))$ acts on $\mathbb{C}P^n$ transitively, so we can write $R_{i\bar{j}k\bar{l}} = g_{i\bar{j}}g_{k\bar{l}} + g_{i\bar{l}}g_{k\bar{j}}$ . Therefore, by definition, $\mathbb{C}P^n$ is a manifold of constant bisectional curvature $+1$.

**Remark**: The manifold $\mathbb{C}^n$ is the flat metric and the bisectional curvature vanishes. Moreover for Kahler Manifold $B^n = \{z \in \mathbb{C}^n; |z| < 1\}$ and let

$$\omega_g = \frac{\sqrt{-1}}{2}\partial\bar{\partial}\log(1 - |z|^2)$$

then $R_{i\bar{j}k\bar{l}} = -(g_{i\bar{j}}g_{k\bar{l}} + g_{i\bar{l}}g_{k\bar{j}})$ and $(B^n, g)$ is Kahler manifold of constant bisectional curvature $-1$.

So we are ready to give the Uniformization Theorem as follows.

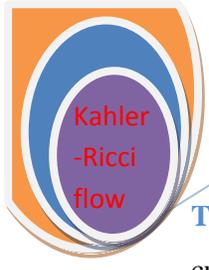

**Theorem (Uniformization Theorem):** If $(M, g)$ is a complete Kahler manifold of constant bisectional curvature $R_{i\bar{j}k\bar{l}} = \lambda(g_{i\bar{j}}g_{k\bar{l}} + g_{i\bar{l}}g_{k\bar{j}})$ for some constant $\lambda$, then its universal covering $\tilde{M}$ is $\mathbb{C}^n$, $\mathbb{C}P^n$ or $B^n$. Moreover, up to scaling, $g$ pulls back to one of the metrics of $\mathbb{C}^n$, $\mathbb{C}P^n$, or $B^n$.

Now we are ready to present an important fact in Kahler-Einstein Manifolds in following theorem.

**Theorem:** Let $(M, g)$ be a Kahler- Einstein manifold. Then the universal covering $(\tilde{M}, g) \cong$ $\mathbb{C}P^n, \mathbb{C}^n \ or \ B^n$ (here $\cong$ mean isometric up to scaling) if and only if

$$\left(\frac{2(n+1)}{n} c_2(M) - c_1(M)^2\right) [\omega]^{n-2} = 0$$

**Proof:** By applying Uniformization Theorem, it suffices to show that there exists a constant $\lambda$ such that $R_{i\bar{j}k\bar{l}} = \lambda(g_{i\bar{j}}g_{k\bar{l}} + g_{i\bar{l}}g_{k\bar{j}})$. But from the definition of Chern class we have

$$\det\left(I + \frac{t\sqrt{-1}}{2\pi}\Omega\right) = I + t\phi_1(\Omega) + t^2\phi_2(\Omega) + \cdots,$$

where $\Omega$ is a matrix valued 2-form $\left(\Omega_i^{\ j}\right)$, which is actually of type $(1,1)$, and can be represented by

$$\Omega_i^{\ j} = g^{j\bar{p}}R_{i\bar{p}k\bar{l}}dz_k \wedge d\bar{z}_l$$

And moreover $c_i(M)$ is represented by $\phi_i(g)$ and $[\phi_i] \in H^{i,i}(M, \mathbb{C}) \cap H^{2i}(M, \mathbb{R})$. In particular $\phi_1(\Omega) = tr\Omega$ represent the first Chern class (the reason is: in general case we have $\det(I + \varepsilon X) = 1 + tr(X)\varepsilon + \mathcal{O}(\varepsilon 2)$). By using properties of trace and determinant,

$$\frac{1}{4\pi^2} tr(\Omega \wedge \Omega) = \frac{1}{4\pi^2}\sum_{k,i}\Omega_i^k \wedge \Omega_i^k = \frac{(\sqrt{-1})^2}{4\pi^2}R_{ip\bar{q}}^k R_{kr\bar{s}}^i dz^p \wedge d\bar{z}^q \wedge dz^r \wedge d\bar{z}^s$$

represents $c_1^2(M) - 2c_2(M)$. So $\phi_1^2(M) - 2\phi_2(M) = \frac{1}{4\pi^2} tr(\Omega \wedge \Omega)$. So we get

$$(c_1^2(M) - 2c_2(M))[\omega]^{n-2} = \left(\frac{\sqrt{-1}}{2\pi}\right)^2 \int_M tr(\Omega \wedge \Omega) \wedge \omega^{n-2}$$

Denote by $R^0(g)$ the traceless part of the curvature $R(g)$. Because $M$ is Kahler-Einstein, we have $(g) = \lambda\omega_g$ . and in local coordinates we get

$$R_{i\bar{j}k\bar{l}}^0 = R_{i\bar{j}k\bar{l}} - \frac{\lambda}{n+1}(g_{i\bar{j}}g_{k\bar{l}} + g_{i\bar{l}}g_{k\bar{j}})$$

In reality the tensor $R^0$ measures how much the metric deviates from having constant bisectional curvature. So if we denote by $|.|$ the norm given by the metric, then its norm square easily computed as

$$|R^0|^2 = |R|^2 + \frac{\lambda^2}{(n+1)^2}(2n^2 + 2n) - \frac{4\lambda}{n+1}R$$

The assumption $R_{i\bar{j}} = \lambda g_{i\bar{j}}$ gives $R = \lambda n$ and $|Ric|^2 = \lambda^2 n$. Then

$$|R^0|^2 = |R|^2 - \frac{2\lambda^2 n}{n+1}$$

so firstly, we show the following equality

$$\frac{1}{n(n-1)4\pi^2}\int_M |R^0|^2\, \omega^n = \left(2c_2(M) - \left(1 - \frac{1}{n+1}\right)c_1^2(M)\right).[\omega]^{n-2}$$

we may weite

$$n(n-1)tr(\Omega \wedge \Omega) \wedge \omega^{n-2} =$$

$$= \sum_{p\neq r}\left(R^k_{ip\bar{p}}R^i_{kr\bar{r}} - R^k_{ip\bar{r}}R^i_{kr\bar{p}}\right)\omega^n$$

$$= \sum_{p,r}\left(R^k_{ip\bar{p}}R^i_{kr\bar{r}} - R^k_{ip\bar{r}}R^i_{kr\bar{p}}\right)\omega^n$$

$$= (|Ric|^2 - |R|^2)\omega^n$$

$$= (\lambda^2 n - |R|^2)\omega^n$$

therefore we get

$$|R^0|^2 \frac{\omega^n}{n(n-1)} = -tr(\Omega \wedge \Omega) \wedge \omega^{n-2} + \lambda^2\left(\frac{1}{n-1} - \frac{2}{(n+1)(n-1)}\right) = -tr(\Omega \wedge \Omega) \wedge \omega^{n-2} + \frac{\lambda^2}{n+1}$$

now notice that, from previous chapter we know

$$\frac{1}{4\pi^2}\int_M \lambda^2\omega^n = \frac{1}{4\pi^2}\int_M (\lambda\omega)^2 \wedge \omega^{n-2} = c_1^2(M).[\omega]^{n-2}$$

and so

$$\frac{1}{n(n-1)4\pi^2}\int_M |R^0|^2\,\omega^n = \left(2c_2(M) - \left(1 - \frac{1}{n+1}\right)c_1^2(M)\right).[\omega]^{n-2}$$

so from this equality we have shown that $M$ is a manifold of constant bisectional curvature(so $R^0 = 0$ because if $R^0 \neq 0$ so $M$ would be biholomorphic to unit ball in $\mathbb{C}^n$) if and only if

$$\left(\frac{2(n+1)}{n}c_2(M) - c_1(M)^2\right)[\omega]^{n-2} = 0$$

so we get the desired result and proof is complte.

**Definition** : The Laplacian acting on tensors is given by

$$\Delta = \frac{1}{2}g^{\alpha\overline{\beta}}\left(\nabla_\alpha\nabla_{\overline{\beta}} + \nabla_{\overline{\beta}}\nabla_\alpha\right)$$

and the Laplacian when acts on functions is $\Delta = g^{\alpha\overline{\beta}}\nabla_\alpha\nabla_{\overline{\beta}} = g^{\alpha\overline{\beta}}\frac{\partial^2}{\partial z^\alpha \partial\overline{z}^\beta}$.

- ➢ **Kähler Ricci flow and monge-Ampere equation**

  The Kahler Ricci flow has tremendous developments in 2000s and becames a rapidly growing topic in geometry. The Kähler –Ricci flow is a second order nonlinear parabolic flow of Kähler metrics. Here we introduce Kähler Ricci flows and show that the Kähler Ricci flow equation is equivalent to parabolic complex monge-Ampere equation which describes our main results.

**Definition**: The Kähler-Ricci flow equation defined by

$$\frac{\partial}{\partial t}g_{\alpha\overline{\beta}}(x,t) = -R_{\alpha\overline{\beta}}(x,t)$$

Also for compact Kähler manifold $(M, J, \omega_0)$ we can define Kähler-Ricci flow equation by $(1,1)$ −form $\omega_0$ as follows.

$$\frac{\partial\omega_t}{\partial t} = -Ric(\omega_t), \qquad\qquad \omega_t|_{t=0} = \omega_0 \qquad (4.3)$$

**Remark**: The first Chern class $c_1(M) = [Ric(\omega_0)]$ depends only on $J$ which is fixed under the flow. Also the Kähler class $[\omega_t]$ under Kähler-Ricci flow $(4.3)$ satisfies the following ODE:

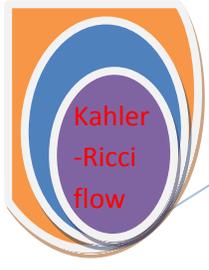

$$\frac{d}{dt}[\omega_t] = -c_1(M) \Longrightarrow [\omega_t] = [\omega_0] - tc_1(M)$$

Firstly, we show that the Kahler Ricci flow equation is equivalent to parabolic complex monge-Ampere equation. First, let $Ric(\omega)$ and $\Omega$ represent the first Chern class and so by $\partial\overline{\partial}$ −lemma tells us that we can find $F$, only depending on $\omega$ and $\Omega$, such that

$$\Omega - \text{Ric}(\omega) = \sqrt{-1}\partial\overline{\partial}F,$$

where $F$ is unique after normalizing to

$$\int_M (e^F - 1)\omega^n = 0$$

It follows from the $\partial\overline{\partial}$ −lemma, any Kahler metric $\widetilde{\omega}$ cohomologous to $\omega$ has the form $\omega + \sqrt{-1}\partial\overline{\partial}\phi$. Suppose function $\phi$ satisfies

$$Ric\big(\omega + \sqrt{-1}\partial\overline{\partial}\phi\big) = \Omega = Ric(\omega) - \sqrt{-1}\partial\overline{\partial}F \qquad (4.4)$$

So, by applying local exression of Ricci curvature in local coordinate , $R_{i\overline{j}} = -\frac{\partial^2}{\partial z^i \partial \overline{z}^j}\log\det g$ , we can write the expression (4.4) as follows,

$$-\sqrt{-1}\partial\overline{\partial}\log\det\left(g_{i\overline{j}} + \frac{\partial^2\phi}{\partial z_i \partial \overline{z}_j}\right) = -\sqrt{-1}\partial\overline{\partial}\log\det(g_{i\overline{j}}) - \sqrt{-1}\ \partial\overline{\partial}F$$

Although this is only locally defined, the following is globally defined

$$\sqrt{-1}\partial\overline{\partial}\log\left(\frac{\det\left(g_{i\overline{j}} + \phi_{i\overline{j}}\right)}{\det(g_{i\overline{j}})}\right) = \sqrt{-1}\ \partial\overline{\partial}F$$

Therefore,

$$\frac{\det\left(g_{i\overline{j}} + \phi_{i\overline{j}}\right)}{\det(g_{i\overline{j}})} = e^F \qquad (4.5)$$

where $F$ is a smooth function on $M$ with $\int_M (e^F - 1)\omega^n = 0$. Equation (4.5) is just a complex Monge-Ampere equation which is also equivalent to

$$\big(\omega + \sqrt{-1}\partial\overline{\partial}\phi\big)^n = e^F\ \omega^n.$$

In general, the exact maximal time $T$ of a Riemannian Ricci flow may not be easy to find. However, fortunately, for Kähler-Ricci flows, the maximal time of existence $T$ is explicitly determined by the initial Kähler class $[\omega_0]$ and the first Chern class.

From (4.3), we know $[\omega_t] = [\omega_0] - t c_1(M)$, so let $\rho \in c_1(M)$, then $\widetilde{\omega}_t = \omega_0 - t\rho$ is in the same Kähler class as the flow metric $\omega_t$. Therefore by $\partial \overline{\partial}$ −lemma, there exists a family of smooth functions $\varphi_t$ on $M$ such that $\omega_t = \widetilde{\omega}_t + \sqrt{-1}\partial\overline{\partial}\varphi_t$. Let $\Omega$ be a time-independent volume form on $M$ such that $\sqrt{-1}\partial\overline{\partial}\Omega = -\rho$. By the same way of previous part, the Kähler Ricci flow equation (4.3) can be rewritten as the following parabolic complex Monge-Ampere equation for $\varphi_t$ as follows

$$\frac{\partial \varphi_t}{\partial t} = \log \frac{\left(\widetilde{\omega}_t + \sqrt{-1}\partial\overline{\partial}\varphi_t\right)^n}{\Omega} \qquad \varphi_0 = 0.$$

Now we state a theorem from Tian and Zhang [], for maximal of existence of $T$.

**Theorem** (Tian-Zhang, []). Let $(M, \omega(t))$ be an Kahler-Ricci flow $\frac{\partial \omega_t}{\partial t} = -Ric(\omega_t)$ on a compact Kahler manifold $M$, with $dim_{\mathbb{C}} = n$, then the maximal existence time $T$ is given by

$$T = \sup\{t : [\omega_0] - t c_1(M) > 0\}$$

## ➤ The normalized Kähler – Ricci flow

Suppose that $(M^n, J, g_0)$ be a closed manifold. We make basic assumption that the first Chern class is a multiple of the Kähler class

$$[\rho_0] = c[\omega_0] \text{ for some } c \in R$$

(for more details see [34]) So it is very important that first Chern class has a sign . i.e., is negative definite, zero, or positive definite

It is very possible that first Chern classes, or indeed any cohomology classes on which we have a concept of order, are neither positive, negative or zero. For example, consider a projective line $P^1$ and a curve $C$ of genus 2 (or greater).The line has a positive first Chern class, while the curve has a negative one. Their product is a compact surface which has neither positive nor negative first Chern class.

We know

$$Vol(M) = Vol_g(M) = \frac{1}{n!}\int_M \omega^n$$

and

$$\int_M R d\mu = \frac{1}{(n-1)!}\int_M \rho \wedge \omega^{n-1}$$

We find that $C = \frac{r}{n}$, where $r := \frac{\int_M R_0 d\mu_{g_0}}{Vol_{g_0}(M)}$ is the average (complex) scalar curvature, so

$$\frac{r}{2\pi n}[\omega_0] = \frac{1}{2\pi}[\rho_0] = C_1(M)$$

So we can define the normalized Kähler Ricci flow as follows

$$\frac{\partial}{\partial t}g_{\alpha\overline{\beta}} = -R_{\alpha\overline{\beta}} + \frac{r}{n}g_{\alpha\overline{\beta}}$$

for $t \in [0, T)$.

**Definition: (Potential function)** We define the potential function $f = f(t)$ by following equation

$$R_{\alpha\overline{\beta}}(g) - \frac{r}{n}g_{\alpha\overline{\beta}} = \partial_\alpha\partial_{\overline{\beta}}f = \nabla_\alpha\nabla_{\overline{\beta}}f$$

**Theorem**([34]): Under the normalized Kähler –Ricci flow we have

1)  $\frac{\partial}{\partial t}d\mu = (r - R)d\mu$

2)  $\frac{\partial R}{\partial t} = \Delta R + \left|R_{\alpha\overline{\beta}}\right|^2 - \frac{r}{n}R$

3)  The normalized Kähler Ricci flow preserves the Volume.

Proof: By applying the definition of normalized Kähler Ricci flow we have

$$\frac{\partial}{\partial t}\log\det g_{\gamma\overline{\delta}} = g^{\gamma\overline{\delta}}\frac{\partial}{\partial t}g_{\gamma\overline{\delta}} = r - R$$

hence by using definition of $d\mu$ we get

$$\frac{\partial}{\partial t} d\mu = (r - R) d\mu$$

For 2) the evalution of the Ricci tensor is

$$\frac{\partial}{\partial t} R_{\alpha\overline{\beta}} = -\partial_\alpha \partial_{\overline{\beta}} \left( \frac{\partial}{\partial t} \log \det g_{\gamma\overline{\delta}} \right) = \partial_\alpha \partial_{\overline{\beta}} R - \frac{\partial}{\partial t} g_{\alpha\overline{\beta}} . R_{\alpha\overline{\beta}}$$

So from the definition of normalized Kähler Ricci flow equation we get the desired result.

For 3) since

$$\int_M (r - R) d\mu = 0$$

So the normalized Kähler Ricci flow preserves the Volume and proof is complete.

# Chapter 5

# Calabi-Yau Conjecture

Our main goal in this section is to present a complete proof of the Calabi conjecture. In the middle of the 70s, Calabi conjecture was solved by S. T. Yau in case the first Chern class is vanishing and Aubin and Yau, independently, in case the Chern class is negative. The uniqueness in these two cases was done by E. Calabi himself in the 50s. In reality finding a distinguished canonical metric on a smooth manifold is one of the fundamental problems in the theory of geometric analysis and here canonical means that the metric depends on the complex structure and is unique up to biholomorphic automorphism. H. Poincaré's Uniformization theorem settles this problem for Riemann surfaces and Calabi conjecture is generalization of Poincaré's Uniformization theorem in higher dimension. In fact Poincaré's Uniformization theorem says that, there is a unique metric with constant curvature in each Kahler class on a Riemann surface. After Calabi tried to extend this fundamental problem to a compact Kahler manifold and he conjectured existence of Kähler-Einstein metrics on a compact Kahler manifold with its first Chern class definite. S.T. Yau, reformulated and reduced Calabi conjecture to a so-called Monge-Ampere equation. He reduced the proof of the existence of solution of Monge-Ampere equation:

$$\det\left(g_{i\bar{j}} + u_{i\bar{j}}\right) = \det\left(g_{i\bar{j}}\right) f(z), \qquad (*)$$

(where $f$ is a smooth positive function on $(M, g)$ ) to a priori estimate by applying Schauder theory and continuity method. This work by Yau opened a vast field for the study of complex Monge-Ampere type

(*). Moreover complex Monge-Ampere is a powerful tool in understanding geometry and topology in Kahler setting. Also Yau studied the generalized form of this equation when the right hand side function $f(z)$ may degenerate or have poles. Later Tian and Yau, solved equation (*) on complete non-compact Kahler manifolds and obtained some strong applications to algebraic geometry. We can now formulate the famous conjecture which was posed by Eugenio Calabi in 1954 and eventually proved by Shing-Tung Yau, [14].

**The Calabi conjecture:**

Let $M$ be a compact, Complex manifold, and $g$ a Kähler metric on $M$ with Kähler form $\omega$, Then for each real, closed (1,1)-form $\rho'$ on $M$ such that $[\rho'] = 2\pi C_1(M)$ in $H^2(M, R)$ there exists a unique Kähler metric $g'$ on $M$ with Kähler form $\omega'$, such that $[\omega] = [\omega']$ in $H^2(M, R)$ and the Ricci form of $g'$ is $\rho'$. In the specific case $C_1(M) = 0$, we can take $\rho' = 0$.

The Calabi conjecture closely related to the question of which Kahler manifolds have Kahler-Eninstein metrics. So firstly, the main propose of this chapter is proving following problem with complete details which is famous to second type of calabi conjecture as follows

**Second type of Calabi conjecture:**

If a compact Kahler manifold has a negative, zero, or positive first Chern class, then it has a kahler-Einstein metric in the same class as its kahler metric, unique, up to rescaling. This was proved for negative first chern classes independently by Thierry Aubin and S.T. Yau. Moreover when the Chern class is zero it was proved by Yau as an easy consequence of the Calabi conjecture.

The Calabi conjecture can be reduced to solving a complex Monge- Ampere equation as follows.

**Another version of the Calabi Conjecture:**

Let $M$ be a compact, complex manifold of the dimension $n$, and $g$ a Kähler metric on $M$, with the Kähler form $\omega$, Let $f$ be a smooth real function on $M$, and define $A > 0$, by $A \int_M e^f dV_g = Vol_g(M)$. Then there exists a unique smooth real function $\varphi$ such that

(a) $\omega + dd^c\varphi$ is a positive (1,1)-form

(b) $\int_M \varphi dV_g = 0$

(c) $(\omega + dd^c\varphi)^n = Ae^f \omega^n$

Moreover, in local holomorphic coordinates $z_1, \ldots, z_n$ condition (c) can be expressed in the following way :

$$\det\left(g_{\alpha\overline{\beta}} + \frac{\partial^2\varphi}{\partial z_\alpha \partial \overline{z_\beta}}\right) = Ae^f \det(g_{\alpha\overline{\beta}})$$

and Yau proved such an equation can be solved by using the continuity method. For proving the Yau-Calabi conjecture we can concentrate on the nonlinear Elliptic and parabolic PDE aspects of the subject. Yau proves the existence of a geometric structure using PDE's, giving importance to the idea that deep insights into geometry can be obtained by studying solutions of such equations. One of main objects in Kähler Geometry is finding Kähler Einstein metrics on a compact manifold. Recall that a necessary condition for a Kähler manifold to admit such metrics is that the first Chern class $C_1(M)$ has a definite sign. In this chapter the idea of Cao will be discussed in detail. By flowing any Kähler metric on a compact Kähler manifold with either $C_1(M) = 0$ or $C_1(M) \leq 0$ , by the Ricci flow, we obtain the unique Kähler –Einstein metric in the same Kähler class as the starting metric.

Richard Hamilton in [1], proved that for any compact 3-manifold with positive Ricci curvature one can deform the initial metric along the heat flow defined by the equation

$$\frac{\partial g_{ij}}{\partial t} = -2R_{ij} + \frac{2}{3}rg_{ij} \qquad (5.1)$$

to an Einstein metric of positive scalar curvature.

We consider the complex version of Hamilton's equation (5.1), of the following type

$$\begin{cases} \frac{\partial \tilde{g}_{i\overline{j}}}{\partial t} = -\tilde{R}_{i\overline{j}} + T_{i\overline{j}} \\ \tilde{g}_{i\overline{j}}(0) = g_{i\overline{j}} \end{cases} \qquad (5.2)$$

where $\tilde{R}_{i\overline{j}}$ denotes the Ricci tensor of the metric $\tilde{g}_{i\overline{j}}$ and after we will prove that $T_{i\overline{j}} = \tilde{R}_{i\overline{j}}(\infty)$.

At first we prove that the solution for equation (5.2), exists for all time, and converges to a limiting metric $\tilde{g}_{i\overline{j}}(\infty)$, as well as show that the derivatives $\frac{\partial \tilde{g}_{i\overline{j}}}{\partial t}(t)$ converges uniformly to a constant as $t$ approaches infinity. Then $\tilde{g}(\infty)$ is the Kähler- Einstein metric which we want.

The beauty of Kähler geometry is that sometimes, all our study relies on the Kähler potential.

Let $\Omega = \frac{\sqrt{-1}}{2\pi} T_{i\bar{j}}\, dz^i \wedge d\bar{z}^j$, be a fixed representative of the first Chern class $C_1(M)$ and denote as before , the Ricci form of a Kähler metric $g$ by

$$\text{Ric} = \frac{\sqrt{-1}}{2\pi} R_{i\bar{j}}\, dz^i \wedge d\bar{z}^j$$

**Lemma** 5.1($\partial\bar{\partial}$ −Lemma [34]) For a global closed $(p,q)$ −form $\alpha$ trivial in cohomology, there is a global $(p-1, q-1)$ −form $\eta$ such that $= \sqrt{-1}\partial\bar{\partial}\eta$ .

Since both $\Omega$ and Ricci form lie in the same cohomology class $C_1(M)$, the $\partial\bar{\partial}$ −Lemma tells us that there exists a smooth function $f$ on $M$ such that

$$T_{i\bar{j}} - R_{i\bar{j}} = \frac{\partial^2 f}{\partial z^i \partial \bar{z}^j}$$

Therefore , if we let (we will see in main theorem that why we assumed this)

$$\tilde{g}_{i\bar{j}} = g_{i\bar{j}} + \frac{\partial^2 u}{\partial z^i \partial \bar{z}^j}$$

Where $u(t)$ is defined on $M \times [0, T)$ with $u(0) = 0$, $0 < T \leq \infty$.

So from equation (5.2)

$$\frac{\partial^2}{\partial z^i \partial \bar{z}^j}\left(\frac{\partial u}{\partial t}\right) = -\tilde{R}_{i\bar{j}} + R_{i\bar{j}} + \frac{\partial^2 f}{\partial z^i \partial \bar{z}^j}$$

and consequently, according to well-known  relation

$$R_{i\bar{j}} = -\frac{\partial^2}{\partial z^i \partial \bar{z}^j} \log \det(g_{i\bar{j}})$$

We get

$$\frac{\partial^2}{\partial z^i \partial \bar{z}^j}\left(\frac{\partial u}{\partial t}\right) = \frac{\partial^2}{\partial z^i \partial \bar{z}^j}\left(\log det\left(g_{i\bar{j}} + \frac{\partial^2 u}{\partial z^i \partial \bar{z}^j}\right) - \log det\left(g_{i\bar{j}}\right)\right) + \frac{\partial^2 f}{\partial z^i \partial \bar{z}^j}$$

or equivalently

$$\partial\bar{\partial}\left(\frac{\partial u}{\partial t}\right) = \partial\bar{\partial} \log det\left(g_{i\bar{j}} + \frac{\partial^2 u}{\partial z^i \partial \bar{z}^j}\right) - \partial\bar{\partial} \log \det\left(g_{i\bar{j}}\right) + \partial\bar{\partial} f$$

**Lemma**([**22-24**]): If $(M, J, g)$ is compact Kähler manifolds and $f: M \to R$ be of $C^2$ −class such that $\partial \overline{\partial} f = 0$. Then $f$ is constant .

Simplifying even more, from previous lemma the evaluation equation for the Kähler potential $u$ is

$$\frac{\partial u}{\partial t} = \log \det \left( g_{i\overline{j}} + \frac{\partial^2 u}{\partial z^i \partial \overline{z}^j} \right) - \log \det \left( g_{i\overline{j}} \right) + f + \varphi(t) \qquad (5.3)$$

where $\varphi(t)$ is a smooth function in $t$ and from (5.3) it satisfying the compatibility condition

$$\int exp \left( \frac{\partial u}{\partial t} - f \right) dV = exp \big( \varphi(t) \big) Vol(M)$$

where $dV$ is the volume element of the metric $g_{i\overline{j}}$

Now since (5.3), is a nonlinear parabolic equation we identify the maximal existence time for a smooth solution of the Kähler Ricci flow.

# ✛ Maximal existence time for a smooth solution of the Kähler Ricci flow

Let $\omega = \omega(t)$ be a solution of the Kähler Ricci flow

$$\begin{cases} \frac{\partial}{\partial t} \omega = -Ric\,(\omega) \\ \omega(0) = \omega_0 \end{cases} \qquad (5.4)$$

As long as the solution exists, the cohomology class $[\omega(t)]$ evolves by

$$\begin{cases} \frac{d}{dt} [\omega(t)] = -C_1(M) \\ [\omega(0)] = [\omega_0] \end{cases}$$

and solving the ordinary differential equation gives

$$[\omega(t)] = [\omega_0] - tC_1(M).$$

Immediately we see that a necessary condition for the Kähler- Ricci flow to exist for $t \in [0, \acute{t})$ is that $[\omega_0] - tC_1(M) > 0$ .

Actually this necessary condition is sufficient. If we define

$$T = Sup \{t > 0 \mid [\omega_0] - tC_1(M) > 0\}$$

Then we have following theorem which was proved by Tian . For finding another proof of the following theorem see [39] .

**Theorem**: There exists a unique maximal solution $g(t)$ of the Kähler Ricci flow for for $t \in [0, T)$.

Proof: without loss of generality fix $\acute{T} < T$. We will show that there exists a solution to (5.4) on $[0, \acute{T})$. To do this, we need to choose metrics $\widehat{\omega}_t$ in the cohomology class $[\omega_0] - tC_1(M)$. Since $[\omega_0] - \acute{T}C_1(M)$ is a Kähler class, there exists a Kähler form $\eta$ in $[\omega_0] - \acute{T}C_1(M)$. We choose our family of metrics $\widehat{\omega}_t$ to be the linear path of metrics between $\omega_0$ and $\eta$. Namely, define

$$\chi = \frac{1}{\acute{T}}(\eta - \omega_0) \in -C_1(M),$$

and

$$\widehat{\omega}_t = \omega_0 + t\chi = \frac{1}{\acute{T}}\left((\acute{T} - t)\omega_0 + t\eta\right) \in [\omega_0] - tC_1(M)$$

Fix a volume form $\Omega$ on $M$ with

$$\frac{\sqrt{-1}}{2\pi} \partial\overline{\partial} \log \Omega = \chi = \frac{\partial}{\partial t} \widehat{\omega}_t \in -C_1(M)$$

We now consider the parabolic complex monge-Amere equation, for $\varphi = \varphi(t)$ a real valued function on $M$.

$$\frac{\partial\varphi}{\partial t} = \log \frac{\left(\widehat{\omega}_t + \frac{\sqrt{-1}}{2\pi}\partial\overline{\partial}\varphi\right)^n}{\Omega}, \quad \widehat{\omega}_t + \frac{\sqrt{-1}}{2\pi}\partial\overline{\partial}\varphi > 0 \quad , \varphi(0) = 0 \quad (5.5)$$

This equation is equivalent to the Kähler Ricci-flow (5.4), Indeed, given a smooth solution $\varphi$ of (5.5) on $[0, \acute{T})$, we can obtain a solution $\omega = \omega(t)$ of (5.4) on $[0, \acute{T})$ as follows. Define $\omega(t) = \widehat{\omega}_t + \frac{\sqrt{-1}}{2\pi}\partial\overline{\partial}\varphi$ and observe that $\omega(0) = \widehat{\omega}_0 = \omega_0$ and

$$\frac{\partial}{\partial t}\omega = \frac{\partial}{\partial t}\widehat{\omega}_t + \frac{\sqrt{-1}}{2\pi}\partial\overline{\partial}\left(\frac{\partial\varphi}{\partial t}\right) = -Ric(\omega)$$

as required. Conversely, suppose that $\omega = \omega(t)$ solves (5.4) on $[0,\acute{T})$. Then since $\widehat{\omega}_t \in [\omega(t)]$, we we can apply the $\partial\overline{\partial}$ −Lemma to find a smooth family of potential functions $\widehat{\varphi}(t)$ such that $\omega(t) = \widehat{\omega}_t + \frac{\sqrt{-1}}{2\pi}\partial\overline{\partial}\widehat{\varphi}(t)$ then

$$\frac{\sqrt{-1}}{2\pi}\partial\overline{\partial}\log\omega^n = \frac{\partial}{\partial t}\omega = \frac{\sqrt{-1}}{2\pi}\partial\overline{\partial}\log\Omega + \frac{\sqrt{-1}}{2\pi}\partial\overline{\partial}\left(\frac{\partial\widehat{\varphi}}{\partial t}\right)$$

So $\frac{\partial\widehat{\varphi}}{\partial t} = \log\frac{\omega^n}{\Omega} + c(t)$ for some smooth function $c:[0,\acute{T})\to R$. Now set $\varphi(t) = \widehat{\varphi}(t) - \int_0^t c(s)ds - \widehat{\varphi}(0)$, noting that since $\omega(0) = \omega_0$ the function $\widehat{\varphi}(0)$ is constant . It follows that $\varphi = \varphi(t)$ solves (5.5), so it suffices only we study (5.5). Since the linearization of the right hand side of (5.5) is the Laplace operator $\Delta_{g(t)}$ which is strictly elliptic (because the linearization of the Monge-Ampere operator is the Laplacian which is strictly elliptic), it follows that (5.5) is a strictly parabolic nonlinear partial differential equation for . So It then follows from "standard parabolic theory" for nonlinear parabolic equations that there is a unique maximal solution on time interval $[0,T_{Max})$.also we can show that $T_{Max} < T'$ (for more details see[34])But if we use of more details only following theorem is sufficient for studying the equation (5.5).

**Theorem** (S.T.Yau[26 and 14]) Suppose $F \in C^k(M)$ ($k \geq 3$). is positive and $\Omega$ is a smooth volume form on $M$. Then

$$\left(\omega + \frac{\sqrt{-1}}{2\pi}\partial\overline{\partial}u\right)^n = \exp\{F + cu\}\,\Omega, \quad c \geq 0 \text{ is a constant} \quad (\beth)$$

has a unique solution $u \in C^{k+1,\alpha}(M)$ for any $\alpha \in [0,1)$.

Proof: At first we prove solutions of this equation are unique up to constants. To see this suppose that $u,v \in C^2(M)$ and $\omega_u{}^n = \omega_v{}^n$ then

$$0 = \omega_u{}^n - \omega_v{}^n = \int_0^1 \frac{d}{dt}(t\omega_u + (1-t)\omega_v)^n dt$$

$$= \left\{n\int_0^1(t\omega_u + (1-t)\omega_v)^{n-1}dt\right\}\wedge\frac{\sqrt{-1}}{2\pi}\partial\overline{\partial}(u-v)$$

(*)

Since $\omega_u$ and $\omega_v$ are known and positive definite, $(t\omega_u + (1-t)\omega_v)^{n-1}$ is positive as a $(n-1,n-1)$ −form $\forall\, t \in [0,1]$. Furthermore, (*) is a linear elliptic equation. By the maximum principle for linear elliptic equations, $sup_M u - v = sup_{\partial M} u - v$, otherwise it achieves its maximum in the interior which implies $u - v$ is constant. $\partial M = \emptyset$ so $u - v$ is constant. If $u$ and $v$ are mean zero, then the integral of $u - v$ over $M$ is zero, so necessarily $u = v$ and we proved uniqueness.

To prove the existence of a solution, a priori estimates up to $C^{2,\alpha}$ are needed but we will prove this fact after and use of this fact freely here. With the $C^2$-estimate for solutions (⊐) , in hand we can now describe the continuity approach used to prove the existence of a solution. The method of continuity associates a family of Monge-Ampere equations,

$$\left(\omega + \frac{\sqrt{-1}}{2\pi}\partial\bar{\partial}u_t\right)^n = A(t)\exp\{tF + cu\}\,\omega^n, \quad \text{for c} \geq 0 \text{ and t} \in [0,1] \quad (*_t)$$

to equation (⊐), Where $A(t) = Vol(M).\left(\int_M e^{tF+cu}\omega^n\right)^{-1}$ is a compatibility constant for each $t$, so it is necessary that integrals over $M$ of both sides of $(*_t)$ be equal for each t $\in [0,1]$ . Note the equation at $t = 1$ is (⊐), the equation we want to solve. By showing $S = \{$t $\in [0,1] : (*_t)$has a solution $\}$ is nonempty, open, and closed the existence of a solution is estanlished. $S \neq \emptyset$ because $u_0 \in R$ a constant solves the equation at $t = 0$ , and for proving why $S$ is open and closed see [35 in chapter Calabi-Yau theorem] .

## Short time existence on semi-linear parabolic PDE

We present an affirmative reason for short time existence on non-linear parabolic PDE. In reality this is kind of meta-theorem and has several version. For instance the fully non-linear case $\frac{\partial u}{\partial t} = \det(\nabla^2 u)$ with $u(t = 0,.)$ convex, it is very complicate and only for quasi-linear case like, $\frac{\partial u}{\partial t} = \text{div}(|\nabla|^{p-1}\nabla u)$, has been treated by several scientific papers.

But proving short time existence on the semi-parabolic linear case, specially $\frac{\partial u}{\partial t} - \Delta u$, which is our main focus is not too hard. So in the case

$$\frac{\partial u}{\partial t} + Lu = f(u, \nabla u)$$

where $L$ is an elliptic linear operator, like Laplacian $-\Delta$, has a rather simple philosophy, So firstly, we look at the Cauchy-Lipschitz-Picard theorem.

**Proposition :** Let $\mathcal{H}$ is a Hilbert space, $u_0 \in \mathcal{H}$, $f \in C^0([0, +\infty[, \mathcal{H})$ and function $A : \mathcal{H} \to \mathcal{H}$, be linear continuous and $\|Au - Av\|_{\mathcal{H}} \leq L\|u - v\|_{\mathcal{H}}$, so

$$\begin{cases} \exists\,! \, u \in C^1([0,+\infty[,\mathcal{H}) \text{ such that} \\ \dfrac{\partial u}{\partial t}(t) + Au(t) = f(t) \; \forall \; t \geq 0 \\ u(0) = u_0 \end{cases}$$

has unique solution.

And moreover if $a \in X$, ($X$ is Banach space) is the initial data, one can rewrites the Cauchy problem as an integral equation

$$u(t) = e^{-tL}a + \int_0^t e^{(s-t)L} f(u(s), \nabla u(s))\,ds := Nu(t)$$

when $> 0$ , is small enough and Banach space $X$ is appropriate, one can proves that $N$ is a contraction in some ball $B(a;r)$ or $C(0,T;X)$. Then there is a unique fixed point $u$; and this is the local solution.

## 🔶 Long time existence on semi-linear parabolic PDE:

In order to show long-time existence, we need to develop some a periori estimate, of the solution up to third order. Throught this subsection $u$ denote the solution to the initial value problem.

$$\begin{cases} \dfrac{\partial}{\partial t}u = \log \det\left(g_{i\bar{j}} + \dfrac{\partial^2 u}{\partial z^i \partial \bar{z}^j}\right) - \log\det\left(g_{i\bar{j}}\right) + f \\ u(x,0) = 0 \quad t = 0 \end{cases} \qquad (5.6)$$

On the maximal time interval $[0,T)$, such that $\tilde{g}_{i\bar{j}} = g_{i\bar{j}} + \dfrac{\partial^2 u}{\partial z^i \partial \bar{z}^j}$ is positive definite and hence defines a Kähler metric on $M$ for any time $t \in [0,T)$.

By differentiating equation (5.6), we get

$$\frac{\partial}{\partial t}\left(\frac{\partial}{\partial t}u\right) = \tilde{g}^{i\bar{j}} \frac{\partial^2}{\partial z^i \partial \bar{z}^j}\left(\frac{\partial}{\partial t}u\right) = \hat{\Delta}\left(\frac{\partial}{\partial t}u\right)$$

where $\tilde{\Delta}$ is the Laplace operator of $\tilde{g}_{i\bar{j}}$ and applying the maximum principle (for more details see PDE's chapter) , we find out that

$$max_M \left|\frac{\partial}{\partial t}u\right| \leq max_M |f|$$

Lemma : Let $u_{min} = \inf_{M\times[0,T)} u$ , then there exist constants $c_1, c_2 > 0$ such that

$$0 < n + \Delta_g u \leq c_1 e^{c_2\,(u(t)-u_{min})} \qquad (5.7)$$

for all $t \in [0,T)$ .

proof . The first inequality (5.7) follows from the fact that $\tilde{g}_{i\bar{j}}$ is positive definite and because

$$\tilde{g}_{i\bar{j}} = g_{i\bar{j}} + \frac{\partial^2 u}{\partial z^i \partial \overline{z}^j}$$

So $g^{i\bar{j}}\tilde{g}_{i\bar{j}} = g^{i\bar{j}}g_{i\bar{j}} + g^{i\bar{j}}\frac{\partial^2 u}{\partial z^i \partial \overline{z}^j}$ therefore as $g^{i\bar{j}}g_{i\bar{j}} = n$ , $n + \Delta u$ is the trace of $\tilde{g}_{i\bar{j}}$ with respect to $g_{i\bar{j}}$ .

For the second inequality, we use of a lemma of S.-T.-Yau

**Lemma (Yau inequality I [14])** We have

$$\left(\tilde{\Delta} - \frac{\partial}{\partial t}\right)\left(exp(-c_0 u)(n + \Delta u)\right)$$

$$\geq -\exp(-c_0 u)\left(\Delta f + n^2 inf_{i \neq 1}\left(R_{i\bar{i}1\bar{1}}\right)\right) - c_0 exp(-c_0 u)\left(n - \frac{\partial u}{\partial t}\right)(n + \Delta u)$$

$$+ \left(c_0 + inf_{i \neq 1}\left(R_{i\bar{i}1\bar{1}}\right)exp(-c_0 u)\right)exp\left(-\frac{\partial u}{\partial t} + f\right)(n + \Delta u)^{\frac{n}{n-1}}$$

(5.8)

where $R_{i\bar{i}1\bar{1}}$ is the bisectional curvature of the metric $g_{i\bar{j}}$ and $c_0$ is a positive constant such that $c_0 + inf_{i \neq 1}\left(R_{i\bar{i}1\bar{1}}\right) > 0$.

For any given $t \in (0, T)$, we assume function $exp(-c_0 u)(n + \Delta u)$ achieves its maximum at point $(p, t_0)$, with $t_0 > 0$ , on $[0, T] \times M$. Then Yau showed at this point the left hand side of (5.8) is non-positive and therefore we have

$$0 \geq -\left(\Delta f + n^2 inf_{i \neq 1}\left(R_{i\bar{i}1\bar{1}}\right)\right) - c_0\left(n - \frac{\partial u}{\partial t}\right)(n + \Delta u)$$

$$+ \left(c_0 + inf_{i \neq 1}\left(R_{i\bar{i}1\bar{1}}\right)\right)exp\left(-\frac{\partial u}{\partial t} + \frac{f}{n-1}\right)(n + \Delta u)^{\frac{n}{n-1}}$$

So because we know $max_M \left|\frac{\partial}{\partial t}u\right| \leq max_M |f|$, therefore by using this fact in previous inequality we get

$$(n + \Delta u)^{\frac{n}{n-1}} \leq \acute{c}(1 + (n + \Delta u)) \quad (5.9)$$

where $\acute{c}$ is a positive constant independent of $t$. So from (5.9) we get some elementary computation we get

$$n + \Delta u(p, t_0) \le c_1$$

Hence $exp(-c_0 u)(n + \Delta u) \le c_1 exp(-c_0 u(p, t_0))$ and it follows that

$$n + \Delta u \le c_1 exp\left(c_0\left(u - inf_{M \times [0,T)} u\right)\right)$$

So we obtain the desired result and the proof is complete ∎

Also we try to present an equivalent inequality like Yau inequality with complete detail in its proof and next we show $C^2$ −estimate by uniform estimation of Laplacian.

At first we start with some lemmas which will be useful later.

**Lemma**: If $\varphi \in C^2(M)$, the following assertions are equivalent

1. $\varphi$ is admissible

2. $YC(\varphi) := \frac{\det(\omega + i\partial\bar{\partial}\varphi)}{\det(\omega)} > 0$

Let $\omega$ is a kahler form and $f: M \to \mathbb{R}$ be in $C^\infty$ class such that $\int_M e^f \omega^n = \int_M \omega^n$ then, there exists an unique function $\varphi: M \to \mathbb{R}$ of $C^\infty$ class such that

1. $\int_M \varphi \omega^n = 0$

2. (Yau − Calabi Equation): $\left(\omega + i\partial\bar{\partial}\varphi\right)^n = e^f \omega^n$

Indeed , the solution $\varphi$ is admissible and the Yau-Calabi equation can be written as follows

$$YC(\varphi) := \frac{\det(\tilde{g})}{\det(g)} = e^f$$

If $\varphi$ be a function of $C^5$ class, admissible on metric $g$, and solution of the Calabi-Yau equation $YC(\varphi) = e^f$ $(f \in C^3)$ with

$$f \le F_0 \text{ and } \Delta f \le F_1 \qquad (F_0 \text{ and } F_1 \text{ are constants})$$

Then such estimates gives a periori estimates on $\Delta\varphi$. So we start with a theorem originally from S.T.Yau.

**Theorem (S.T.Yau inequality II) :** There exists a constant $c$ depends on the curvature on metric $g$, and constants $A, B$ and $C$ with $C > 0$ depending to $n, F_0, F_1$ and curvature such that

$$(n - \Delta\varphi)\tilde{\Delta}(\text{Log}(n - \Delta\varphi) - c\varphi) \leq A + B(n - \Delta\varphi) - C(n - \Delta\varphi)^{\frac{n}{n-1}}$$

Note that, here, $\Delta$ is complex Laplacian on metric $g$ and $\tilde{\Delta}$ is complex Laplacian on metric $\tilde{g}$.

Proof : Calabi-Yau equation in local view can be written as

$$\log\left(\frac{\det(\tilde{g})}{\det(g)}\right) = \log(\det(\tilde{g})) - \log(\det(g)) = f$$

By first time defferentiating we have : $\partial_{\bar{k}}(\log(\det(g))) = \frac{\partial_{\bar{k}}(\det(g))}{\det(g)}$, or $\partial_{\bar{k}}(\det(g)) = d(\det(g)).\partial_{\bar{k}} = (d(\det)_g \circ d g).\partial_{\bar{k}} = d(\det)_g.\partial_{\bar{k}}g$, which implies that

$$\partial_{\bar{k}}(\det(g) = \det(g)tr(g^{-1}\partial_{\bar{k}}g)$$

Therefore, $\partial_{\bar{k}}(\log(\det(g))) = tr(g^{-1}\partial_{\bar{k}}g) = g^{a\bar{b}}\partial_{\bar{k}}g_{a\bar{b}}$, hence

$$\partial_{\bar{k}}f = \tilde{g}^{a\bar{b}}\partial_{\bar{k}}\tilde{g}_{a\bar{b}} - g^{a\bar{b}}\partial_{\bar{k}}g_{a\bar{b}} = \left(\tilde{g}^{a\bar{b}} - g^{a\bar{b}}\right)\partial_{\bar{k}}g_{a\bar{b}} + \tilde{g}^{a\bar{b}}\partial_{\bar{k}a\bar{b}}\varphi$$

By second time differentiating, we get

$$\partial_{j\bar{k}}f = \left(\partial_j\tilde{g}^{a\bar{b}} - \partial_j g^{a\bar{b}}\right)\partial_{\bar{k}}g_{a\bar{b}} + \left(\tilde{g}^{a\bar{b}} - g^{a\bar{b}}\right)\partial_{j\bar{k}}g_{a\bar{b}} + \partial_j\tilde{g}^{a\bar{b}}\partial_{\bar{k}a\bar{b}}\varphi + \tilde{g}^{a\bar{b}}\partial_{j\bar{k}a\bar{b}}\varphi,$$

Therefore $\Delta f$, can be written as follows

$$\Delta f = -g^{j\bar{k}}\left(\partial_j\tilde{g}^{a\bar{b}} - \partial_j g^{a\bar{b}}\right)\partial_{\bar{k}}g_{a\bar{b}} - \left(\tilde{g}^{a\bar{b}} - g^{a\bar{b}}\right)g^{j\bar{k}}\partial_{j\bar{k}}g_{a\bar{b}} - g^{j\bar{k}}\partial_j\tilde{g}^{a\bar{b}}\partial_{\bar{k}a\bar{b}}\varphi - \tilde{g}^{a\bar{b}}g^{j\bar{k}}\partial_{j\bar{k}a\bar{b}}\varphi$$

Let $x \in M$, when $x$ is belong to $g$ −normal chart , the previous expression reduces to

$$\Delta f = 0 - \sum_{j=1}^{m}\tilde{g}^{a\bar{b}}\partial_{j\bar{j}}g_{a\bar{b}} + \sum_{j,a=1}^{m}\partial_{j\bar{j}}g_{a\bar{a}} - \sum_{j=1}^{m}\partial_j\tilde{g}^{a\bar{b}}\partial_{\bar{j}a\bar{b}}\varphi - \sum_{j=1}^{m}\tilde{g}^{a\bar{b}}\partial_{j\bar{j}a\bar{b}}\varphi$$

$$= -\sum_{j=1}^{m}\tilde{g}^{a\bar{b}}\left(\partial_{j\bar{j}}g_{a\bar{b}} + \partial_{j\bar{j}a\bar{b}}\varphi\right) + \sum_{j,a=1}^{m}\partial_{j\bar{j}}g_{a\bar{a}} - \sum_{j=1}^{m}\partial_j\tilde{g}^{a\bar{b}}\partial_{\bar{j}a\bar{b}}\varphi$$

But we know $\partial_j\tilde{g}^{a\bar{b}} = -\tilde{g}^{a\bar{s}}\tilde{g}^{l\bar{b}}\partial_j\tilde{g}_{l\bar{s}}$, so when $x$ is in $g$-normal chart , it can be considered as $\partial_j\tilde{g}^{a\bar{b}} = -\tilde{g}^{a\bar{s}}\tilde{g}^{l\bar{b}}\partial_{jl\bar{s}}\varphi$ , hence we get

$$\Delta f = -\sum_{j=1}^{m}\tilde{g}^{a\bar{b}}\left(\partial_{j\bar{j}}g_{a\bar{b}} + \partial_{j\bar{j}a\bar{b}}\varphi\right) + \sum_{j,a=1}^{m}\partial_{j\bar{j}}g_{a\bar{a}} + \sum_{j=1}^{m}\tilde{g}^{a\bar{s}}\tilde{g}^{l\bar{b}}\partial_{jl\bar{s}}\varphi\partial_{\bar{j}a\bar{b}}\varphi$$

which can be written (by using curvature formula in $g$-normal chart, $R_{j\bar{k}b\bar{a}} = \partial_{j\bar{k}} g_{b\bar{a}}$)

$$\Delta f = -\sum_{j=1}^{m} \tilde{g}^{a\bar{b}}(R_{j\bar{j}a\bar{b}} + \partial_{j\bar{j}a\bar{b}}\varphi) + \sum_{j,a=1}^{m} R_{j\bar{j}a\bar{a}} + \sum_{j=1}^{m} \tilde{g}^{a\bar{s}}\tilde{g}^{l\bar{b}}\partial_{jl\bar{s}}\varphi\partial_{\bar{j}a\bar{b}}\varphi \qquad (\mathfrak{H})$$

Now we calculate $\tilde{\Delta}\Delta\varphi$. We know $\Delta\varphi = -g^{a\bar{b}}\partial_{a\bar{b}}\varphi$, So

$$\tilde{\Delta}\Delta\varphi = \tilde{g}^{j\bar{k}}\partial_{j\bar{k}}(g^{a\bar{b}}\partial_{a\bar{b}}\varphi) = \tilde{g}^{j\bar{k}}\partial_j(\partial_{\bar{k}}g^{a\bar{b}}\partial_{a\bar{b}}\varphi + g^{a\bar{b}}\partial_{\bar{k}a\bar{b}}\varphi)$$

$$= \tilde{g}^{j\bar{k}}(\partial_{j\bar{k}}g^{a\bar{b}}\partial_{a\bar{b}}\varphi + \partial_{\bar{k}}g^{a\bar{b}}\partial_{ja\bar{b}}\varphi + \partial_j g^{a\bar{b}}\partial_{\bar{k}a\bar{b}}\varphi + g^{a\bar{b}}\partial_{j\bar{k}a\bar{b}}\varphi)$$

But we know, $\partial_j \tilde{g}^{a\bar{b}} = -\tilde{g}^{a\bar{s}}\tilde{g}^{l\bar{b}}\partial_j \tilde{g}_{l\bar{s}}$, so when $x$, is in $g$-normal chart, we get $\partial_j \tilde{g}^{a\bar{b}} = 0$, and even we have $\partial_{\bar{k}}\tilde{g}^{a\bar{b}} = 0$ . Therefore when $x$ is in $g$ −normal chart we have.

$$\tilde{\Delta}\Delta\varphi = \tilde{g}^{j\bar{k}}\left(\partial_{j\bar{k}}g^{a\bar{b}}\partial_{a\bar{b}}\varphi + \sum_{a=1}^{m}\partial_{j\bar{k}a\bar{a}}\varphi\right)$$

By derivating on both sides of $\partial_{\bar{k}}g^{a\bar{b}} = -g^{a\bar{s}}g^{l\bar{b}}\partial_{\bar{k}}g_{l\bar{s}}$, we obtain, $\partial_{j\bar{k}}g^{a\bar{b}} = -\partial_j g^{a\bar{s}}g^{l\bar{b}}\partial_{\bar{k}}g_{l\bar{s}} - g^{a\bar{s}}\partial_j g^{l\bar{b}}\partial_{\bar{k}}g_{l\bar{s}} - g^{a\bar{s}}g^{l\bar{b}}\partial_{j\bar{k}}g_{l\bar{s}}$ , and , when $x$ is in $g$-normal chart , we know $\partial_{j\bar{k}}g^{a\bar{b}} = -\partial_{j\bar{k}}g_{b\bar{a}} = -R_{j\bar{k}b\bar{a}}$ , So we get

$$\tilde{\Delta}\Delta\varphi = -\sum_{a,b=1}^{m} \tilde{g}^{j\bar{k}}R_{j\bar{k}b\bar{a}}\,\partial_{a\bar{b}}\varphi + \sum_{a=1}^{m} \tilde{g}^{j\bar{k}}\,\partial_{j\bar{k}a\bar{a}}\varphi \qquad (\mathfrak{H}\mathfrak{H})$$

By combinating the equalities $(\mathfrak{H})$, $(\mathfrak{H}\mathfrak{H})$, we obtain,

$$\Delta f + \tilde{\Delta}\Delta\varphi = \sum_{j,a=1}^{m} R_{j\bar{j}a\bar{a}} - \sum_{j=1}^{m} \tilde{g}^{a\bar{b}}R_{j\bar{j}a\bar{b}} - \sum_{a,b=1}^{m} \tilde{g}^{j\bar{k}}R_{j\bar{k}b\bar{a}}\,\partial_{a\bar{b}}\varphi + \sum_{j=1}^{m} \tilde{g}^{a\bar{s}}\tilde{g}^{l\bar{b}}\partial_{jl\bar{s}}\varphi\partial_{\bar{j}a\bar{b}}\varphi$$

$$+ \underbrace{\sum_{a=1}^{m} \tilde{g}^{j\bar{k}}\partial_{j\bar{k}a\bar{a}}\varphi - \sum_{j=1}^{m} \tilde{g}^{a\bar{b}}\partial_{j\bar{j}a\bar{b}}\varphi}_{0}$$

When $x$ is in $g$-normal chart, always we have: $\partial_{a\bar{b}}\varphi = \varphi_{a\bar{b}}$, $\partial_{jl\bar{s}}\varphi = \varphi_{jl\bar{s}}$ and $\partial_{\bar{j}a\bar{b}}\varphi = \varphi_{\bar{j}a\bar{b}}$ .So

$$\Delta f + \tilde{\Delta}\Delta\varphi = \underbrace{\sum_{j,a=1}^{m} R_{j\bar{j}a\bar{a}} - \sum_{j=1}^{m} \tilde{g}^{a\bar{b}}R_{j\bar{j}a\bar{b}} - \sum_{a,b=1}^{m} \tilde{g}^{j\bar{k}}R_{j\bar{k}b\bar{a}}\,\varphi_{a\bar{b}}}_{:=\mathcal{A}} + \sum_{j=1}^{m} \tilde{g}^{a\bar{s}}\tilde{g}^{l\bar{b}}\varphi_{jl\bar{s}}\varphi_{\bar{j}a\bar{b}} \qquad (\varkappa)$$

Now here, we present the deagonalization of $g$ −normal chart in $x$ with $\tilde{g}$ in $x$ . If $g$ and $\tilde{g}$ are diagonal in $x$, then, $[\varphi_{a\bar{b}}]_{1 \le a,b \le m}$ , is diagonal, too. And we have $\tilde{g}^{a\bar{a}} = \frac{1}{\tilde{g}_{a\bar{a}}} = \frac{1}{1+\varphi_{a\bar{a}}}$ ; therefore we have

$$\mathcal{A} = \sum_{j,a=1}^{m} R_{j\bar{j}a\bar{a}} - \sum_{j,a=1}^{m} \frac{1}{1+\varphi_{a\bar{a}}} R_{j\bar{j}a\bar{a}} - \sum_{j,a=1}^{m} \frac{\varphi_{a\bar{a}}}{1+\varphi_{j\bar{j}}} R_{j\bar{j}a\bar{a}}$$

$$= \sum_{j,a=1}^{m} \frac{\varphi_{a\bar{a}}}{1+\varphi_{a\bar{a}}} R_{j\bar{j}a\bar{a}} - \sum_{a,j=1}^{m} \frac{\varphi_{a\bar{a}}}{1+\varphi_{j\bar{j}}} R_{j\bar{j}a\bar{a}} = \sum_{j,a=1}^{m} \frac{\varphi_{j\bar{j}} - \varphi_{a\bar{a}}}{(1+\varphi_{a\bar{a}})(1+\varphi_{j\bar{j}})} R_{j\bar{j}a\bar{a}} \, \varphi_{a\bar{a}} \quad \text{(٦)}$$

By exchanging $j$ and $a$ and using symmetric property of curvature： $R_{j\bar{j}a\bar{a}} = R_{a\bar{a}j\bar{j}}$ , we get

$$\mathcal{A} = \sum_{j,a=1}^{m} R_{j\bar{j}a\bar{a}} \, \varphi_{j\bar{j}} \frac{\varphi_{a\bar{a}} - \varphi_{j\bar{j}}}{(1+\varphi_{a\bar{a}})(1+\varphi_{j\bar{j}})} \qquad \text{(٦٦)}$$

By summing the equalities of (٦) and (٦٦) we obtain

$$2\mathcal{A} = \sum_{j,a=1}^{m} R_{j\bar{j}a\bar{a}} \frac{2\varphi_{j\bar{j}}\varphi_{a\bar{a}} - \varphi_{j\bar{j}}^2 - \varphi_{a\bar{a}}^2}{(1+\varphi_{a\bar{a}})(1+\varphi_{j\bar{j}})}$$

So

$$\mathcal{A} = -\frac{1}{2} \sum_{j,a=1}^{m} R_{j\bar{j}a\bar{a}} \frac{(\varphi_{j\bar{j}} - \varphi_{a\bar{a}})^2}{(1+\varphi_{a\bar{a}})(1+\varphi_{j\bar{j}})} = \frac{1}{2} \sum_{j \neq a \in \{1, \dots m\}} (-R_{j\bar{j}a\bar{a}}) \frac{(\tilde{g}_{j\bar{j}} - \tilde{g}_{a\bar{a}})^2}{\tilde{g}_{j\bar{j}}\tilde{g}_{a\bar{a}}}$$

Hence, if $\tilde{g}$, is deagonalizing of $g$ −normal chart in $x$, we have

$$\mathcal{A} \ge \left[ \inf_{j \neq a \in \{1, \dots m\}} (-R_{j\bar{j}a\bar{a}}) \right] \times \sum_{j,a=1}^{m} \frac{(\tilde{g}_{j\bar{j}} - \tilde{g}_{a\bar{a}})^2}{2\tilde{g}_{j\bar{j}}\tilde{g}_{a\bar{a}}}$$

Moreover,

$$\sum_{j,a=1}^{m} \frac{(\tilde{g}_{j\bar{j}} - \tilde{g}_{a\bar{a}})^2}{2\tilde{g}_{j\bar{j}}\tilde{g}_{a\bar{a}}} = \sum_{j,a=1}^{m} \left( \frac{\tilde{g}_{j\bar{j}}}{2\tilde{g}_{a\bar{a}}} + \frac{\tilde{g}_{a\bar{a}}}{2\tilde{g}_{j\bar{j}}} - 1 \right) = \left( \sum_{j,a=1}^{m} \frac{\tilde{g}_{j\bar{j}}}{\tilde{g}_{a\bar{a}}} \right) - m^2 = \left( \sum_{j,a=1}^{m} \tilde{g}^{a\bar{a}} \tilde{g}_{j\bar{j}} \right) - m^2$$

Which implies that,

$$\mathcal{A} \geq \left[\inf_{j \neq a \in \{1,\dots m\}}(-R_{j\bar{J}a\bar{a}})\right] \times \left(\sum_{j,a=1}^{m} \tilde{g}^{a\bar{a}}\tilde{g}_{j\bar{J}}\right) - m^2 \qquad (\varkappa\varkappa)$$

So by combining the equations of $(\varkappa)$ and $(\varkappa\varkappa)$, we obtain

$$\tilde{\Delta}\Delta\varphi \geq \underbrace{-\Delta f + \sum_{j=1}^{m} \tilde{g}^{a\bar{a}}\,\tilde{g}^{b\bar{b}}\varphi_{jb\bar{a}}\varphi_{\bar{J}b\bar{a}}}_{T_1} + \underbrace{\left[\inf_{j \neq a \in \{1,\dots m\}}(-R_{j\bar{J}a\bar{a}})\right]}_{T_2} \times \left(\left(\sum_{j,a=1}^{m} \tilde{g}^{a\bar{a}}\tilde{g}_{j\bar{J}}\right) - m^2\right)$$

Finally, we calculate $\tilde{\Delta}(Log(m - \Delta\varphi) - c\varphi)$, that $c > 0$ is a constant which subsequently choosen.

$$\tilde{\Delta}(Log(m - \Delta\varphi) - c\varphi) = -\tilde{g}^{j\bar{k}}\partial_j\partial_{\bar{k}}(Log(m - \Delta\varphi) - c\varphi)$$

$$= -\tilde{g}^{j\bar{k}}\partial_j\left(-(m - \Delta\varphi)^{-1}\partial_{\bar{k}}(\Delta\varphi) - c\partial_{\bar{k}}(\varphi)\right)$$

$$= -\tilde{g}^{j\bar{k}}\left(-(m - \Delta\varphi)^{-1}\partial_{j\bar{k}}(\Delta\varphi) + (m - \Delta\varphi)^{-2}\left(-\partial_j(\Delta\varphi)\right)\partial_{\bar{k}}(\Delta\varphi) - c\partial_{j\bar{k}}(\varphi)\right)$$

$$= -(m - \Delta\varphi)^{-1}\tilde{\Delta}(\Delta\varphi) - c\tilde{\Delta}\varphi + (m - \Delta\varphi)^{-2}\tilde{g}^{j\bar{k}}\partial_j(\Delta\varphi)\partial_{\bar{k}}(\Delta\varphi)$$

Hence in $x$, when $\tilde{g}$ is deagonalizing of $g$ −normal chart, we have

$$(m - \Delta\varphi)\tilde{\Delta}(\log(m - \Delta\varphi) - c\varphi) = -\tilde{\Delta}\Delta\varphi - c\tilde{\Delta}\varphi(m - \Delta\varphi) + \underbrace{(m - \Delta\varphi)^{-1}\tilde{g}^{j\bar{J}}\partial_j(\Delta\varphi)\partial_j(\Delta\varphi)}_{:=\mathfrak{B}}$$

Or, we can write $\Delta\varphi = -\sum_{a=1}^{m}\partial_{a\bar{a}}\varphi = -\sum_{a=1}^{m}(\tilde{g}_{a\bar{a}} - 1) = m - \sum_{a=1}^{m}\tilde{g}_{a\bar{a}}$, therefore

$$m - \Delta\varphi = \sum_{a=1}^{m}\tilde{g}_{a\bar{a}}$$

Moreover , $\partial_j(\Delta\varphi) = \partial_j\left(-g^{a\bar{b}}\partial_{a\bar{b}}\varphi\right) = -g^{a\bar{b}}\partial_{ja\bar{a}}\varphi = -\sum_{a=1}^{m}\partial_{ja\bar{a}}\varphi$, and $\partial_j(\Delta\varphi) = \overline{\partial_{\bar{J}}(\Delta\varphi)}$ (the function $\varphi$ and so $\Delta\varphi$ is real valued): so from this we get

$$\mathfrak{B} = \frac{\sum_{j=1}^{m}\tilde{g}^{j\bar{J}}\left|\sum_{a=1}^{m}\partial_{ja\bar{a}}\varphi\right|^2}{\sum_{a=1}^{m}\tilde{g}_{a\bar{a}}} = \frac{\sum_{j=1}^{m}\tilde{g}^{j\bar{J}}\left|\sum_{a=1}^{m}(\sqrt{\tilde{g}^{a\bar{a}}}\partial_{ja\bar{a}}\varphi)(\sqrt{\tilde{g}_{a\bar{a}}})\right|^2}{\sum_{a=1}^{m}\tilde{g}_{a\bar{a}}}$$

and by Cauchy-Schwarz, inequality, we deduce that

$$\mathfrak{B} \leq \frac{\sum_{j=1}^{m}\tilde{g}^{j\bar{J}}\left(\sum_{a=1}^{m}\left|\sqrt{\tilde{g}^{a\bar{a}}}\partial_{ja\bar{a}}\varphi\right|^2\right) \times \left(\sum_{a=1}^{m}\left|\sqrt{\tilde{g}_{a\bar{a}}}\right|^2\right)}{\sum_{a=1}^{m}\tilde{g}_{a\bar{a}}}$$

Therefore

$$\mathfrak{B} \leq \sum_{j,a=1}^{m} \tilde{g}^{j\bar{j}} \tilde{g}^{a\bar{a}} \left| \partial_{ja\bar{a}} \varphi \right|^2 \leq \sum_{j,a,b=1}^{m} \tilde{g}^{j\bar{j}} \tilde{g}^{a\bar{a}} \left| \partial_{jb\bar{a}} \varphi \right|^2 \quad (\text{ because } \tilde{g}^{k\bar{k}} \geq 0)$$

or

$$\sum_{j,a,b=1}^{m} \tilde{g}^{j\bar{j}} \tilde{g}^{a\bar{a}} \left| \partial_{jb\bar{a}} \varphi \right|^2$$

$$= \sum_{j,a,b=1}^{m} \tilde{g}^{j\bar{j}} \tilde{g}^{a\bar{a}} \partial_{jb\bar{a}} \varphi \partial_{\bar{j}\bar{b}a} \varphi = \sum_{j,a,b=1}^{m} \tilde{g}^{b\bar{b}} \tilde{g}^{a\bar{a}} \partial_{bj\bar{a}} \varphi \partial_{\bar{b}\bar{j}a} \varphi = \sum_{j,a,b=1}^{m} \tilde{g}^{b\bar{b}} \tilde{g}^{a\bar{a}} \partial_{jb\bar{a}} \varphi \partial_{\bar{j}\bar{b}a} \varphi$$

$$= \sum_{j=1}^{m} \tilde{g}^{b\bar{b}} \tilde{g}^{a\bar{a}} \varphi_{jb\bar{a}} \varphi \varphi_{\bar{j}a\bar{b}} \varphi = T_1$$

We know, $\partial_{jb\bar{a}} \varphi = \varphi_{jb\bar{a}}$ so, this show that

$$\mathfrak{B} \leq T_1$$

Moreover , $\tilde{\Delta}\varphi = -\tilde{g}^{j\bar{j}} \partial_{j\bar{j}}\varphi = -\sum_{j=1}^{m} \tilde{g}^{j\bar{j}} (\tilde{g}_{j\bar{j}} - 1) = -\sum_{j=1}^{m} (1 - \tilde{g}^{j\bar{j}}) = \sum_{j=1}^{m} \tilde{g}^{j\bar{j}} - m$, therefore by combining previous inequalities and last equality, we get

$$(m - \Delta\varphi)\tilde{\Delta}(\log(m - \Delta\varphi) - c\varphi)$$

$$\leq \Delta f - T_1 - T_2 \left( \left( \sum_{j,a=1}^{m} \tilde{g}^{a\bar{a}} \tilde{g}_{j\bar{j}} \right) - m^2 \right) - c \left( \sum_{j=1}^{m} \tilde{g}^{j\bar{j}} - m \right)(m - \Delta\varphi) + T_1$$

$$\leq (\Delta f + m^2 T_2) + cm(m - \Delta\varphi) - (c + T_2)(m - \Delta\varphi) \left( \sum_{j=1}^{m} \tilde{g}^{j\bar{j}} \right)$$

For continuing the proof, we need to following Corollary,

**Corollary** : Let $m \in \mathbb{N}, m \geq 2$. If $a_1, a_2, \dots a_m > 0$, therefore we have following inequality:

$$\sum_{i=1}^{m} \frac{1}{a_i} \geq \left( \sum_{i=1}^{m} a_i \right)^{\frac{1}{m-1}} \left( \prod_{i=1}^{m} \frac{1}{a_i} \right)^{\frac{1}{m-1}}$$

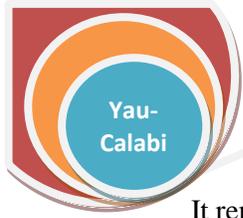

It remains to show the term $\sum_{j=1}^m \tilde{g}^{j\bar{j}}$, and for this, we use of previous corollary and apply corollary for $\tilde{g}^{j\bar{j}}$ and we obtain

$$\sum_{j=1}^m \tilde{g}^{j\bar{j}} = \sum_{j=1}^m \frac{1}{\tilde{g}_{j\bar{j}}} \geq \left(\sum_{j=1}^m \tilde{g}_{j\bar{j}}\right)^{\frac{1}{m-1}} \left(\prod_{j=1}^m \frac{1}{\tilde{g}_{j\bar{j}}}\right)^{\frac{1}{m-1}} = (m - \Delta\varphi)^{\frac{1}{m-1}} e^{-f\frac{1}{m-1}}$$

Since, $e^f = \frac{\det(\tilde{g})}{\det(g)} = \prod_{j=1}^m \tilde{g}_{j\bar{j}}$ .

Note that if $\mathfrak{A}_x(M)$, be unitary basis of $T_x^{\mathbb{C}}(M)$, then $\mathfrak{A}(M) = \bigcup \mathfrak{A}_x(M)$ is $\mathfrak{A}_{2m}(\mathbb{C})$ principal fiber of base M, where $\mathfrak{A}_{2m}(\mathbb{C}) := \{A \in M_{2m}(\mathbb{C}) : A^t\bar{A} = I\}$ is unitary group of order $2m$ $\mathfrak{A}(M)$ is compact manifold. In addition the function

$$\mathfrak{A}(M) \to \mathbb{C}$$

$$(e_1, \ldots, e_m, \bar{e_1}, \ldots, \bar{e_m}) \to T_2 := \inf_{j \neq a \in \{1, \ldots, m\}} \left(-R_{j\bar{j}a\bar{a}}\right)$$

is continuous, so $T_2$ is bounded on $\mathfrak{A}(M)$, in particular there exists $c > 0$ and $d > 0$ such that for each point and any choice of base unit $T_2 > 1 - c$ i.e, $T_2 + c > 1$ and $T_2 \leq d$. Hence

$$(m - \Delta\varphi)\tilde{\Delta}(\log(m - \Delta\varphi) - c\varphi) \leq \underbrace{(\Delta f + m^2 T_2)}_{A_x} + \underbrace{cm}_{B}(m - \Delta\varphi) - \underbrace{e^{\frac{-f}{m-1}}}_{C_x}(m - \Delta\varphi)^{\frac{m}{m-1}}$$

$A_x$ and $C_x$ are functions in $x$ . Moreover , $\Delta f \leq F_1$ and $T_2 \leq d$ thus, $A_z \leq F_1 + m^2 d := A$ and $f \leq F_0$ so $e^{\frac{-f}{m-1}} \geq e^{\frac{-F_0}{m-1}} := C > 0$ , we deduce that :

$$(m - \Delta\varphi)\tilde{\Delta}(\log(m - \Delta\varphi) - c\varphi) \leq A + B(m - \Delta\varphi) - C(m - \Delta\varphi)^{\frac{m}{m-1}}$$

$A$ depends to $F_1$, $m$ and curvature, $B$ depends to $m$ and curvature, and $C$ depends to $F_0$ and $m$. So proof is complete.

Now we consider a simple upper bound for $\Delta\varphi$.

**Theorem** :Let $\varphi$ is admissible function of class $C^2$ , then we have

$$\Delta\varphi < m$$

Proof: In fact, $\Delta\varphi = -g^{a\bar{b}}\partial_{a\bar{b}}\varphi = -g^{a\bar{b}}(\bar{g}_{a\bar{b}} - g_{a\bar{b}}) = tr(g^{-1}g) - tr(g^{-1}\bar{g}) = m - \underbrace{tr(g^{-1}\bar{g})}_{>0}$ ∎

So we try to find a lower bound for $\Delta\varphi$ as follows.

**Theorem**: Let $\varphi$ is a admissible function of class $C^5$ and be solution of Calabi-Yau equation $YC(\varphi) = e^f$ $(f \in C^3)$ with $f \leq F_0$ and $\Delta f \leq F_1$. Then there exists constant $c > 0$ depended to curvature and a constant $C_1 > 0$ depended to $m, F_0, F_1$ and curvature, such that

$$m - C_1 e^{c(\sup \varphi - \inf \varphi)} \leq \Delta\varphi < m$$

Proof: Assume that at point $x_0 \in M$ the function $\log(m - \Delta\varphi) - c\varphi$, take its maximum, then we have $\tilde{\Delta}(\log(m - \Delta\varphi) - c\varphi)(x_0) \geq 0$, so by applying S.T.Yau inequality II, we get $A + B(m - \Delta\varphi(x_0)) - C(m - \Delta\varphi(x_0))^{\frac{m}{m-1}} \geq 0$. In addition, the function $\wp: ]0, +\infty[ \to \mathbb{R}$ , $x \to A + Bx - Cx^{\frac{m}{m-1}}$ is continuous and tends to $-\infty$, when $x$ tend to $+\infty$ (because $C > 0$), so there is $C_1 > 0$ (depended to $A, B$ and $C$, thus $m, F_0, F_1$ and curvature) such that $\wp(x) < 0$ for all $x > C_1$ . So $m - \Delta\varphi(x_0) > 0$ and $\wp(m - \Delta\varphi(x_0)) \geq 0$ ,thus $m - \Delta\varphi(x_0) \leq C_1$ . But $\log(m - \Delta\varphi) - c\varphi \leq \log(m - \Delta\varphi(x_0)) - c\varphi(x_0)$ , so by applying exponential function on both sides, we get $m - \Delta\varphi \leq C_1 e^{c(\varphi - \varphi(x_0))} \leq C_1 e^{c(\sup \varphi - \inf \varphi)}$ $(c > 0)$, so we conclude $\Delta\varphi \geq m - C_1 e^{c(\sup \varphi - \inf \varphi)}$. ∎

## ✦ Zero order estimate

For obtaining the zero estimate of $u$ under normalization again we use of a lemma of S.-T.-Yau . Also San proved the same result in the context of open manifolds.

**Lemma (Yau[14])**: There exists positive constants $C_2$ and $C_3$ such that

$$\sup_{M \times [0,T)} v \leq C_2 \quad , \quad \sup_{M \times [0,T)} \int |v| dV \leq C_3$$

Where $v = u - \frac{1}{Vol(M)} \int u dV$

Now, we are going to get a lower bound for the function $v$. The strategy is to use the Nash-Moser iteration process. We want to bound $L^p$ −norms of the function $v$ by lower $L^p$ −norms, inductively.

**Lemma** : There exists a constant $C_4 > 0$ such that

$$\sup_{M \times [0,T)} |v| \leq C_4$$

**Proof**: The volume forms of the metrics $g$ and $\hat{g}$ are given, respectively, by

$$dV = \det\left(g_{i\bar{j}}\right) \wedge_{i=1}^{n} \left(\frac{\sqrt{-1}}{2} dz^i \wedge d\overline{z}^j\right) = \frac{\omega^n}{n!}$$

$$d\hat{V} = \det\left(\hat{g}_{i\bar{j}}\right) \wedge_{i=1}^{n} \left(\frac{\sqrt{-1}}{2} dz^i \wedge d\overline{z}^j\right) = \frac{\hat{\omega}^n}{n!}$$

Where $\omega = \frac{\sqrt{-1}}{2} g_{i\bar{j}} dz^i \wedge d\overline{z}^j$ and $\hat{\omega} = \frac{\sqrt{-1}}{2} \hat{g}_{i\bar{j}} dz^i \wedge d\overline{z}^j$ .

But we have

$$\frac{\partial u}{\partial t} = \log det\left(g_{i\bar{j}} + \frac{\partial^2 u}{\partial z^i \partial \overline{z}^j}\right) - \log \det\left(g_{i\bar{j}}\right) + f$$

So $\left(\frac{\partial u}{\partial t} - f\right) = \log \det \tilde{g}_{i\bar{j}\,i\bar{j}} - \log \det g_{i\bar{j}}$ . Therefore $exp\left(\frac{\partial u}{\partial t} - f\right) = \frac{\det \tilde{g}_{i\bar{j}\,i\bar{j}}}{\det(g_{i\bar{j}})}$, Then we get

$d\tilde{V} = exp\left(\frac{\partial u}{\partial t} - f\right) dV$, so by this fact it follows that, for $p > 1$, we define $\tilde{v} = v - C_2 - 2$ ,

$$\frac{-1}{n!} \int_M \frac{(-\tilde{v})^{p-1}}{p-1} (\omega^n - \hat{\omega}^n) = -\int_M \frac{(-\tilde{v})^{p-1}}{p-1} (dV - d\hat{V}) = \int_M \frac{(-\tilde{v})^{p-1}}{p-1} \left(exp\left(\frac{\partial u}{\partial t} - f\right) - 1\right) dV \quad (5.10)$$

where we renormalized $v$ so that $\tilde{v} < -1$ which is certainly possible by previous lemma . On the other hand

$$-\int_M \frac{(-\tilde{v})^{p-1}}{p-1} (\omega^n - \hat{\omega}^n) = -\int_M \frac{(-\tilde{v})^{p-1}}{p-1} \left(\omega^n - \left(\omega + \frac{\sqrt{-1}}{2}\partial\overline{\partial}\tilde{v}\right)^n\right)$$

$$= \int_M \frac{(-\tilde{v})^{p-1}}{p-1} \left(\left(\frac{\sqrt{-1}}{2}\partial\overline{\partial}\tilde{v}\right) \wedge \sum_{j=1}^{n-1} \omega^j \wedge \hat{\omega}^{n-j-1}\right)$$

$$= \int_M (-\tilde{v})^{p-2} \left(\frac{\sqrt{-1}}{2}\partial\tilde{v} \wedge \overline{\partial}\tilde{v}\right) \wedge \sum_{j=1}^{n-1} \omega^j \wedge \hat{\omega}^{n-j-1}$$

where the last equality follows from the following remark

**Remark**: If $h: R \rightarrow R$ is an increasing function, then

$$\int_M h\omega^{n-1} \wedge \partial\overline{\partial}u = -\int_M \partial h \wedge \omega^{n-1} \wedge \overline{\partial}u$$

Proof: By applying Stokes Theorem we have

$$0 = \int_M \partial\left(h\omega^{n-1} \wedge \overline{\partial}u\right) = \int_M \partial h \wedge \omega^{n-1} \wedge \overline{\partial}u + \int_M h\omega^{n-1} \wedge \partial\overline{\partial}u$$

So proof is complete ∎

Now, we come back to the proof. We showed that

$$-\int_M \frac{(-\tilde{v})^{p-1}}{p-1}(\omega^n - \widetilde{\omega}^n) = \int_M (-\tilde{v})^{p-2}\left(\frac{\sqrt{-1}}{2}\partial v \wedge \overline{\partial}v\right) \wedge \sum_{j=1}^{n-1} \omega^j \wedge \widetilde{\omega}^{n-j-1}$$

$$\geq \int_M (-\tilde{v})^{p-2}\left(\frac{\sqrt{-1}}{2}\partial\tilde{v} \wedge \overline{\partial}\tilde{v}\right) \wedge \omega^{n-1}$$

$$(5.11)$$

where the last inequality follows from the fact that the terms of $\frac{\sqrt{-1}}{2}\left(\partial\tilde{v} \wedge \overline{\partial}\tilde{v}\right) \wedge \omega^j \wedge \widetilde{\omega}^{n-j-1}$ are all nonnegative so from (5.10) and (5.11) we get

$$\int_M (-\tilde{v})^{p-2}|\nabla\tilde{v}|^2 dV \leq n! \int_M \frac{(-\tilde{v})^{p-1}}{p-1}\left(exp\left(\frac{\partial u}{\partial t} - f\right) - 1\right) dV \qquad (5.12)$$

where $|\nabla\tilde{v}|^2 = g^{i\overline{j}}\frac{\partial\tilde{v}}{\partial z^i}\frac{\partial\tilde{v}}{\partial\overline{z}^j}$ and also we used of this fact that $\partial\tilde{v} \wedge \overline{\partial}\tilde{v} = \frac{\partial\tilde{v}}{\partial z^i} \wedge \frac{\partial\tilde{v}}{\partial\overline{z}^j} dz^i \wedge d\overline{z}^j$ . Since

$$(-\tilde{v})^{p-2}|\nabla\tilde{v}|^2 = 4p^{-2}\left|\nabla(-\tilde{v})^{\frac{p}{2}}\right|^2$$

It follows from (5.12) that

$$\int_M \left|\nabla(-\tilde{v})^{\frac{p}{2}}\right|^2 dV \leq C\frac{p^2}{p-1}\int (-\tilde{v})^{p-1} dV$$

and hence

$$\left\|(-\tilde{v})^{\frac{p}{2}}\right\|_{H^1}^2 = \int_M \left|\nabla(-\tilde{v})^{\frac{p}{2}}\right|^2 dV + \int_M (-\tilde{v})^p dV$$

(Denote as usual by $H^1$ the first Sobolev Space and $\int_{H^1} f^2 = \int_{L^2} f^2 + \int_{L^2} |\nabla f|^2$ .) Then

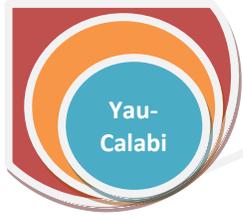

$$\left\| (-\tilde{v})^{\frac{p}{2}} \right\|_{H^1}^2 \leq C \frac{p^2}{p-1} \int_M (-\tilde{v})^p \, dV \leq Cp \int_M (-\tilde{v})^p \, dV \qquad (5.13)$$

For $p > 1$ and sufficiently large. Also inequality (5.13) holds for $p = 1$, simply replacing $\frac{(-\tilde{v})^{p-1}}{p-1}$ by $\log(-v)$ in the argument.

Now because $\left\| \tilde{v} - \frac{1}{Vol(M)} \int_M \tilde{v} \, dV \right\|$ is bounded (by previous lemma and definition of $\tilde{v}$), we apply the Sobolev inequality as follows

$$\left\| (-\tilde{v})^{\frac{p}{2}} \right\|_{L^{\frac{2n}{n-1}}}^2 \leq C \left\| (-\tilde{v})^{\frac{p}{2}} \right\|_{H^1}^2 \qquad (5.14)$$

Putting (5.13), (5.14) together we see that

$$\| \tilde{v} \|_{L^{\frac{pn}{n-1}}}^p \leq Cp \| \tilde{v} \|_{L^p}^p \quad \text{for } p = 1 \text{ or } p \gg 1.$$

Let $p = \gamma^j$ in (5.14), where $\gamma = \frac{n}{n-1}$ and $j = 0,1,2,...$, then by induction on $j$ we obtain

$$\| \tilde{v} \|_{L^{\gamma^{j+1}}} < C^{\sum_{k=0}^{j} \frac{1}{\gamma^k}} . \gamma^{\sum_{k=0}^{j} \frac{k}{\gamma^k}} . C_3$$

Letting $j \to \infty$, we have

$$\| \tilde{v} \|_{L^\infty} \leq C_4$$

So $\sup_{M \times [0,T)} |\tilde{v}| \leq C_4$ and therefore $\sup_{M \times [0,T)} |v| \leq C_4$ and proof is complete ∎ .

**Proposition :** There exists a constant $C_5$ such that

$$0 < n + \Delta v \leq C_5$$

Proof: By applying the previous Lemma and this fact that $\inf -u = -\sup u$ we get

$$0 < n + \Delta v = n + \Delta u \leq C_1 \exp\left( C_0 \left( u - \inf_{M \times [0,T)} u \right) \right) \leq C_1 \exp\left( C_0 \left( v - \inf_{M \times [0,T)} v \right) \right) \leq C_5$$

So the proof is complete ∎ .

## 🔆 First order estimate

For the first order estimate, we have the following result originally from Schauder theory

**Lemma** (Schauder estimate [23]) :There exists a constant $C_6$ such that

$$\sup_{\text{M}\times[0,\text{T})} |\nabla v| \leq C_6 \left( \sup_{\text{M}\times[0,\text{T})} |\Delta v| + \sup_{\text{M}\times[0,\text{T})} |v| \right)$$

**Remark**: In reality constant $C_6$ depends on the bounds of the coefficients of $\Delta = g^{i\bar{j}} \frac{\partial^2 u}{\partial z_i \partial \bar{z}_j}$ , so it depends to the bounds on the metric and because we will prove in next Proposition that the metrics $\tilde{g}_{i\bar{j}}(t)$ are uniformly equivalent to the metric $g_{i\bar{j}}$ , so by this reason we have uniform bound independent of $t$.

**Corollary**: There exists a constant $C_7$ such that

$$\sup_{\text{M}\times[0,\text{T})} |\nabla v| \leq C_7$$

Proof: By applying previous lemmas we obtain the desired estimate, so proof is complete ∎ .

## 🔆 Second order estimate

At first we show that the metrics $\tilde{g}_{i\bar{j}}(t)$ are uniformly equivalent to the metric $g_{i\bar{j}}$

**Proposition:** The metrics $\tilde{g}_{i\bar{j}}(t)$ are uniformly equivalent to the metric $g_{i\bar{j}}$

Proof: According to previous Corollary we know that $A = n + \Delta u \leq C_5$ and also $n = \tilde{g}^{i\bar{j}} g_{i\bar{j}} + \Delta_{\tilde{g}(t)} u$ because

$$\tilde{g}^{i\bar{j}} g_{i\bar{j}}(t) + \Delta_{\tilde{g}(t)} u = \tilde{g}^{i\bar{j}} g_{i\bar{j}}(t) + \tilde{g}^{i\bar{j}} \frac{\partial^2 u}{\partial z^i \partial \bar{z}^j} = \tilde{g}^{i\bar{j}} \left( g_{i\bar{j}}(t) + \frac{\partial^2 u}{\partial z^i \partial \bar{z}^j} \right) = \tilde{g}^{i\bar{j}} \tilde{g}_{i\bar{j}}(t) = n$$

So suppose that $\lambda_i$ denote the eigenvalues of $\tilde{g}_{i\bar{j}}(t)$ with respect to $g_{i\bar{j}}$ then we know

$$\tilde{g}_{i\bar{j}} = g_{i\bar{j}} + \frac{\partial^2 u}{\partial z^i \partial \bar{z}^j}$$

So

$$g^{i\bar{j}}\tilde{g}_{i\bar{j}} = g^{i\bar{j}}g_{i\bar{j}} + g^{i\bar{j}}\frac{\partial^2 u}{\partial z^i \partial \bar{z}^j} = n + \Delta_{g(t)}u$$

So $A = \sum_{i=1}^{n} \lambda_i$ and the eigenvalues of the matrix $\left(\frac{\partial^2 u}{\partial z^i \partial \bar{z}^j}\right)$ with respect to $g_{i\bar{j}}(t)$ are $\lambda_i - 1$ . But we know $A \le C_5$ so $\sum_{i=1}^{n} \lambda_i \le C_5$ therefore as long as $g_{i\bar{j}}$ is positive –definite, $\lambda_i \le C_5$ , for each $i$. But the product of the eigenvalues of the matrix $\left(\frac{\partial^2 u}{\partial z^i \partial \bar{z}^j}\right)$ with respect to $g_{i\bar{j}}$ is $\det(\tilde{g}_{i\bar{j}}(t))$ with respect to $\det(g_{i\bar{j}}(t))$ (product of eigenvalues of a matrix is exactly determinant of matrix) which means $\frac{\det(\tilde{g}_{i\bar{j}}(t))}{\det(g_{i\bar{j}}(t))}$ . But this quantity is exactly (as we proved it before) $exp\left(\frac{\partial u}{\partial t} - f\right)$ which is bounded. So this means that $A \le \frac{\partial^2 u}{\partial z^i \partial \bar{z}^j} \le B$ for some positive constants $A$ and  . So we estimated all the second order derivatives and the proof is complete ■

#### 🞂 Third order estimate

For proving the third order estimate for the function $v$ we use of the Yau's Lemma

**Lemma** (S.-T.-Yau[14]) Let we have the quantity

$$S = g^{i\bar{r}}g^{\bar{j}s}g^{k\bar{t}}v_{i\bar{j}k}v_{\bar{r}s\bar{t}}$$

then there exists positive constants $C_8$, $C_9$ and $C_{10}$ such that

$$\left(\hat{\Delta} - \frac{\partial}{\partial t}\right)(S + C_8\Delta v) \ge C_9 S - C_{10} \qquad (5.15)$$

At the maximum point $p(t)$ of $S + C_8\Delta v$ at time $t$, the inequality (5.15) shows that $0 \ge C_9 S - C_{10}$

Therefore $C_9(S + C_8\Delta v) \le C_{10} + C_8 C_9 \Delta v$ at $(t)$. Since we have already estimated $\Delta v$ so it follows that $sup_{M\times[0,T)} S + C_8 \Delta v$  is bounded. Yau by this lemma proved that this fact gives us the estimate for all the third order derivative of  .

#### 🞂 **Long –time existence**

For proving long-time existence, we need some information about parabolic Schauder estimate. In this subsection we introduce parabolic Hölder and Sobolev Space and then we give a parabolic Schauder estimate for parabolic heat equation. At first we have to say that the reason for introducing these spaces is the heat operators maps between parabolic Hölder and Sobolev spaces . Good reference on parabolic Hölder and Sobolev spaces is Krylov [36]. We first recall $C^k-$Spaces, Hölder spaces, Sobolev spaces of maps $u: I \to X$ where $I \subset R$ is an open and bounded interval and $X$ is Banach Spaces. For $k \in N$ we define $C^k_{loc}(I; X)$ to be the space of $k-times$ continuously differentiable maps $u: I \to X$. We define the $C^k-$norm by

$$\|u\|_{C^k} = \sum_{j=0}^{k} \sup_{t \in I} \left\| \partial_t^j u(t) \right\|_X \quad \text{for } u \in C^k_{loc}(I; X)$$

whenever it is finite , and we define

$$C^k(I, X) = \left\{ u \in C^k_{loc}(I; X); \ \|u\|_{C^k} < \infty \right\}$$

Moreover for $\alpha \in (0,1]$ we define the $C^{k,\alpha}-$norm by

$$\|u\|_{C^{k,\alpha}} = \|u\|_{C^k} + \sup_{t \neq s \in I} \frac{\left\| \partial_t^k u(t) - \partial_t^k u(s) \right\|_X}{|t-s|^\alpha}$$

for $u \in C^k_{loc}(I; X)$. By $C^{k,\alpha}_{loc}(I; X)$ we denote the space of maps $u \in C^k_{loc}(I; X)$ with finite $C^{k,\alpha}-$norm on every $\subset I$ , and we define

$$C^{k,\alpha}(I, X) = \left\{ u \in C^{k,\alpha}_{loc}(I; X); \ \|u\|_{C^{k,\alpha}} < \infty \right\}$$

then $C^k(I, X)$ and $C^{k,\alpha}(I, X)$ are both Banach spaces.

Next we recall again the Sobolev spaces of maps $u: I \to X$. Let $k \in N$ and $\in [1, \infty)$ . For a $k-times$ weakly differentiable map $u: I \to X$ we define the $W^{k,p}-$norm by

$$\|u\|_{W^{k,p}} = \left( \sum_{j=0}^{k} \int_I \left\| \partial_t^j u(t) \right\|_X^p dt \right)^{\frac{1}{p}}$$

whenever it is finite. We denote by $W^{k,p}_{loc}(I; X)$ the space of $k-times$ weakly differentiable maps $u: I \to X$ with finite $W^{k,p}$-norm on every $J \subset I$, and we define

$$W^{k,p}(I,X) = \left\{ u \in W^{k,p}_{loc}(I;X); \; \|u\|_{W^{k,p}} < \infty \right\}$$

then $W^{k,p}(I,X)$ is a Banach space. If $k = 0$, then we write $L^p_{loc}(I;X)$ and $L^p(I;X)$ instead of $W^{0,p}_{loc}(I;X)$ and $W^{0,p}(I;X)$, respectively. Now, we are ready to define parabolic $C^k-$ spaces and parabolic Hölder spaces. Let $(M,g)$ be a Reimann manifold and $k, l \in \mathbb{N}$ with $2k < l$. We define

$$C^{k,l}(I \times M) = \bigcap_{j=0}^{k} C^j\big(I; C^{l-2j}(M)\big)$$

then $C^{k,l}(I \times M)$ is a Banach space with the norm given by

$$\|u\|_{C^{k,l}} = \sum_{i,j} \sup_{(t,x) \in I \times M} \left| \partial_t^i \nabla^j u(t,x) \right|$$

for $u \in C^{k,l}(I \times M)$. where the sum is taken over $i = 1,2,\dots,k$ and $j = 1,\dots,l$ with $2i + j \leq l$. If $\alpha \in (0,1)$, then we define the parabolic Hölder space $C^{k,l,\alpha}(I \times M)$ by

$$C^{k,l,\alpha}(I \times M) = \bigcap_{j=0}^{k} C^{j,\frac{\alpha}{2}}\big(I, C^{l-2j}(M)\big) \cap C^j\big(I, C^{l-2j,\alpha}(M)\big)$$

then $C^{k,l,\alpha}(I \times M)$ is a Banach space with norm given by

$$\|u\|_{C^{k,l,\alpha}} = \sum_{i,j} \left\{ \sup_{(t,x) \in I \times M} \left| \partial_t^i \nabla^j u(t,x) \right| + \sup_{x \in M} \big[ \partial_t^i \nabla^j u(.,x) \big]_{\frac{\alpha}{2}} + \sup_{t \in I} \big[ \partial_t^i \nabla^j u(t,.) \big]_{\alpha} \right\}$$

for $u \in C^{k,l,\alpha}(I \times M)$, where the sum is taken over $i = 1,2,\dots,k$ and $j = 1,\dots,l$ with $2i + j \leq l$. Thus a function $u: I \times M \to R$ lies in $C^{k,l,\alpha}(I \times M)$ if and only if all derivatives of the form $\partial_t^i \nabla^j u$ with $i \leq k$, $j \leq l$, and $2i + j \leq l$ exist and are Hölder continuous in time with Hölder exponent $\alpha/2$ and Hölder continuous on $M$ with Hölder exponent $\alpha$.

Finally we define parabolic Sobolev space, Let $k, l \in N$ with $2k \leq l$, and $p \in [1,\infty)$, Then we define the parabolic Sobolev space $W^{k,l,p}(I \times M)$ by

$$W^{k,l,p}(I \times M) = \bigcap_{j=0}^{k} W^{j,p}(I; W^{l-2j,p}(M))$$

Then $W^{k,l,p}(I \times M)$ is a Banach space with norm given by

$$\|u\|_{W^{k,l,p}} = \left( \sum_{i,j} \int_I \int_M |\partial_t^i \nabla^j u(t,.)|^p \, dV_g \, dt \right)^{\frac{1}{p}}$$

for $u \in W^{k,l,p}(I \times M)$ where the sum is taken over $i = 1, 2, \dots, k$ and $j = 1, \dots, l$ with $2i + j \leq l$. Thus $W^{k,l,p}(I \times M)$ is the space of functions $u: (0,T) \times M \to R$, such that all weak derivatives of the form $\partial_t^i \nabla^j u$ with $i \leq k, j \leq l$, and $2i + j \leq l$ lie in $W^{0,0,p}(I \times M)$. Note that $W^{0,0,p}(I \times M) = L^p(I \times M)$.

## 🌱 Schauder estimate on parabolic heat equation

Let $\Omega \subset R^m$ be an open and bounded domain and $T > 0$. Let $a^{ij}, b^j, c: (0,T) \times \Omega \to R$ be continuous with $a^{ij} = a^{ji}$ for $i, j = 1, \dots m$ and define a linear differential operator $L$ acting on function $u \in C^{1,2}((0,T) \times \Omega)$ by

$$Lu = \frac{\partial u}{\partial t} - \frac{\partial}{\partial x^i}\left( a^{ij}(t,x) \frac{\partial u}{\partial x^j} \right) - b^i(t,x) \frac{\partial u}{\partial x^i} - c(t,x)u. \qquad (5.16)$$

the functions $a^{ij}$, $b^i$, and $c$ are called the coefficients of $L$. We assume that (5.16) is parabolic. This means that there exists a constant $\lambda > 0$, such that

$$\lambda^{-1}|\xi|^2 \leq a^{ij}(t,x)\xi_i\xi_j \leq \lambda|\xi|^2$$

for $(t,x) \in (0,T) \times \Omega$ and $\xi = (\xi_1, \xi_2, \dots, \xi_m) \in R^m$.

From now on let us assume that the coefficients of $L$ are smooth on $(0,T) \times \Omega$. Then we have the following Schauder estimates for solutions of $Lu = f$.

**Thesorem(Schauder estimate I [36-37]):** Let $k \in N$ with $k \geq 2$ and $\alpha \in (0,1)$. Let $u \in C^{1,k,\alpha}((0,T) \times \Omega)$, $f \in C^{0,k-2,\alpha}((0,T) \times \Omega)$, and assume that $Lu = f$. Then for every $\hat{\Omega} \subset \Omega$, there exists a constant $C > 0$ depending only on $\hat{\Omega}, k, \alpha, \lambda$, and the $C^{0,k}$ −norm of the coefficients of $L$ on $(0,T) \times \Omega$, such that

$$\|u\|_{C^{1,k,\alpha}} \leq C\left( \|f\|_{C^{0,k-2,\alpha}} + \|u\|_{C^{0,0}} \right)$$

where the norm on the left side is on $(0,T) \times \hat{\Omega}$ and the norm on the right side is on $(0,T) \times \Omega$.

Now we assume the $L^p$ −estimate for second order linear parabolic equations , which can be found in Krylov[36].

**Thesorem(Schauder estimate II [36-37]):** Let $k \in N$ with $k \geq 2$ and $p \in (1,\infty)$ . Let $u \in W^{1,2,p}\big((0,T) \times \Omega\big)$, $f \in W^{0,k-2,p}\big((0,T) \times \Omega\big)$, and assume that $Lu = f$ . Then $u \in W^{1,k,p}\big((0,T) \times \hat{\Omega}\big)$ for every $\hat{\Omega} \subset \Omega$ . Moreover for every $\hat{\Omega} \subset \Omega$ there exists a constant $C > 0$ depending only on $\hat{\Omega}$ ,$k, p, \lambda$, and the $C^{0,k}$-norm of the coefficients of $L$ on $(0,T) \times \Omega$ such that

$$\|u\|_{W^{1,k,p}} \leq C\big(\|f\|_{W^{0,k-2,p}} + \|u\|_{W^{0,0,p}}\big)$$

where the norm on the left side is on $(0,T) \times \hat{\Omega}$ and the norm on the right side is on $(0,T) \times \Omega$ .

Now we present another simple version of parabolic Schauder estimate.

**Definition**: Let $Q_T = \Omega \times (0,T)$ then, the parabolic boundary $\partial_p Q_T$ is defined as:

$$\partial_p Q_T = (\partial\Omega \times (0,T]) \cup (\Omega \times \{0\})$$

**Definition (parabolic Hölder space):** We say that $u \in C^{\alpha}(Q_T)$ if

$$|u(x_1,t_1) - u(x_2,t_2)| \leq C\left(|x_1 - x_2| + \sqrt{|t_1 - t_2|}\right)^{\alpha}$$

Also we say $u \in C^{2+\alpha}(Q_T)$ if

$$u, u_t, D_X u, D_X^2 u \in C^{\alpha}(Q_T)$$

**Theorem:** (Schauder estimate) If $u$ is a smooth solution of

$$\begin{cases} \dfrac{\partial u}{\partial t} - \Delta u = f, & in \ Q_T \\ u = g, & on \quad \partial_p Q_T \end{cases}$$

then

$$\|u\|_{C^{2+\alpha}(Q_T)} \leq C\left(\|u\|_{C^0(Q_T)} + \|f\|_{C^{\alpha}(Q_T)} + \|g\|_{C^{2,\alpha}(\partial Q_T)}\right)$$

Now we are ready to study on Long-time existence. By applying first, second and third order estimate we are in position to prove the Long-time existence of Kähler Ricci potential $u$.

**Theorem:** Let $u$ be the solution of

$$\begin{cases} \frac{\partial u}{\partial t} = \log \det \left( g_{i\bar{j}} + \frac{\partial^2 u}{\partial z^i \partial \bar{z}^j} \right) - \log \det \left( g_{i\bar{j}} \right) + f \\ u(x,t) = 0, \qquad\qquad t = 0 \end{cases} \tag{5.17}$$

On the maximal time interval $0 \leq t < T$ and let $v$ be the normalization of $u$ ($v = u - Ave_M u$).

Then the $C^\infty$ −norm (informally , because we don't have $C^\infty$ −norm) of $v$ are uniformly bounded for all $t \in (0,T)$ and consequently $= \infty$ . More over there exists a time sequence $t_n \to \infty$ such that $v(x, t_n)$ converges in $C^\infty$ topology to a smooth function $v_\infty(x)$ on $M$ as $\to \infty$ .

**Proof**: By differentiating (5.17) with respect to $z^k$ we obtain

$$\left( \tilde{\Delta} - \frac{\partial}{\partial t} \right) \left( \frac{\partial u}{\partial z^k} \right) = g^{i\bar{j}} \frac{\partial}{\partial z^k} \left( g_{i\bar{j}} \right) - \hat{g}^{i\bar{j}} \frac{\partial}{\partial z^k} \left( g_{i\bar{j}} \right) \tag{5.18}$$

then the coefficients of operator $\left( \tilde{\Delta} - \frac{\partial}{\partial t} \right)$ are bounded in $C^{0,\alpha}$ norm (because the third derivatives of $v$ are bounded, we have $C^{2,\alpha}$ bounded on $v$ . So the coefficients of operator $\tilde{\Delta} - \frac{\partial}{\partial t}$ are bounded in $C^{0,\alpha}$ norm) and also the right hand side of (5.18) also has estimate in $C^{0,\alpha}$ norm for all $0 < \alpha < 1$ . If we explain more, this is because of "third order estimate for $S$" . It implies that the Laplacian of $v$ (with respect to the fixed metric) is bounded in $C^1$ . So in particular in $C^\alpha$ for any $0 < \alpha < 1$ .

From the parabolic Schauder theory we know then $\frac{\partial u}{\partial z^k}$ has uniform $C^{2,\alpha}$ estimate and similarly for $\frac{\partial u}{\partial \bar{z}^k}$ . So the coefficients of $\hat{\Delta} - \frac{\partial}{\partial t}$ and the right hand side of (5.18) has uniform $C^{1,\alpha}$ estimate. Apply the parabolic Schauder estimate again we see that $\frac{\partial u}{\partial z^k}$ and $\frac{\partial u}{\partial \bar{z}^k}$ have uniform $C^{3,\alpha}$ estimate. By iteration we conclude that the $C^\infty$ −norm of $v(x,t)$ are uniformly bounded for all $t \in (0,T)$ and consequently we can select a time sequence $t_n \to \infty$ so that $v(x, t_n)$ converge to a smooth function $v_\infty(x)$ as $t \to \infty$ . On the other hand, since $\frac{\partial u}{\partial t}$ is uniformly bounded in t , the function $u$ cannot blow up in finite time. So the solution $u$ exists for all time, therefore proof is complete ∎

### The uniform convergence of the potential $u(t)$:

This section will be devoted to the proof of uniform convergence of the normalized potential $v$, as well as to show that $\frac{\partial u}{\partial t}$ converges to a constant as $t \to \infty$.

In [38], Cao used a slight generalization of the Li-Yau Harnack , inequality developed in [14]. We state Cao's version here.

**Theorem (Cao[38])** :Let $M$ be an $n$-dimensional compact manifold and let $g_{i\bar{j}}(t)$ , $t \in [0, \infty)$, be a family of Kähler metrics satisfying

     ***i.***     $C g_{i\bar{j}}(0) \leq g_{i\bar{j}}(t) \leq C^{-1} g_{i\bar{j}}(0)$

     ***ii.***    $\left| \frac{\partial}{\partial t} g_{i\bar{j}} \right| (t) \leq C g_{i\bar{j}}(0)$

     ***iii.***   $R_{i\bar{j}}(t) \geq -k g_{i\bar{j}}(0)$

Where $C$ and $k$ are positive constants independent of . Let $\Delta_t$ denote the Laplace operator of the metric $g_{ij}(t)$ .If $\varphi(x, t)$ is a positive solution for the equation

$$\left( \Delta_t - \frac{\partial}{\partial t} \right) \varphi(x, t) = 0$$

on $M \times [0, \infty)$, then for any $\alpha > 1$, we have.

$$\sup_{x \in M} \varphi(x, t_1) \leq \inf \varphi(x, t_2) \left( \frac{t_2}{t_1} \right)^{\frac{n}{2}} exp \left( \frac{1}{4(t_2 - t_1)} C_2^2 d^2 + \left( \frac{n\alpha k}{2(\alpha - 1)} + C_2 C_3 (n + A) \right) (t_2 - t_1) \right)$$

where $d$ is the diameter of $M$ measured by the metric $g_{ij}(0)$ , $A = sup \| \nabla^2 \log \varphi \|$ and $0 < t_1 < t_2 < \infty$.

Let $F = \frac{\partial}{\partial t} u$, where $u$ is the Kähler Ricci potential. Then , $F$ is the solution to

$$\begin{cases} \left( \Delta_t - \frac{\partial}{\partial t} \right) F = 0 \\ F(x, 0) = f(x) \end{cases} \qquad (5.19)$$

By the maximum principle for the parabolic equation, we have that

$$\sup_{x \in M} F(x, t_2) < \sup_{x \in M} F(x, t_1) < \sup_{x \in M} f(x) \quad (A)$$

where $t_2 > t_1 > 0$,

$$\inf_{x \in M} F(x, t_2) > \inf_{x \in M} F(x, t_1) > \inf_{x \in M} f(x) \quad (B)$$

In Cao's theorem because

$$\hat{g}_{i\bar{j}} = g_{i\bar{j}} + \frac{\partial^2 u}{\partial z^i \partial \bar{z}^j}$$

and $\frac{\partial^2 u}{\partial z^i \partial \bar{z}^j}$ is bounded, so we can replace $\tilde{g}_{i\bar{j}}$ with $g_{i\bar{j}}$ in condition $i, ii$ and $iii$ in previous theorem.

Let us define auxiliary functions

$$\varphi_n(x, t) = \sup_{x \in M} F(x, n-1) - F(x, (n-1)+t)$$

$$\psi_n(x, t) = F(x, (n-1)+t) - \inf_{x \in M} F(x, n-1)$$

$$\omega(t) = \sup_{x \in M} F(x, t) - \inf_{x \in M} F(x, t)$$

which are all positive functions and satisfy in (5.19) . Applying Cao's theorem to these functions for times $t_1 = \frac{1}{2}, t_2 = 1$, we obtain

$$\omega(n-1) + \omega\left(n - \frac{1}{2}\right) \le \gamma\big(\omega(n-1) - \omega(n)\big)$$

where $\gamma$ is a constant independent of $n$, and hence

$$\omega(n) \le \delta\omega(n-1), \qquad , \quad \delta = \frac{\gamma - 1}{\gamma} < 1$$

By induction we obtain

$$\omega(n) \le \delta^n (\sup_{x \in M} f - \inf_{x \in M} f) \qquad (5.20)$$

Since the oscillation function is decreasing on $t$ (from (A) and (B) we can get this fact), we conclude from (5.20) for $\delta = e^{-a}$

$$\omega(t) \le C e^{-at}$$

(Note that here we used of this fact that for $t \in R^+ \Longrightarrow \exists n \in N, n \le t < n+1$)

Now, we define the normalized derivative

$$\varphi = \frac{\partial u}{\partial t} - \frac{1}{Vol(M)}\int_M \frac{\partial u}{\partial t}\, d\tilde{V}$$

Notice that

$$d\tilde{V} = \det\left(g_{i\bar{j}} + \frac{\partial^2 u}{\partial z^i \partial \bar{z}^j}\right)\bigwedge_{i=1}^{n}\left(\frac{\sqrt{-1}}{2}\, dz^i \wedge d\bar{z}^j\right)$$

So

$$\frac{\partial}{\partial t}(d\tilde{V}) = \frac{\partial}{\partial t}\left(\det\left(g_{i\bar{j}} + \frac{\partial^2 u}{\partial z^i \partial \bar{z}^j}\right)\right)\bigwedge_{i=1}^{n}\left(\frac{\sqrt{-1}}{2}\, dz^i \wedge d\bar{z}^j\right)$$

$$= \frac{\partial}{\partial t}\log\det\left(g_{i\bar{j}} + \frac{\partial^2 u}{\partial z^i \partial \bar{z}^j}\right)d\hat{V} = \tilde{\Delta}\left(\frac{\partial u}{\partial t}\right)d\tilde{V}$$

The evolution equation for $\varphi$ is given by

$$\frac{\partial \varphi}{\partial t} = \frac{\partial^2 u}{\partial t^2} - \frac{1}{Vol(M)}\int_M \frac{\partial^2 u}{\partial t^2}\, d\tilde{V} - \frac{1}{Vol(M)}\int_M \frac{\partial u}{\partial t}\tilde{\Delta}\left(\frac{\partial u}{\partial t}\right)d\tilde{V} = \tilde{\Delta}\left(\frac{\partial u}{\partial t}\right) - \frac{1}{Vol(M)}\int \frac{\partial u}{\partial t}\tilde{\Delta}\left(\frac{\partial u}{\partial t}\right)d\tilde{V}$$

(5.21)

(Where here we used the fact that $\frac{\partial}{\partial t}\left(\frac{\partial u}{\partial t}\right) = \tilde{\Delta}\left(\frac{\partial u}{\partial t}\right) - \frac{\partial u}{\partial t}$ )) . Consider the quantity.

$$E = \frac{1}{2}\int_M \varphi^2 d\tilde{V}$$

Then because $\frac{\partial}{\partial t}\left(d\tilde{V}\right) = \tilde{\Delta}\left(\frac{\partial u}{\partial t}\right)d\tilde{V}$ and (5.21) we get,

$$\frac{dE}{dt} = \int_M \varphi \frac{\partial \varphi}{\partial t}\, d\tilde{V} + \frac{1}{2}\int_M \varphi^2 \tilde{\Delta}\left(\frac{\partial u}{\partial t}\right)d\tilde{V}$$

$$= \int_M \left(\frac{\partial u}{\partial t} - \frac{1}{Vol(M)}\int_M \frac{\partial u}{\partial t}\, d\tilde{V}\right)\left(\tilde{\Delta}\left(\frac{\partial u}{\partial t}\right) - \frac{1}{Vol(M)}\int_M \frac{\partial u}{\partial t}\tilde{\Delta}\left(\frac{\partial u}{\partial t}\right)\right)d\tilde{V}$$

$$+ \frac{1}{2}\int_M \varphi^2 \tilde{\Delta}\left(\frac{\partial u}{\partial t}\right)d\tilde{V} = \int_M \frac{\partial u}{\partial t}\tilde{\Delta}\left(\frac{\partial u}{\partial t}\right)d\tilde{V} + \frac{1}{2}\int_M \varphi^2 \tilde{\Delta}\left(\frac{\partial u}{\partial t}\right)d\tilde{V}$$

but we know in general, (when $M$ is compact we have $\int_M \vartheta \Delta \vartheta = - \int_M \nabla \vartheta \nabla \vartheta$ ), so the last equality is equal to

$$\int_M (-1 + \varphi) \left| \hat{\nabla} \frac{\partial u}{\partial t} \right|^2 d\hat{V} \le \frac{-1}{2} \int_M \left| \hat{\nabla} \frac{\partial u}{\partial t} \right|^2 d\hat{V} = \frac{-1}{2} \int_M \left| \hat{\nabla} \varphi \right|^2 d\hat{V}$$

Where $\left| \hat{\nabla}( \ ) \right|^2 = \hat{g}^{ij}( \ )_i ( \ )_j$ and the last inequality follows from the fact that $\inf_{x \in M} \varphi(x, t) < 0$. So

$$\sup_{x \in M} \varphi(x, t) < \sup \varphi(x, t) - \inf \varphi(x, t) = \omega(t)$$

But for $t$ large enough because $\omega(t) < Ce^{-at}$, so we can take $\omega(t) < \frac{1}{2}$, so

$$\sup_{x \in M} \varphi(x, t) < \omega(t) < \frac{1}{2}$$

Poincare inequality applied to $\varphi$ (because $\int_M \varphi = 0$)tells us that

$$\int_M \left| \hat{\nabla} \varphi \right|^2 d\hat{V} \ge \lambda_1(t) \int_M \varphi^2 d\hat{V}$$

(In reality if $f \in W_0^{1,2}(A)$ then $\|f\|_{L^2(A)} \le C \|\nabla f\|_{L^2(A)}$. The constant $C^{-2}$ is the infimum of the quotient $\min \frac{\int_A |\nabla f|^2 dx}{\int_A |f|^2}$ is the first eigenvalue of Laplacian of  ) where $\lambda_1(t)$ is the first eigenvalue of the operator $\hat{\Delta}$ at time $t$, As we proved it before , because the metrics $\tilde{g}_{i\bar{j}}(t)$ are uniformly equivalent to  $g_{i\bar{j}}(t)$ so there exists a constant $C_5 > 0$, such that $\lambda_1(t) > C_5$ for all $t$ and so this implies

$$\frac{dE}{dt} \le -C_5 E,$$

This implies that there exists a positive constant $C_6$ depending only on function $f$ such that

$$E \le C_6 e^{-C_5 t}$$

Since $\tilde{g}_{i\bar{j}}(t)$ and $g_{i\bar{j}}(t)$ are uniformly equivalent then $\sqrt{\det \tilde{g}_{i\bar{j}}(t)}$ and $\sqrt{\det g_{i\bar{j}}(t)}$ are uniformly equivalent so the volume forms $d\tilde{V}$ are uniformly equivalent to $dV$ and we also have

$$\int_M \varphi^2 dV \le \acute{C}_6 \, e^{-C_5 t} \quad (5.22)$$

**Theorem** : As $t \to \infty$, $v(x,t)$ converges to the function $v_\infty(x)$ in $C^\infty$ −Topology and $\frac{\partial u}{\partial t}$ converge to a constant in $C^\infty$ −Topology as $t \to \infty$

Proof: We will show that the family $v(x,t)$ is Cauchy in $L^1$ −norm to some function.

**Claim**: $v(x,t)$ is Cauchy in the $L^1$ −norm

**Proof of Claim**:  For any $0 < s < \acute{s}$ we have according to (5.22) and $(t) \leq C_4 e^{-at}$ ,

$$\int_M |v(x,\acute{s}) - v(x,s)| dV$$

$$\leq \int_M \int_s^{\acute{s}} \left| \frac{\partial v}{\partial t}(x,t) \right| dt dV$$

$$= \int_s^{\acute{s}} \int_M \left| \frac{\partial u}{\partial t} - \frac{1}{Vol(M)} \int_M \frac{\partial u}{\partial t} \right| dV dt \leq \int_s^{\acute{s}} \int_M |\varphi| dV dt + \int_s^{\acute{s}} \int_M \frac{1}{Vol(M)} \left| \int_M \frac{\partial u}{\partial t} d\hat{V} \right|$$

$$- \int_M \frac{\partial u}{\partial t} dV \right| dV \leq Vol(M)^{\frac{1}{2}} \int_s^{\infty} \left( \int_M \varphi^2 dV \right)^{\frac{1}{2}} dt + \frac{1}{Vol(M)} \int_s^{\infty} \omega(t) dt$$

$$\leq C_7 \int_s^{\infty} e^{-C_5 2t} dt + C_8 \int_s^{\infty} e^{-at} dt$$

So we proved that $v(x,t)$, when $t \to \infty$, are Cauchy in $L^1$-norm so because $L^1$ is complete, hence $v(x,t)$ converges in $L^1$ norm to some function $\acute{v}_\infty$.  Also we know that for some time sequence (according to our previous results) $t_n \to \infty$,  $v(x,t)$ converges to the smooth function $v(x)$ in $C^\infty$ −topology as $n \to \infty$. These together imply that  $v_\infty(x) \equiv \acute{v}_\infty(x)$ . Therefore $v(x,t)$ converges in $L^1$ norm to $v_\infty(x)$.

Claim : $v(x,t)$ converges to $v_\infty(x)$ in $C^\infty$ −topology.

Proof : suppose this is not true, so for some integer $r$, and $\varepsilon > 0$, there is a sequence $t_n$ with

$$\| v(x,t_n) - v_\infty(x) \|_{C^r} > \varepsilon.$$

But $v(x,t_n)$ are bounded in $C^\infty$ −topology so by applying Arzelà–Ascoli theorem, there is a subsequence which we again it by $v(x,t_n)$ such that $v(x,t_n)$ converge in $C^\infty$ −topology to a smooth function $\hat{v}_\infty(x) \neq v_\infty(x)$ . This is a contradiction because $v(x,t_n)$ do converge to $v_\infty(x)$ in $L^1$ norm. Hence as $t \to \infty$, $v(x,t)$ converge to $v_\infty(x)$ in $C^\infty$ −topology. This proves the first statement and consequently it follows from equation 5.17, that $\frac{du}{dt}$ converge to $\frac{du}{dt}(x,\infty)$ which is equal to $\log \det \left( g_{i\bar{j}} + \frac{\partial v_\infty}{\partial z^i \partial \bar{z}^j} \right) - \log \det \left( g_{i\bar{j}} \right) + f$ in $C^\infty$ −topology as $t \to \infty$. But also from $\omega(t) \leq C_4 e^{-at}$ we know that $\frac{du}{dt}(x,\infty)$ must be a constant function on . So we obtain the desired result and the proof of this theorem is complete. ∎

### 🔸 Convergence to a Kähler-Einstein Metric

At this point, we proved all the ingredients to complete the proof of the main theorem of this memoire, so we start to proof the main theorem.

Based on the works in previous sections we now have the following:

**Main theorem:** Let $M$ be a compact Kähler manifold of complex dimension $n$, with the Kähler metric $g_{i\bar{j}} dz^i \wedge d\bar{z}^j$. Then for any closed $(1,1)$ −form $\frac{\sqrt{-1}}{2} T_{i\bar{j}} dz^i \wedge d\bar{z}^j$ which represents the first Chern class $C_1(M)$ of $M$ one can deform the initial metric by the heat equation

$$\frac{\partial \tilde{g}_{i\bar{j}}}{\partial t} = -\tilde{R}_{i\bar{j}} + T_{i\bar{j}}$$

to another Kähler metric $\overline{g}_{i\bar{j}}$ which is in the same Kähler class as $g_{i\bar{j}}$ so that $T_{i\bar{j}}$ is the Ricci tensor of $\overline{g}_{i\bar{j}}$.

**Corollary:** If the first chern class of $M$ is equal to zero then one can deform the initial metric in the negative Ricci direction to a Ricci flat metric.

**Proof of the main theorem:** Let $R_{i\bar{j}}$ be the Ricci tensor . Then as we explained before the $(1,1)$ −form $\frac{\sqrt{-1}}{2} R_{i\bar{j}} dz^i \wedge d\bar{z}^j$ represents $C_1(M)$ . Since we assume that $\frac{\sqrt{-1}}{2} T_{i\bar{j}} dz^i \wedge d\bar{z}^j$ also represents $C_1(M)$ we know that

$$T_{i\bar{j}} - R_{i\bar{j}} = \frac{\partial^2 f}{\partial z^i \partial \bar{z}^j}$$

for some real valued smooth function $f$ on $M$. According to the first theorem of this chapter we can find a smooth function $u(x,t)$ on $M \times [0, \infty)$ which solves the problem

$$\begin{cases} \frac{\partial u}{\partial t} = \log \det \left( g_{i\bar{j}} + \frac{\partial^2 u}{\partial z^i \partial \bar{z}^j} \right) - \log \det \left( g_{i\bar{j}} \right) + f \\ \quad\quad u(x,t) = 0, \quad\quad\quad\quad t = 0 \end{cases} \quad (5.23)$$

and that

$$\tilde{g}_{i\bar{j}} = g_{i\bar{j}} + \frac{\partial^2 u}{\partial z^i \partial \bar{z}^j}$$

define a family of Kähler metrics on $M$. Moreover we know (as proved before) that as $t \to \infty$, $\tilde{g}_{i\bar{j}}$ converge in $C^\infty$ −topology to the limit metric $\tilde{g}_{i\bar{j}}(\infty)$ $(\tilde{g}_{i\bar{j}}(\infty)$ is a metric because $\tilde{g} > Cg(0)$ so it is positive definite)and that $\frac{\partial \tilde{g}_{i\bar{j}}}{\partial t}$ converges uniformly to zero .By well-known formula the Ricci tensor $\tilde{R}_{i\bar{j}}$ of $\tilde{g}_{i\bar{j}}$ is given by

$$\tilde{R}_{i\bar{j}} = -\frac{\partial^2}{\partial z^i \partial \bar{z}^j} \log \det \left( g_{i\bar{j}} + \frac{\partial^2 u}{\partial z^i \partial \bar{z}^j} \right)$$

Differentiating equation (5.23) we have

$$\frac{\partial^2}{\partial z^i \partial \bar{z}^j} \left( \frac{\partial u}{\partial t} \right) = \frac{\partial^2}{\partial z^i \partial \bar{z}^j} \log \det \left( g_{i\bar{j}} + \frac{\partial^2 u}{\partial z^i \partial \bar{z}^j} \right) - \frac{\partial^2}{\partial z^i \partial \bar{z}^j} \log \det \left( g_{i\bar{j}} \right) + \frac{\partial^2 f}{\partial z^i \partial \bar{z}^j}$$

i.e., the metric $\tilde{g}_{i\bar{j}}$ satisfy the equation

$$\frac{\partial \tilde{g}_{i\bar{j}}}{\partial t} = -\tilde{R}_{i\bar{j}} + T_{i\bar{j}}$$

Letting $t \to \infty$ we conclude that $T_{i\bar{j}} = \tilde{R}_{i\bar{j}}(\infty)$. This completes the proof of this theorem ∎ .